\documentclass[reqno]{amsart}
\usepackage{amsmath,amsfonts,amssymb,amsthm}
\usepackage[usenames,dvipsnames]{xcolor}
\usepackage{enumitem}
\usepackage{comment}
\usepackage{graphicx}
\usepackage[outdir=./]{epstopdf}
\usepackage[colorlinks=true]{hyperref} 
\usepackage{url}
\usepackage[utf8]{inputenc}
\hypersetup{linkcolor=Violet,citecolor=PineGreen}
\newtheorem{teor}{Theorem}[section]

\newtheorem{prop}[teor]{Proposition}
\newtheorem{coro}[teor]{Corollary}
\theoremstyle{definition}
\newtheorem{defi}[teor]{Definition}

\newtheorem{hipo}[teor]{Hypothesis}

\newtheorem{nota}[teor]{Remark}
\newtheorem{notas}[teor]{Remarks}
\numberwithin{equation}{section}

\newcommand{\R}{\mathbb R}
\newcommand{\T}{\mathcal T}

\newcommand{\Z}{\mathbb{Z}}
\newcommand{\Q}{\mathbb{Q}}

\newcommand{\N}{\mathbb{N}}

\newcommand{\mB}{\mathcal{B}}

\newcommand{\mU}{\mathcal{U}}

\newcommand{\mR}{\mathcal{R}}

\newcommand{\U}{\mathcal{U}}

\newcommand{\ep}{\varepsilon}

\newcommand{\mI}{\mathcal{I}}

\newcommand{\W}{\Omega}
\newcommand{\w}{\omega}
\newcommand{\lb}{\lambda}

\newcommand{\wma}{\wit{\mathfrak{a}}}
\newcommand{\wmr}{\wit{\mathfrak{r}}}

\newcommand{\ma}{\mathfrak{a}}
\newcommand{\mb}{\mathfrak{b}}
\newcommand{\mr}{\mathfrak{r}}

\newcommand{\G}{\Gamma}

\newcommand{\wit}{\widetilde}
\newcommand{\wih}{\widehat}

\newcommand{\n}[1]{\left\|#1\right\|}

\newcommand{\lsm}{\left[\begin{smallmatrix}}
\newcommand{\rsm}{\end{smallmatrix}\right]}

\begin{document}
\title[Critical transitions in quadratic ODEs]
{Critical transitions in piecewise uniformly continuous
concave quadratic ordinary differential equations}
\author[I.P. Longo]{Iacopo P. Longo}
\author[C. N\'{u}\~{n}ez]{Carmen N\'{u}\~{n}ez$^*$}
\author[R. Obaya]{Rafael Obaya}
\address[I.P. Longo]{Technische Universit\"{a}t M\"{u}nchen,
Forschungseinheit Dynamics, Zentrum Mathematik, M8,
Boltzmannstra{\ss}e 3, 85748 Garching bei M\"{u}nchen, Germany.}
\address[C. N\'{u}\~{n}ez and R. Obaya]{Departamento de Matem\'{a}tica Aplicada,
Universidad de Valladolid, Calle Doctor Mergelina s/n, 47011 Valladolid, Spain.}
\email[Iacopo Longo]{longoi@ma.tum.de}
\email[Carmen N\'{u}\~{n}ez]{carmen.nunez@.uva.es}
\email[Rafael Obaya]{rafoba@wmatem.eis.uva.es}
\thanks{All the authors were partly supported by
Ministerio de Ciencia, Innovaci\'{o}n y Universidades under project
RTI2018-096523-B-I00 and by the University of Valladolid under project PIP-TCESC-2020.
I.P.~Longo was also partly supported by the European Union's Horizon 2020
research and innovation programme under the Marie Skłodowska-Curie grant
agreement No 754462 and by TUM International Graduate
School of Science and Engineering (IGSSE)}
\date{}

\begin{abstract}
A critical transition for a system modelled by a
concave quadratic scalar ordinary differential equation
occurs when a small variation of the coefficients changes dramatically the dynamics,
from the existence of an attractor-repeller pair of hyperbolic solutions to the
lack of bounded solutions. In this paper, a tool to analyze this phenomenon for
asymptotically nonautonomous ODEs with bounded uniformly continuous or bounded
piecewise uniformly continuous
coefficients is described, and used to determine the occurrence of critical transitions
for certain parametric equations. Some numerical experiments contribute to
clarify the applicability of this tool.
\end{abstract}
\maketitle

\section{Introduction}
Substantial and irreversible changes in the output of a system upon a negligible
change in the input are referred to as {\em critical transitions\/} or {\em tipping points}.
Motivated by current exceptional challenges in nature and society (\cite{glad}, \cite{sche}),
the study of the several mechanisms leading to a critical transition has experienced
a renewed scientific thrust. In recent years, for example, it has been observed that a
time-dependent transition connecting a {\em past\/} dynamical system to a {\em future\/}
one can give rise to critical transitions
when the transition dynamics \lq\lq fails to connect
the limit ones\rq\rq~(\cite{aspw}).
This type of phenomenon has
been identified in several real scenarios including ecology (\cite{svht}, \cite{vahf}),
climate (\cite{aajqw}, \cite{aspw}, \cite{lodi}, \cite{walc}), biology (\cite{hill}),
and quantum mechanics (\cite{kato}), among others.
\par
Frequently in the literature (see for example \cite{aspw}, \cite{kijo}, \cite{rdcb}, \cite{okwi}),
the evolution of the system from
the past to the future is modeled by an asymptotically autonomous differential equation.
An asymptotically nonautonomous version of this theory has been considered
recently for the first time in \cite{lnor}, where also the past and future
systems are time-dependent: this reference deals with scalar
quadratic differential equations of the type
\begin{equation}\label{1.ecu}
 y' =-\big(y-\G(t)\big)^2+p(t)
\end{equation}
with $\G(t):=(2/\pi)\,\arctan(ct)$ for $c>0$ and
$p\colon\R\to\R$ bounded and uniformly continuous. There are two main
reasons for this choice. First, the global dynamics induced by a
quadratic differential equation is basically described by the presence
or the absence of a (classical) attractor-repeller pair of (bounded)
hyperbolic solutions. In consequence,
these equations offer a solid structure to formulate and study the possible occurrence of
critical transitions: small changes in the coefficients may cause an attractor-repeller
pair to disappear. In fact, quadratic differential equations have been identified
as prototype models for the so-called {\em rated-induced tipping}
(which we will describe below) since the very beginning (\cite{awvc}), and have been
further studied in this context (\cite{aspw}, \cite{hart}, \cite{risi1}).
Second, quadratic differential equations appear as mathematical models
in many different areas of applied sciences, which makes this formulation
interesting by itself. For instance: several model in mathematical finance respond
to this type of equations (\cite{bass}, \cite{botg});
the relation \eqref{1.ecu} is also the
Riccati equation of a two-dimensional linear hamiltonian system and the possible presence of
the attractor-repeller pair is related with the existence of an exponential dichotomy
of this linear equation, which in turn determines the existence of a local attractor or
the lack of bounded solutions in some associated nonlinear models (\cite{copp2},
\cite{jonnf});
and equations \eqref{1.ecu} are simple models of concave differential equations,
which appear often in applications and share a common dynamical description
given by the presence or absence of an attractor-repeller pair (\cite{chue}, \cite{nuos4}).
\par
In this paper, with the aim to contribute to a more robust mathematical theory of
critical transitions, we go deeper in the theoretical and numerical analysis
initiated in \cite{lnor}, which is now extended to equations \eqref{1.ecu} with much more
general coefficients, as well as to more general types of critical transitions.
When $\G$ and $p$ are arbitrary measurable
functions belonging to the Banach space $L^\infty(\R,\R)$, \eqref{1.ecu} fits in
the class of Carath\'{e}odory differential equations, which have well-known regularity
properties. We analyze the case where these coefficients are bounded
and piecewise uniformly continuous functions. A highly technical and
far from trivial extension of the methods used in \cite{lnor}
allows us to show that the description of the dynamical possibilities
there given remains valid in this extended framework. In particular,
the bifurcation analysis for
$y'=-(y-\G(t))^2+p(t)+\lb$ associates a certain real value
$\lb^*(\G,p)$ to \eqref{1.ecu}, in such a way that \eqref{1.ecu}
admits an attractor-repeller pair of hyperbolic solutions if
$\lb^*(\G,p)<0$ ({\sc case A}), it admits bounded but no hyperbolic solutions if $\lb^*(\G,p)=0$
({\sc case B}), and bounded solutions do not exist if $\lb^*(\G,p)>0$ ({\sc case C}).
\par
This description shows that the map $(\G,p)\mapsto\lb^*(\G,p)$ is a strong tool to
analyze the occurrence of critical transitions as $\G$ and $p$ vary.
However, and despite its locally Lipschitz character for the $L^\infty$-norm,
$\lb^*$ is not a continuous function on $\G$ and $p$
when the $L^1_{\rm loc}$-topology is taken on the
set of considered coefficients. Looking for this continuity is one of the
most challenging problem in this paper: it
forces us to be more restrictive in the choices
of $\G$ and $p$. More precisely, we extend the results of \cite{lnor},
to the case that: the coefficients $\G$ and $p$ are bounded and piecewise uniformly
continuous functions (BPUC, for short), with an at most countable
set of discontinuity points;
the asymptotic limits $\G(\pm\infty)$ exist
and are finite; and the equation $y'=-y^2+p(t)$ has an
attractor-repeller pair, which implies this same dynamical structure for the
past and future systems $y'=-(y-\G(\pm\infty))^2+p(t)$.
These three hypotheses
will be in force in the next paragraphs. They ensure the existence of: a
local pullback attractor for \eqref{1.ecu} which \lq\lq connects with the
attractor for the past\rq\rq~as time decreases, meaning that the distance
between both maps goes to 0 as $t\to-\infty$; and of a local pullback
repeller for \eqref{1.ecu} which \lq\lq connects with the repeller
for the future\rq\rq~as time increases. When this local pullback attractor
and repeller are globally defined and different, they form an attractor-repeller pair
which, in addition, connects those of the past and the future, and we are in {\sc case A}:
this is the situation usually called ({\em end-point}) {\em tracking}.
If the local pullback attractor is globally defined and coincides with the
local pullback repeller, then they provide a unique bounded solution, and
we are in {\sc case B}. And the only remaining possibility
is that none of them is globally defined, which corresponds to {\sc case C},
and is sometimes called {\em tipping}. When a small variation of $\G$ and $p$ changes the
dynamics  from {\sc case A} to {\sc case C} (from tracking to tipping),
we have a critical transition.
\par
In this paper, we analyze the occurrence of critical transitions as a
parameter $c$ varies for two different types of one-parametric equations of the
shape \eqref{1.ecu}, which now we write as $y'=-(y-\G_c(t))^2+p(t)$.
For both models, the function $\lb^*(\G_c,p)$ varies continuously with the parameter $c$,
and the most basic type of critical transition (which we call {\em transversal\/})
occurs when its graph crosses the vertical axis: this
means a change from {\sc case A} to {\sc case C} at a particular {\em tipping value\/} $c_0$ of the parameter.
In particular, as expected, the dynamics fits in {\sc case B} for $c=c_0$.
The previous description of these cases shows the link between this type of tipping points and
a simple nonautonomous saddle-node bifurcation pattern \cite{nuob6}: a transversal
critical transition occurs when the attractor-repeller pair collides in just one bounded solution.
Such a collision has been explored analytically and numerically in several contexts: in
one-dimensional systems (\cite{aspw}, \cite{kulo});
in higher-dimensional systems (\cite{alas}, \cite{wixt}, \cite{xie},
\cite{risi2}); in set-valued dynamical systems \cite{cari}; in random dynamical systems
(\cite{hart}); in regards to early-warning signals (\cite{risi1, risi2});
and in the nonautonomous formulation (\cite{lnor}).
There are other
points of connection between the two considered cases. For instance, a large enough
transition
$\G_c(+\infty)-\G_c(-\infty)$ guarantees the occurrence of critical transitions,
while a decreasing function $\G_c$ makes this occurrence impossible.
The role played by the size of the coefficients of the model
in the occurrence of tipping points is a key question, which appears implicit in
several works, as \cite{rdcb}, \cite{okwi} and \cite{altw}.
\par
For our first model, $\G_c(t):=c\,\G(t)$ for a $C^1$ function $\G$
(always with finite asymptotic limits), and $p$ is a BPUC function.
An in-depth analysis of the
map $c\mapsto\wih\lb(c):=\lb^*(c\,\G,p)$ shows its continuity as well as some
fundamental monotonicity properties. This allows us to prove that, if $\G$
has a local increasing point, then $\wih\lb(c)>0$ if $c$ is large enough. Since,
by hypothesis, $\wih\lb(0)<0$, at least a critical transition occurs.
In addition, there is a unique zero of $\wih\lb$ (a unique critical transition)
if $\G$ is nondecreasing.
\par
Our second model fits in a rate-induced tipping pattern, as in almost
all the afore-mentioned references.
In this case, we take $\G_c(t):=\G(c\,t)$ for a fixed $\G$,
so that $c$ determines the speed of the transition from the past system to the future system,
which are common for all $c>0$. As before, $p$ is assumed to be BPUC;
and now we include the analysis of bounded piecewise constant transition
functions $\G$.
These models seem to be physically reasonable. When the rate $c$ tends
to infinite, the transition function tends to a new piecewise constant function,
and hence the limit equation is included in the theoretical formulation. The function
$\lb_*(c):=\lb^*(\G_c,p)$ varies continuously with respect to $c$ on $\R^+\cup\{\infty\}$.
From this continuity, it is posible to deduce the tracking when the rate $c$
is small and also the occurrence of tracking or tipping when it is large enough,
based on the analysis equation corresponding to $c=\infty$. In addition, if
the piecewise constant function $\G^h$ is defined by coinciding with an initially fixed
continuous $\G$ at the discrete set $\{jh\,|\;j\in\Z\}$, and $\G_c^h(t):=\G^h(c\,t)$, then
the function $\lb_*(c,h):=\lb^*(\G_c^h,p)$ varies continuously with $c\in[0,\infty)\cup\{\R\}$
and $h\in[0,\infty)$, and hence the properties of
the continuous case can be understood by taking limits as $h$ tends to $0$.
These facts, combined with a simple numerical analysis and with an easy characterization
of $\lb_*(\infty,h)$, allow us to show interesting tipping phenomena for a quite simple
example (as its possibly revertible character)
and to explain the concept of {\em partial tipping} in our setting.
The occurrence of tipping points in piecewise constant
transition functions is also analyzed in \cite{altw} and \cite{lodi}.
\par
The paper is organized as follows. Section \ref{2.sec} extends to the most
general situation considered in the paper some dynamical properties previously known for
quadratic differential equation with continuous coefficients. An important part of the
(highly technical) proofs is postponed to Appendix \ref{appendix}. Section \ref{3.sec}
starts an in-depth
study of the bifurcation function $\lb^*(\G,p)$ and includes the analysis of the first
model above mentioned. The last two sections of the paper concern the occurrence of
rate-induced tipping for the second model. Section \ref{4.sec} deals with the
case where the functions $\G$ is continuous, whereas in Section \ref{5.sec}
the transition function is taken piecewise constant. The
phenomenon of partial tipping is described in Section \ref{4.sec}.
Appendix \ref{appendix2}, which completes the paper,
justifies the accuracy of the numerical examples included in the previous sections.
\section{General results for concave quadratic scalar ODEs}\label{2.sec}
Throughout the paper, $L^\infty(\R,\R)$ is the Banach space of
essentially bounded functions $q\colon\R\to\R$ endowed
with the norm $\n{q}$ given by the inferior of the set of
real numbers $k\ge 0$ such that the Lebesgue measure of
$\{\,t\in\R\:|\;|q(t)|>k\,\}$ is zero.
\par
Let us consider the nonautonomous concave quadratic scalar equation
\begin{equation}\label{2.ecucon}
 x'=-x^2+q(t)\,x+p(t)\,,
\end{equation}
where $q,p$ belong to $L^\infty(\R,\R)$.
Later on, we will have to be more restrictive in the choice of $q$ and $p$,
but we will first establish some general properties.
Throughout this section, $t\mapsto x(t,s,x_0)$ represents
the unique maximal solution of \eqref{2.ecucon}
satisfying $x(s,s,x_0)=x_0$, defined for $t\in\mI_{s,x_0}=(\alpha_{s,x_0},\beta_{s,x_0})$
with $-\infty\le\alpha_{s,x_0}<s<\beta_{s,x_0}\le\infty$.
Recall that, in this setting, a solution is an absolutely continuous function
on each compact interval of $\mI_{s,x_0}$
which satisfies \eqref{2.ecucon} at Lebesgue almost every $t\in\mI_{s,x_0}$;
and that $\mI_{s,x_0}=\R$ if $x(t,s,x_0)$ is bounded.
The results establishing the existence and properties of this
unique maximal solution can be found in~\cite[Chapter 2]{cole}.
Recall also that the real map $x$,
defined on an open subset of $\R\times\R\times\R$ containing
$\{(s,s,x_0)\,|\;s,x_0\in\R\}$, satisfies $x(s,s,x_0)=x_0$ and
$x(t,l,x(l,s,x_0))=x(t,s,x_0)$ whenever all the involved terms
are defined.
In fact, these results hold for Carath\'{e}odory differential
equations of more general type. For instance, those of the form \eqref{2.ecucon}
with $q,p\in L^1_{\rm loc}(\R,\R)$,
where $L^1_{\rm loc}(\R,\R)$ is the space of Borel functions $b\colon\R\to\R$
which are integrable on compact intervals (which, as explained in
Appendix \ref{appendix}, is a complete metric space).
\subsection{Hyperbolic solutions and their persistence}
Let $q,p$ belong to $L^\infty(\R,\R)$.
A bounded solution $\wit b\colon\R\to\R$ of~\eqref{2.ecucon} is said to be
{\em hyperbolic\/} if the corresponding variational equation $z'=(-2\,\wit b(t)+q(t))\,z$
has an exponential dichotomy on $\R$. That is (see~\cite{copp1}), if
there exist $k_b\ge 1$ and $\beta_b>0$ such that either
\begin{equation}\label{2.masi}
 \exp\int_s^t(-2\,\wit b(l)+q(l))\,dl\le k_b\,e^{-\beta_b(t-s)} \quad
 \text{whenever $t\ge s$}
\end{equation}
or
\begin{equation}\label{2.menosi}
 \exp\int_s^t (-2\,\wit b(l)+q(l))\,dl\le k_b\,e^{\beta_b(t-s)} \quad
 \text{whenever $t\le s$}
\end{equation}
holds. If~\eqref{2.masi} holds, the hyperbolic solution $\wit b$ is
{\em (locally) attractive}, and if \eqref{2.menosi} holds, $\wit b$ is {\em (locally) repulsive}.
In both cases, we call $(k_b,\beta_b)$ a (non-unique)
{\em dichotomy constant pair\/} for the solution $\wit b$ (or for the
equation $z'=(-2\,\wit b(t)+q(t))\,z$).
\begin{prop}\label{2.proppersiste}
Assume that~\eqref{2.ecucon} has an attractive (resp.~repulsive)
hyperbolic solution $\wit b_{q,p}$. Then, this hyperbolic
solution is persistent in the following sense: for any $\ep>0$ there exists
$\delta_\ep>0$ such that, if $\bar q,\bar p\in L^\infty(\R,\R)$
satisfy $\n{\bar q-q}<\delta_\ep$ and $\n{\bar p-p}<\delta_\ep$,
then also the perturbed differential equation
\[
 x'=-x^2+\bar q(t)\,x+\bar p(t)
\]
has an attractive (resp.~repulsive) hyperbolic solution
$\wit b_{\bar q,\bar p}$ which satisfies $\|\wit b_{q,p}-\wit b_{\bar q,\bar p}\|<\ep$.
In addition, there exists a common dichotomy constant pair for the variational
equations $z'=(-2\,\wit b_{\bar q,\bar p}(t)+\bar q(t))\,z$ corresponding to
all the functions $\bar q$ and $\bar p$ which satisfy
$\n{\bar q-q}<\delta_\ep$ and $\n{\bar p-p}<\delta_\ep$.
\end{prop}
\begin{proof}
The proof follows step by step that of \cite[Proposition 3.2]{lnor}.
Note that given $s\in  L^\infty(\R,\R)$, the equation
\begin{equation}\label{2.ecus}
 y'=(-2\,\wit b_{q,p}(t)+q(t))\,y+s(t)
\end{equation}
has a (unique) bounded solution, given by
$t\mapsto \int_{-\infty}^t u(t)\,u^{-1}(l)\,s(l)\,dl$
for $u(t):=\exp\int_{0}^t(-2\wit b(l)+q(l))\,dl$.
This allows us to define the operator $T$ on the Banach
space of real bounded continuous functions on $\R$ as in~\cite{lnor},
and repeat the whole argument used there.
\end{proof}
The next result shows the persistence also of those solutions for which the
variational equation has exponential dichotomy not in the whole of $\R$,
but in a half-line. We represent by
\lq\lq$\sup\text{ess}_{t\in\mI}$\rq\rq~the restriction of
the $L^\infty$-norm to an interval $\mI$, and by
$L^\infty(\mI,\R)$ the corresponding Banach space.
\begin{prop}\label{2.proppers}
Let $\widehat q,\widehat p\colon(-\infty,t_*]\to\R$ belong to $L^\infty((-\infty,t_*],\R)$,
where $t_*\in\R$. Assume that the equation
\begin{equation}\label{2.ecuquad}
 x'=-x^2+\widehat q(t)\,x+\widehat p(t)
\end{equation}
has a bounded solution $\wit b_{\widehat p,\widehat q}\colon(-\infty,t_*]\to\R$ satisfying
\[
 \exp\int_s^t(-2\,\wit b_{\widehat q,\widehat p}(l)+\widehat q(l))\,dl\le k\,e^{-\beta(t-s)} \quad
 \text{whenever $t_*\ge t\ge s$}
\]
for some constants $k\ge 1$ and $\beta>0$.
Given $\ep>0$, there exists $\delta_\ep>0$ such that,
if $\bar q,\bar p\colon(-\infty,t_*]\to\R$ belong to $L^\infty((-\infty,t_*],\R)$
and satisfy $\sup\text{\rm ess}_{t\le t_*}|\bar q(t)-\widehat q(t)|<\delta_\ep$ and
$\sup\text{\rm ess}_{t\le t_*}|\bar p(t)-\widehat p(t)|<\delta_\ep$, then the equation
\[
 x'=-x^2+\bar q(t)\,x+\bar p(t)
\]
has a solution $\wit b_{\bar q,\bar p}\colon(-\infty,t_*]\to\R$
such that $\sup_{t\le 0}|\wit b_{\widehat q,\widehat p}(t)-\wit b_{\bar q,\bar p}(t)|<\ep$ and
\begin{equation}\label{2.att2}
 \exp\int_s^t(-2\,\wit b_{\bar q,\bar p}(l)+\bar q(l))\,dl\le \bar k\,e^{-\bar\beta(t-s)}
\quad \text{whenever $t_*\ge t\ge s$}
\end{equation}
for some constants $\bar k\ge 1$ and $\bar \beta>0$.
\par
Let now $\widehat q,\widehat p\colon[t_*,\infty)\to\R$ belong to $L^\infty([t_*,\infty),\R)$,
and assume that the equation~\eqref{2.ecuquad}
has a bounded solution $\wit b_{\widehat p,\widehat q}\colon[t_*,\infty))\to\R$ satisfying
\[
 \exp\int_s^t(-2\,\wit b_{\widehat q,\widehat p}(l)+\widehat q(l))\,dl\ge k\,e^{\beta(t-s)} \quad
 \text{whenever $t_*\le t\le s$}
\]
for some constants $k\ge 1$ and $\beta>0$.
Then, the conclusions are analogous to those of the first case.
\end{prop}
\begin{proof}
The proof is almost identical to that of Proposition~\ref{2.proppersiste}.
The differences, in the first case, are that now we work just on $(-\infty,t_*]$,
and that~\eqref{2.ecus} may have solutions different from
$\wit b_{\widehat p,\widehat q}$ which are bounded in this interval.
Nevertheless, we can define the operator $T$ by the same expression,
acting now on the Banach space of the bounded continuous real functions on
$(-\infty,t_*]$; and the argument of~\cite{lnor} works.
The proof of the second case is analogous.
\end{proof}
\subsection{Concavity, and the sets of half-bounded and bounded solutions}
Let $q,p$ belong to $L^\infty(\R,\R)$.
The concavity on $x$ of the function giving rise to \eqref{2.ecucon}
ensures the concavity with respect to the state of the corresponding solutions:
\begin{prop}\label{2.propsemi}
As long as the involved terms are defined, we have
\[
\begin{split}
 x(t,s,\rho\, x_1+(1-\rho)\,x_2)&> \rho\, x(t,s,x_1)+(1-\rho)\,x(t,s,x_2)
 \quad\text{if $\rho\in(0,1)$ and $t> s$}\,,\\
 x(t,s,\rho\, x_1+(1-\rho)\,x_2)&< \rho\, x(t,s,x_1)+(1-\rho)\,x(t,s,x_2)
 \quad\text{if $\rho\in(0,1)$ and $t< s$}\,.
\end{split}
\]
\end{prop}
\begin{proof}
We rewrite the equation as $x'=f(t,x)$. Then,
since $f$ is strictly concave in its second argument,
$f(t,\rho\,x_1+(1-\rho)\,x_2)>\rho\,f(t,x_1)+(1-\rho)\,f(t,x_2)$
if $\rho\in(0,1)$. This inequality and the comparison result for
Carath\'{e}odory equations given in \cite[Theorem 2]{olop} (based on
the previous results of \cite{cafi}) prove the assertions.
\end{proof}
The concavity has also fundamental consequences on the properties of the sets
\[
\begin{split}
 \mB^-&:=\Big\{(s,x_0)\in\R^2\,\Big|\;\sup_{t\in(\alpha_{s,x_0},s]}x(t,s,x_0)
 <\infty\Big\}\,,\\
 \mB^+&:=\Big\{(s,x_0)\in\R^2\,\Big|\;\inf_{t\in[s,\beta_{s,x_0})}x(t,s,x_0)
 >-\infty\Big\}\,,\\
\end{split}
\]
which may be empty.
We fix $\ep>0$ and $m>0$ large enough to satisfy
$-m^2+|q(t)|\,m+|p(t)|\le-\ep$ for all $t\in\R$, which yields
$-x^2+p(t)x+q(t)\le-\ep$ for all $t\in\R$ and $|x|\ge m$. Then,
for all $(s,x_0)\in\R^2$,
$\liminf_{t\to(\alpha_{s,x_0})^+}x(t,s,x_0)>-m$ and
$\limsup_{t\to(\beta_{s,x_0})^-}x(t,s,x_0)<m$.
In other words, {\bf any solution remains upper bounded as time increases
and lower bounded as time decreases}. We will use this property repeatedly
in the paper without further reference.
In particular, $\alpha_{s,x_0}=-\infty$ for all $(s,x_0)\in\mB^-$ and
$\beta_{s,x_0}=\infty$ for all $(s,x_0)\in\mB^+$;
and $\mB:=\mB^-\cap\mB^+$ is the (possibly empty) set of
pairs $(s,x_0)$ giving rise to
(globally defined) bounded solutions of~\eqref{2.ecucon}.
\begin{nota}\label{2.notaBV}
Recall that, given a continuous function $f\colon[a,b]\to\R$ of bounded variation
(as is the case with any monotonic continuous function),
there exists a finite Borel measure $\mu$ such that $f(x)-f(a)=\mu([a,x))$.
The Radon-Nikodym decomposition of $\mu$ with respect to the Lebesgue
measure $l$, $\mu=\mu_{ac}+\mu_s$, provides the {\em singular part} of $f$,
$f_s(x):=\mu_s([a,x))$. In addition, $f$ is differentiable at
$l$-a.e.~$t\in[a,b]$ and $f'$ is $L^1$ with respect to $l$. Moreover,
if $f$ is nondecreasing, then $f'(t)\ge 0$ whenever it exists, and
$f(x)-f(a)=\int_a^x f'(t)\,dt+f_s(x)$, with $f_s$ nondecreasing and
with $f_s'(t)=0$ for $l$-a.e.~$t\in[a.b]$. Finally,
$f$ is absolutely continuous on $[a,b]$ if and only if $f_s\equiv 0$.
(See e.g.~\cite[Exercises 1.13 and 1.12, and Theorem 6.10]{rudi}.)
In particular, any bounded solution of a Carath\'{e}odory equation
satisfies the initial conditions of Theorem \ref{2.teoruno}(v).
\end{nota}
\begin{teor}\label{2.teoruno}
Let $\mB^\pm,\mB$ and $m$ be
the sets and constant above defined.
\begin{itemize}
\item[\rm(i)] If $\mB^-$ is nonempty, then there exist
a set $\mR^-$ coinciding with $\R$ or with
a negative open half-line and a maximal solution
$a\colon\mR^-\to(-\infty,m)$ of~\eqref{2.ecucon} such that,
if $s\in\mR^-$, then
$x(t,s,x_0)$ remains bounded as $t\to-\infty$ if and only if $x_0\le a(s)$; and if
$\sup\mR^-<\infty$, then $\lim_{t\to(\sup\mR^-)^-} a(t)=-\infty$.
\item[\rm(ii)] If $\mB^+$ is nonempty, then there exist
a set $\mR^+$ coinciding with $\R$ or with
a positive open half-line and a maximal solution
$r\colon\mR^+\to(-m,\infty)$ of~\eqref{2.ecucon} such that,
if $s\in\mR^+$, then $x(t,s,x_0)$ remains bounded as $t\to+\infty$
if and only if $x_0\ge r(s)$; and
if $\inf\mR^+>-\infty$, then $\lim_{t\to(\inf\mR^+)^+} r(t)=\infty$.
\item[\rm(iii)] Let $x$ be a solution defined on a maximal interval $(\alpha,\beta)$.
If it satisfies $\liminf_{t\to\beta^-}x(t)=-\infty$, then $\beta<\infty$; and if
$\limsup_{t\to\alpha^+}x(t)=\infty$, then $\alpha>-\infty$.
In particular, any globally defined solution is bounded.
\item[\rm(iv)] The set $\mB$ is nonempty if and only if $\mR^-=\R$ or $\mR^+=\R$, in
which case both equalities hold, $a$ and $ r$ are
globally defined and bounded solutions of~\eqref{2.ecucon}, and
$\mB=\{(s,x_0)\in\R^2\,|\; r(s)\le x_0\le a(s)\}\subset\R\times[-m,m]$.
\item[\rm(v)]
Let the function $b\colon\R\to\R$ be bounded, continuous,
of bounded variation and with nonincreasing singular part on every compact interval of
$\R$. Assume that $b'(t)\le-b^2(t)+q(t)\,b(t)+p(t)$
for almost all $t\in\R$. Then, $\mB$ is nonempty, and $r\le b\le a$.
If, in addition, there exists $ t_0\in\R$ such that  
$b'(t_0)<-b^2(t_0)+q(t_0)\,b(t_0)+p(t_0)$, then $r<a$.
And, if $b'(t)<-b^2(t)+q(t)\,b(t)+p(t)$ for almost all $t\in\R$,
then $r<b<a$.
\end{itemize}
\end{teor}
\begin{proof}
The proofs of (i)-(iv) repeat step by step those
of \cite[Theorem 3.1]{lnor}.
The unique required change is in (iii), where we substitute \lq\lq for all $t\ge s_0$"~by
\lq\lq for Lebesgue a.a.~$t\ge s_0$".
Let us prove (v). The comparison theorem for Carath\'{e}odory equations
(see \cite[Theorem 2]{olop}) yields $x(t,s,b(s))\ge b(t)$ for all
$s\in\R$ and
$t\ge s$, so that $(s,b(s))\in\mB^+$; and
$x(t,s,b(s))\le b(t)$ for all $t\le s$,
so that $(s,b(s))\in\mB^-$.
Consequently, $(s,b(s))\in\mB$ for all $s\in\R$:
$\mB$ is nonempty, and $r\le b\le a$.
If, in addition, there is $t_0\in\R$ with
$b'(t_0)<-b^2(t_0)+q(t_0)\,b(t_0)+p(t_0)=(d/dt)x(t,t_0,b(t_0))|_{t=t_0}$,
then an easy contradiction argument shows that there exists $t_1>t_0$ such
that $b(t_1)<x(t_1,t_0,b(t_0))$. Hence, $x(t,t_0,b(t_0))$ and $x(t,t_1,b(t_1))$
are different bounded solutions of~\eqref{2.ecucon}.
Hence, $(t_1,b(t_1)),(t_1,x(t_1,t_0,b(t_0))\in\mB$,
which ensures that $r<a$. Finally, under the last
assumption in (v), we can adapt the argument in \cite{olop}
to prove that $x(t,s,b(s))>b(t)$ whenever
$t>s$ and $x(t,s,b(s))<b(t)$ whenever $t<s$. Hence,
$a(t)=a(t,t-1,a(t-1))\ge x(t,t-1,b(t-1))>b(t)$ and
$r(t)\le x(t,t+1,b(t+1))<b(t)$ for any $t\in\R$,
which completes the proof.
\end{proof}
\begin{nota}\label{2.notabasta}
Note that~\eqref{2.ecucon} has a bounded solution
if and only if there exist times $t_1\le t_2$
(which can be equal) such that
the solutions $a$ and $r$ defined in Theorem~\ref{2.teoruno}
are respectively defined at least on $(-\infty,t_2]$ and $[t_1,\infty)$,
and $a(t)\ge r(t)$ for $t\in[t_1,t_2]$. The \lq\lq only if\rq\rq~follows
from Theorem~\ref{2.teoruno}(iv).
To check the \lq\lq if\rq\rq, we assume that, despite the
described situation, $a$ is unbounded. Then, it is not globally defined
and, since it is upper bounded,
its graph goes to $-\infty$ (that is, it has a vertical asymptote)
at a certain time to the right of $t_2$; but, if so,
this graph intersects that of $r$, impossible. Note also
that the inequality $a(t)>r(t)$ for $t\in[t_1,t_2]$
is equivalent to the existence of at least two bounded solutions.
\end{nota}
\subsection{Occurrence of an attractor-repeller pair}\label{2.subsecAR}
As said before, the main results in this paper require us to be more exigent
with the properties assumed on the coefficients of the quadratic equation
\eqref{2.ecucon}. Let $\Delta\subset\R$ be a {\em disperse} set, i.e.,
$\Delta=\{a_j\in\R\,|\; j\in\Z\}$ with $\inf_{j\in\Z}(a_{j+1}-a_j)>0\,.$
We denote by $BPUC_\Delta(\R,\R)$ the set of bounded real functions which are
defined and uniformly continuous on $\R\!-\!\Delta$. More precisely,
$q\colon\R\!-\!\Delta\to\R$ belongs to $BPUC_\Delta(\R,\R)$ if and only if
\begin{itemize}
\item[\bf c1] there is $c>0$ such that $|q(t)|<c$ for all $t\in\R\!-\!\Delta$;
\item[\bf c2] for all $\ep>0$, there is $\delta=\delta(\ep)>0$
such that,
if $t_1,t_2\in (a_j,a_{j+1})$ for some $j\in\Z$ and $t_2-t_1<\delta$, then $|q(t_2)-q(t_1)|<\ep$.
\end{itemize}
The \lq\lq$P$"~in the notation makes reference to the piecewise continuity of
$q$: it is clear that, if $q\in BPUC_\Delta(\R,\R)$, then the lateral limits
$q(a_j^+):=\lim_{t\to a_j^+}q(t)$ and $q(a_j^-):=\lim_{t\to a_j^-}q(t)$
exist for all $j\in\Z$, although possibly $q(a_j^+)\neq q(a_j^-)$. We will
assume that any function of $BPUC_\Delta(\R,\R)$ is
defined and right-continuous on the whole real line.
This assumption causes no difference in our results, but slightly
simplify the description of some of their proofs.
\begin{defi}\label{2.defiBPUC}
A {\em bounded\/} function $q\colon\R\to\R$ is {\em piecewise uniformly continuous\/}
({\em BPUC\/} for short)
if there exists a finite number of disperse sets $\Delta_1,\ldots,\Delta_n$ and
functions $q_i\in BPUC_{\Delta_i}(\R,\R)$ for $i=1,\ldots,n$ such that
$q=q_1+\cdots+q_n$.
\end{defi}
Note that the vector space $BPUC(\R,\R)$ of BPUC functions is a subset of
$L^\infty(\R,\R)$, and that the $L^\infty$-norm of a BPUC function
coincides with $\n{q}:=\sup_{t\in\R}|q(t)|$.
Clearly, any bounded and uniformly continuous
function is BPUC.
\par
Many of our results referring \eqref{2.ecucon}
consider BPUC coefficients $q$ and $p$. Let us explain the reason for this restriction.
Theorems \ref{2.teorhyp} and \ref{2.teorlb*}
provide fundamental insight in the dynamics of~\eqref{2.ecucon}: they
extend several properties proved in \cite{lnor} for bounded and
uniformly continuous functions $q,p$ to the BPUC case.
As in that paper, the construction of the hull $\W_r$ in $L^1_{\rm loc}(\R,\R^2)$
for $r:=(q,p)$ (i.e., the closure in $L^1_{\rm loc}(\R,\R^2)$ of the set of shifts
$r_t(s):=r(s+t)$), as well as of continuous flows on $\W_r$ and on $\W_r\times\R$,
are crucial tools: these constructions, standard for nonautonomous differential equations,
allow us to use techniques from topological dynamics. The
definitions of hull and flows, and the proofs of their properties,
are more technical in the present setting
of BPUC coefficients than in that of \cite{lnor}. The point is that
taking $r\colon\R\to\R^d$ with any number $d$ of BPUC component functions guarantees
the compactness of $\W_r$ and the continuity of the flows.
In order to avoid drawing focus away from the
objective of this work, we prefer to postpone a more detailed description
of these quite technical concepts and results, as well as their proofs, to
Appendix~\ref{appendix}.
We point out here that, if $r\colon\R\to\R^d$ is almost-periodic, the
topology used to define $\W_r$ on $L^1_{\rm loc}(\R,\R^d)$
coincides with that of the uniform convergence on $\R$:
see, e.g.,~\cite[Chapter 1]{fink}.
\par
Theorem \ref {2.teorhyp} shows that, if $q,p$ are BPUC functions,
then the solutions $a$ and $r$ associated to \eqref{2.ecucon} by
Theorem~\ref{2.teoruno} are globally defined and uniformly separated
if and only if they are hyperbolic. Its proof is given in Appendix \ref{appendix}.
\begin{defi}\label{2.defus}
Two globally defined solutions $x_1(t)$ and $x_2(t)$ of~\eqref{2.ecucon}
with $x_1\le x_2$ are {\em uniformly separated\/} if $\inf_{t\in\R}(x_2(t)-x_1(t))>0$.
\end{defi}
\begin{teor}\label{2.teorhyp}
Let $q,p\colon\R\to\R$ be BPUC functions,
assume that the equation~\eqref{2.ecucon} has bounded solutions,
and let $a$ and $r$ be the (globally defined) functions provided by
Theorem~{\rm\ref{2.teoruno}}. Then, the following assertions are equivalent:
\begin{itemize}
\item[\rm \rm(a)] The solutions $a$ and $r$ are uniformly separated.
\item[\rm \rm(b)] The solutions $a$ and $r$ are hyperbolic, with $a$ attractive and $r$ repulsive.
\item[\rm \rm(c)] The equation~\eqref{2.ecucon} has two different hyperbolic solutions.
\end{itemize}
In this case,
\begin{itemize}
\item[\rm \rm(i)] let $(k_a,\beta_a)$ and $(k_r,\beta_r)$ be dichotomy constant pairs for
the hyperbolic solutions $a$ and $r$, respectively, and let us choose
any $\bar\beta_a\in(0,\beta_a)$ and any $\bar\beta_r\in(0,\beta_r)$.
Then, given $\ep>0$, there exist $k_{a,\ep}\ge 1$ and
$k_{r,\ep}\ge 1$ (depending also on the choice of $\bar\beta_a$ and
of $\bar\beta_r$, respectively) such that
\[
\begin{split}
 &\quad\qquad |a(t)-x(t,s,x_0)|\le k_{a,\ep}\,e^{-\bar\beta_a(t-s)}|a(s)-x_0|
 \quad\text{if $x_0\ge r(s)+\ep$ and $t\ge s$}\,,\\
 &\quad\qquad |r(t)-x(t,s,x_0)|\le k_{r,\ep}\,e^{\bar\beta_r(t-s)}|r(s)-x_0|
 \quad\text{if $x_0\le a(s)-\ep$ and $t\le s$}\,.
\end{split}
\]
In addition,
\[
\begin{split}
 &\quad\qquad |a(t)-x(t,s,x_0)|\le k_a\,e^{-\beta_a(t-s)}|a(s)-x_0|
 \quad\text{if $x_0\ge a(s)$ and $t\ge s$}\,,\\
 &\quad\qquad |r(t)-x(t,s,x_0)|\le k_r\,e^{\beta_r(t-s)}|r(s)-x_0|
 \quad\text{if $x_0\le r(s)$ and $t\le s$}\,.
\end{split}
\]
\item[\rm \rm(ii)] The equation~\eqref{2.ecucon} does not have more
hyperbolic solutions, and
$a$ and $r$ are the only bounded solutions of
\eqref{2.ecucon} which are uniformly separated.
\end{itemize}
\end{teor}
\begin{defi}\label{2.defiAR}
In the situation described by Theorem~\ref{2.teorhyp}, $(a,r)$ is a
({\em classical\/}) {\em attractor-repeller pair\/} for~\eqref{2.ecucon}.
\end{defi}
Note that the global dynamics in the case of existence of an attractor-repeller
pair is described by Theorems~\ref{2.teoruno} and~\ref{2.teorhyp}.
\par
We include in this subsection the definitions of local pullback attractors and
repellers, which are to some extent related to the classical ones, and which play a
fundamental role in the dynamical description of the next sections: see
e.g.~Remark~\ref{3.notaABC}.
These definitions adapt those given in Section 3.8 of~\cite{klra} to the case of a
(possibly) locally defined solution.
A solution $\bar a\colon(-\infty,\beta)\to\R$ (with $\beta\le\infty$) of~\eqref{2.ecucon} is
{\em locally pullback attractive\/} if there exist $s_0<\beta$ and $\delta>0$ such that,
if $s\le s_0$ and $|x_0-\bar a(s)|<\delta$, then $x(t,s,x_0)$ is defined on
$[s,s_0]$ and, in addition,
\[
 \lim_{s\to-\infty}\max_{x_0\in [\bar a(s)-\delta,\bar a(s)+\delta]}|\bar a(t)-x(t,s,x_0)|=0
 \quad \text{for all $t\le s_0$}\,.
\]
Note that, in our scalar case, this is equivalent to say that, if
$s\le s_0$, then the solutions $x(t,s,a(s)\pm\delta)$ are defined on
$[s,s_0]$ and, in addition,
\[
 \lim_{s\to-\infty}|\bar a(t)-x(t,s,\bar a(s)\pm\delta)|=0
 \quad \text{for all $t\le s_0$}\,.
\]
A solution $\bar r\colon(\alpha,\infty)\to\R$ (with $\alpha\ge-\infty$) of~\eqref{2.ecucon} is
{\em locally pullback repulsive\/} if the solution $\bar r^*\colon(-\infty,-\alpha)\to\R$ of
$y'=-h(-t,y)$ given by $\bar r^*(t)=\bar r(-t)$ is locally pullback attractive.
In other words, it there exist $s_0>\alpha$ and $\delta>0$ such that,
if $s\ge s_0$, then the solutions $x(t,s,\bar r(s)\pm\delta)$ are defined on
$[s_0,s]$ and, in addition,
\[
 \lim_{s\to\infty}|\bar r(t)-x(t,s,\bar r(s)\pm\delta)|=0
 \quad \text{for all $t\ge s_0$}\,.
\]
\subsection{One-parametric variation of the global dynamics.}
Let us now consider the parametric family of equations
\begin{equation}\label{2.ecuconlb}
 x'=-x^2+q(t)\,x+p(t)+\lb\,,
\end{equation}
where $q$ and $p$ are BPUC functions and $\lb$ varies in $\R$.
Let $\mB_\lb$ be the (possibly empty) set of bounded
solutions, and $a_\lb$ and $r_\lb$ the corresponding bounded solutions
provided by Theorem \ref{2.teoruno} when $\mB_\lb$ is nonempty.
The next result, proved in Appendix \ref{appendix},
shows the existence of a bifurcation value
$\lb^*$: for smaller values of the parameter, there are no bounded solutions,
while for greater ones two hyperbolic solutions exist. We will
talk hence about nonautonomous saddle-node bifurcation.
\begin{teor}\label{2.teorlb*}
There exists a unique $\lb^*=\lb^*(q,p)\in[-\|q^2/4+p\|,\|p\|\,]$ such that
\begin{itemize}
\item[\rm \rm(i)] $\mB_\lb$ is empty if and only if $\lb<\lb^*$.
\item[\rm \rm(ii)] If $\lb^*\le\lb_1<\lb_2$,
then $\mB_{\lb_1}\varsubsetneq\mB_{\lb_2}$. More precisely,
\[
 r_{\lb_2}<r_{\lb_1}\le a_{\lb_1}<a_{\lb_2}\,.
\]
In addition,
$\lim_{\lb\to\infty}a_\lb(t)=\infty$ and $\lim_{\lb\to\infty}r_\lb(t)=-\infty$
uniformly on $\R$.
\item[\rm \rm(iii)] $\inf_{t\in\R}(a_{\lb^*}(t)-r_{\lb^*}(t))=0$,
and~\eqref{2.ecuconlb}$_{\lb^*}$ has no hyperbolic solution.
\item[\rm(iv)] If $\lb>\lb^*$, then $a_\lb$ and $r_\lb$ are
uniformly separated and the unique hyperbolic solutions
of~\eqref{2.ecuconlb}$_\lb$.
\item[\rm \rm(v)] $\lb^*(q,p+\lb)=\lb^*(q,p)-\lb$ for any $\lb\in\R$.
\end{itemize}
\end{teor}
\begin{teor}\label{2.teorconlb1}
Let $q,\bar q,p,\bar p\colon\R\to\R$
be BPUC functions which are norm-bounded by a constant $\kappa$,
and let $\lb^*(q,p)$ and $\lb^*(\bar q,\bar p)$
be the constants provided by Theorem~{\rm \ref{2.teorlb*}}. Then, there exists a
constant $m_\kappa$ such that
\[
 |\lb^*(\bar q,\bar p)-\lb^*(q,p)|\le m_\kappa\big(\n{\bar q-q}+\n{\bar p-p}\big)\,.
\]
In particular, the map $\lb^*\colon BPUC\times BPUC\to\R$ is
continuous.
\end{teor}
\begin{proof}
Theorem~\ref{2.teorlb*} ensures that $\lb^*(q,p)$ is bounded by
$\kappa+\kappa^2/4$. Let $m_\kappa\ge 1$ satisfy
$-m_\kappa^2+\kappa\,m_\kappa+\kappa+\kappa^2/4<0$.
Then, $\n{b}\le m_\kappa$ for any bounded solution $b$ of
$x'=-x^2+q(t)\,x+p(t)+\lb^*(q,p)$: see
Theorem~\ref{2.teoruno}. Consequently, at almost all $t\in\R$,
this bounded solution $b$ satisfies
\[
\begin{split}
 b'(t)&=-b^2(t)+\bar q(t)\,b(t)+\bar p(t)+(q(t)-\bar q(t))\,b(t)+(p(t)-\bar p(t))+\lb^*(q,p)\\
 &\le -b^2(t)+\bar q(t)\,b(t)+\bar p(t)+m_\kappa \big(\n{\bar q-q}+\n{\bar p-p}\big)+\lb^*(q,p)\,.
\end{split}
\]
Theorem~\ref{2.teoruno}(v) (see also Remark \ref{2.notaBV}) ensures that
$x'=-x^2+\bar q(t)\,x+\bar p(t)+m_\kappa \big(\n{\bar q-q}+\n{\bar p-p}\big)+\lb^*(q,p)$
has a bounded solution, and hence Theorem~\ref{2.teorlb*}(i) ensures that
$\lb^*(\bar q,\bar p)\le m_\kappa \big(\n{\bar q-q}+\n{\bar p-p}\big)+\lb^*(q,p)$.
The same argument shows that
$\lb^*(q,p)\le m_\kappa \big(\n{\bar q-q}+\n{\bar p-p}\big)+\lb^*(\bar q,\bar p)$,
and both inequalities prove the first assertion. The second one is clear.
\end{proof}
\begin{nota}\label{2.notabif}
The previous result shows that the variation in $\lb$ of the
family \eqref{2.ecuconlb} determines a nonautonomous bifurcation pattern
of saddle-node type: the absence of bounded solutions for $\lb<\lb^*(q,p)$
gives rise to the existence of an attractor-repeller pair for $\lb>\lb^*(q,p)$.
See, e.g.,~\cite{nuob6}, \cite{anja} and \cite{fuhr}.
Note also that the equation corresponding to the bifurcation value
$\lb^*(q,p)$ has either a unique bounded solution or infinitely many ones,
none of them hyperbolic. The first situation is simpler and more common,
but there are well-known examples of the second case: we refer the interested
reader to \cite{nuob6} for the details. Corollary \ref{3.corocrec}(ii)
provides a simple way to get examples of this nontrivial bifurcation pattern.
\end{nota}
\section{A particular case of concave quadratic equations}\label{3.sec}
Let us fix BPUC functions (see Definition \ref{2.defiBPUC})
$\G,p\colon\R\to\R$ such that the asymptotic limits of $\G$,
$\gamma_\pm:=\lim_{t\to\pm\infty}\G(t)$, exist and are finite. These
conditions will be in force in this initial part of Section \ref{3.sec},
whereas in some of the subsections we will impose more or less restrictive
conditions on $\G$ and $p$ which we will describe in due time.
Observe that $-2\,\G$ and $p-\G^2$ are also BPUC functions.
In what follows, we will analyze some general facts concerning
the dynamical possibilities for
\begin{equation}\label{3.ecuini}
 y'=-\big(y-\G(t)\big)^2+p(t)\,,
\end{equation}
whose solution with value $y_0$ at $t=s$ is represented by $y(t,s,y_0)$.
We understand $\G$ as a {\em transition\/} from $\gamma_-$ (in the past) to
$\gamma_+$ (in the future). In this way,
\begin{equation}\label{3.ecu+}
 y' =-(y-\gamma_+)^2+p(t)\,,
\end{equation}
and
\begin{equation}\label{3.ecu-}
 y' =-(y-\gamma_-)^2+p(t)
\end{equation}
play the role of \lq\lq limit"~equations for~\eqref{3.ecuini}
as $t\to\infty$ and as $t\to-\infty$, respectively.
We will refer to them also as {\em future equation} and {\em past
equation}. Note also that the global dynamics of these two equations
is \lq\lq identical" to that of
\begin{equation}\label{3.ecup}
 x'=-x^2+p(t)
\end{equation}
since they are obtained
from this one by the trivial changes of variables.
\begin{defi}\label{3.deficasos}
\begin{itemize}
\item[\rm -] The equation~\eqref{3.ecuini} is in \hypertarget{CA}{{\sc case A}} if
it has two different hyperbolic solutions.
\item[\rm -] The equation~\eqref{3.ecuini} is in \hypertarget{CB}{{\sc case B}} if
it has at least a bounded solution but no hyperbolic ones.
\item[\rm -] The equation~\eqref{3.ecuini} is in \hypertarget{CC}{{\sc case C}} if
it has no bounded solutions.
\end{itemize}
\end{defi}
Clearly, these three cases exhaust the possibilities.
Theorem~\ref{2.teorhyp} proves that \hyperlink{CA}{\sc case A}
is equivalent to the existence of an attractor-repeller pair, which determines
the global dynamics of \eqref{3.ecuini}. We will see below
that much more can be said in any of the three situations
if the next condition (assumed when indicated) holds:
\begin{hipo}\label{3.hipo}
The equation \eqref{3.ecup}
has an attractor-repeller pair $(\wit a,\wit r)$.
\end{hipo}
\begin{nota}\label{3.notahipo}
Hypothesis~\ref{3.hipo} is equivalent to any of these assertions:
$(\wit a+\gamma_+,\wit r+\gamma_+)$ is an
attractor-repeller pair for~\eqref{3.ecu+};
$(\wit a+\gamma_-,\wit r+\gamma_-)$ is an
attractor-repeller pair for~\eqref{3.ecu-};
$(\wit a+\G(0),\wit r+\G(0))$ is an
attractor-repeller pair for $y'=-(y-\G(0))^2+p(t)$.
\end{nota}
The next result and Remark \ref{3.notaABC} below are
fundamental to understand the dynamics of \eqref{3.ecuini}
in \hyperlink{CA}{\sc cases A}, \hyperlink{CB}{B} and \hyperlink{CC}{C}
under Hypothesis~\ref{3.hipo}.
\begin{teor}\label{3.teorexit}
Assume Hypothesis~{\rm \ref{3.hipo}}, and let
$(\wma_\pm,\wmr_\pm):=(\wit a +\gamma_\pm,\wit r+\gamma_\pm)$
be the attractor-repeller pairs for the future and past
equations \eqref{3.ecu+} and \eqref{3.ecu-}.
Then,
\begin{itemize}
\item[\rm(i)] there exist the functions $\ma$ and $\mr$
associated to \eqref{3.ecuini} by Theorem~{\rm \ref{2.teoruno}}.
\item[\rm(ii)] $\lim_{t\to-\infty}|\ma(t)-\wma_-(t)|=0$,
$\lim_{t\to-\infty}|y(t,s,y_0)-\wmr_-(t)|=0$ whenever
$\ma(s)$ exists and $y_0<\ma(s)$,
$\lim_{t\to +\infty}|\mr(t)-\wmr_+(t)|=0$, and
$\lim_{t\to +\infty}|y(t,s,y_0)-\wma_+(t)|=0$ whenever
$\mr(s)$ exists and $y_0>\mr(s)$.
\item[\rm(iii)] The solutions $\ma$ and $\mr$
are respectively locally pullback attractive and locally
pullback repulsive.
\item[\rm(iv)] If $\ma$ and $\mr$ are globally defined
and different, then they are uniformly separated, and hence
$(\wma,\wmr):=(\ma,\mr)$ is an attractor-repeller pair for
\eqref{3.ecuini}.
\item[\rm(v)] If the equation \eqref{3.ecuini}
does not have hyperbolic solutions, then
it has at most one bounded solution $\ma=\mr$.
\end{itemize}
\end{teor}
\begin{proof}
(i) Proposition~\ref{2.proppersiste} applied to
the attractor-repeller pair $(\wma_-,\wmr_-)$ of
\eqref{3.ecu-} states that, given $\ep>0$, there exists
$\delta_-=\delta_-(\ep)>0$ such that if $\n{\Sigma-\gamma^-}\le\delta_-$
then the equation $y'=-(y-\Sigma(t))^2+p(t)$ also has an attractor-repeller pair
$(\wma_\Sigma,\wmr_\Sigma)$ with $\n{\wma_\Sigma-\wma_-}\le\ep$ and
$\n{\wmr_\Sigma-\wmr_-}\le\ep$. It also states that there exists
a common dichotomy pair $(k_\ep,\beta_\ep)$ for all these functions $\Sigma$
which can be assumed to be valid for both hyperbolic solutions.
\par
We choose $t^-=t^-(\ep)<0$ such that
$|\G(t)-\gamma^-|\le\delta_-$ if $t\le t^-$, and define
$\Sigma^-(t)$ as $\G(t)$ on $(-\infty,t^-)$ and
as $\G(t^-)$ on $[t^-,\infty)$. Then,
$\big\|\Sigma^--\gamma^-\big\|\le\delta$, and hence
$y'=-(y-\Sigma^-(t))^2+p(t)$ has an attractor-repeller pair
$(\wma_{\Sigma^-},\wmr_{\Sigma^-})$, with
$\big\|\wma_{\Sigma^-}-\wma_-\big\|\le\ep$ and
$\big\|\wmr_{\Sigma^-}-\wmr_-\big\|\le\ep$.
In particular,
\begin{equation}\label{3.att3}
 \exp\int_s^t(-2\,\wma_{\Sigma^-}(l)+2\,\Sigma^-(l))\,dl\le k_\ep\,e^{-\beta_\ep(t-s)}
\quad \text{whenever $t\ge s$}\,.
\end{equation}
Let us now define $\widehat\ma_{\Sigma^-}$ as the solution of~\eqref{3.ecuini}
with value $\widehat\ma_{\Sigma^-}(t^-)=\wma_{\Sigma^-}(t^-)$. Our goal is to check
that $\widehat\ma_{\Sigma^-}$ coincides with the function $\ma$ of the statement.
Since $\widehat\ma_{\Sigma^-}(t)=\wma_{\Sigma^-}(t)$
for $t\le t^-$, it remains bounded as $t$ decreases, which proves
that $\ma$ exists and that
$\widehat\ma_{\Sigma^-}\le\ma$. To prove the converse inequality, we
take $y_0>\widehat\ma_{\Sigma^-}(t^-)$ and check that $y(t,t^-,y_0)$ is
unbounded as $t$ decreases. This property follows from
\[
 \frac{1}{k_\ep}\:e^{\beta_\ep(t^--t)}\le
 \exp\int_{t^-}^t(-2\,\widehat\ma_{\Sigma^-}(l)+2\,\Sigma^-(l))\,dl\le
 \frac{y(t,t^-,y_0)-\widehat\ma_{\Sigma^-}(t)}{y_0-\widehat\ma_{\Sigma^-}(t^-)}
\]
if $t\le t^-$: the first inequality comes from~\eqref{3.att3},
and the second one can be obtained, for instance, as (3.15) in~\cite{lnor}.
\par
To complete the proof of (i), we
work with $(\wma_+,\wmr_+)$ and use an analogous argument
in order to obtain $t^+$ such that $\mr$ is defined
al least on $[t^+,\infty)$.
\smallskip\par
(ii) We keep the notation established in the proof of (i). There,
we have checked that, given $\ep>0$, there exists $t^-$ such that
$|\ma(t)-\wma_-(t)|=|\wma_{\Sigma^-}(t)-\wma_-(t)|
\le\ep$ if $t\le t^-$, which proves the first
assertion for $\ma$ in this case. On the other hand, if
$y_0<\ma(s)$, then there exists $t_0<t^-$
such that $y(t_0,s,y_0)<\ma(t_0)=\wma_{\Sigma^-}(t_0)$.
Since $y(t,s,y_0)=y(t,t_0,y(t_0,s,y_0))$ solves
$y'=-(y-\Sigma^-(t))^2+p(t)$ for $t\le t_0$, we conclude from
Theorem~\ref{2.teorhyp}(i) that $\lim_{t\to-\infty}|y(t,s,y_0)-\wmr_-(t)|=0$.
The proofs of the two remaining assertions are similar.
\smallskip\par
(iii) Let us take $\ep\in(\,0,\,\inf_{t\in\R}(\wit a(t)-\wit r(t))\,)$.
We have obtained in (i) the time $t^-$ and the functions
$\wma_{\Sigma^-}$ and $\wmr_{\Sigma^-}$
satisfying $\inf_{s\in(-\infty,t^-]}(\ma(s)-\wmr_{\Sigma^-}(s))=
\inf_{s\in(-\infty,t^-]}(\wma_{\Sigma^-}(s)-\wmr_{\Sigma^-}(s))>\ep$.
Hence, Theorem~\ref{2.teoruno}(ii) applied to
$y'=-(y-\Sigma^-(t))^2+p(t)$ ensures that its solutions
$y^-(t,s,\ma(s)\pm\ep)$
are defined for any $t\ge s$ if $s\le t^-$.
Now we fix $t\le t^-$ and take $s\le t$. If $l\in[s,t]$, then
$\ma(l)=\wma_{\Sigma^-}(l)$  and $y(l,s,\ma(s)\pm\ep)$ coincide
with the solutions $y^-(l,s,\wma_{\Sigma^-}(s)\pm\ep)$ of
$y'=-(y-\Sigma^-(t))^2+p(t)$.
Therefore, Theorem~\ref{2.teorhyp}(i) applied to this last equation
and $\ep$ provides, for any $\beta_0\in(0,\beta_\ep)$,
a constant $k_0=k_0(\beta_0,\ep)\ge 1$ (independent of $s$) with
\begin{equation}\label{3.exprate}
 |\ma(t)-y(t,s,\ma(s)\pm\ep)|
 =|\wma_{\Sigma^-}(t)-y^-(t,s,\wma_{\Sigma^-}(s)\pm\ep)|\\
 \le k_0\,e^{-\beta_0(t-s)}\ep\,,
\end{equation}
which is as small as desired if $-s$ is large enough.
This proves (iii) in the case of $\ma$.
The proof for $\mr$ is analogous.
\smallskip\par
(iv) Assume the global existence of $\ma$ and $\mr$, with
$\mr<\ma$. According to (ii), $\lim_{t\to-\infty}|\ma(t)-\wma_-(t)|=0$
and $\lim_{t\to-\infty}|\mr(t)-\wmr_-(t)|=0$, so that their
distance is bounded from below on $(-\infty,0]$. Point (ii) also ensures
$\lim_{t\to +\infty}|\ma(t)-\wma_+(t)|=0$
and $\lim_{t\to +\infty}|\mr(t)-\wmr_+(t)|=0$. Therefore,
$\ma$ and $\mr$ are uniformly separated.
Theorem~\ref{2.teorhyp} proves that they form an attractor-repeller pair.
\smallskip\par
(v) It follows from (iv) that the unique
possibility for the existence of bounded solutions but not of hyperbolic ones
is that $\ma=\mr$, which proves (v).
\end{proof}
\begin{nota}\label{3.notaABC}
Assume the conditions on $\G$ and $p$ described at the beginning of the section,
and Hypothesis~{\rm \ref{3.hipo}}, and let
$(\wma_\pm,\wmr_\pm):=(\wit a +\gamma_\pm,\wit r+\gamma_\pm)$
be the attractor-repeller pairs for the future and past
equations \eqref{3.ecu+} and \eqref{3.ecu-}.
Under the assumed conditions on $\G$ and $p$,
Theorem \ref{3.teorexit}, combined with Theorems \ref{2.teorlb*} and \ref{2.teorhyp},
proves the next statements (among many other properties).
\par
- \hyperlink{CA}{\sc Case A} holds for \eqref{3.ecuini}
if and only if the equation has an attractor-repeller pair $(\wma,\wmr)$
(see Definition \ref{2.defiAR}); or, equivalently, if it has two different
bounded solutions. In this case, this attractor-repeller pair connects
$(\wma_-,\wmr_-)$ to $(\wma_+,\wmr_+)$:
$\lim_{t\to\pm\infty}|\wma(t)-\wma_\pm(t)|=0$ and
$\lim_{t\to\pm\infty}|\wmr(t)-\wmr_\pm(t)|=0$.
This situation is often referred to as {\em end-point tracking}.
In addition, $\wma(t)$ is the unique solution approaching $\wma_-$ as time
decreases, and $\wmr(t)$ is the unique solution approaching $\wmr_+$ as time
increases.
\par
- \hyperlink{CB}{\sc Case B} holds for \eqref{3.ecuini}
if and only if the equation has a unique bounded solution $\mb$. In this case, this
solution is locally pullback attractive and repulsive (see Subsection \ref{2.subsecAR}),
and it connects $\wma_-$ to $\wmr_+$:
$\lim_{t\to-\infty}|\mb(t)-\wma_-(t)|=0$ and
$\lim_{t\to+\infty}|\mb(t)-\wmr_+(t)|=0$.
In addition, no other solution of \eqref{3.ecuini} satisfies
any of these two properties.
\par
- \hyperlink{CB}{\sc Case C} holds
if and only if the equation has no bounded solutions. In this case, there
exists a locally pullback attractive solution $\ma$ which is the
unique solution bounded at $-\infty$ approaching $\wma_-$
as time decreases (i.e., with
$\lim_{t\to-\infty}|\ma(t)-\wma_-(t)|=0$); and it exists a locally pullback
repulsive solution $\mr$ which is the unique solution
bounded at $+\infty$ approaching $\wmr_+$
as time increases (i.e., with $\lim_{t\to+\infty}|\mr(t)-\wmr_+(t)|=0$).
This situation of loss of connection is sometimes referred to as {\em tipping}.
\par
The interested reader can in find
\cite[Figures 1-6]{lnor} some drawings showing the dynamical behavior in
each one of these three cases. (There is a last-version typo there: the
graphs of \hyperlink{CA}{\sc cases A} and \hyperlink{CC}{C} are interchanged).
\par
We also point out that, in the three dynamical cases,
the constants $\beta_0$ and $k_0$ appearing in
\eqref{3.exprate} can be chosen to get
\[
 |\ma(t)-y(t,s,y_0)|
 \le k_0\,e^{-\beta_0(t-s)}|\ma(s)-y_0|
 \quad\text{for $y_0\ge\wmr_-(s)+\ep$ and if $s\le t\le t^-$}\,.
\]
That is, $\ma(t)$
forwardly attracts exponentially fast all the solutions $y(t,s,y_0)$ starting
above $\wmr_-(s)+\ep$ for $s<t^-$ while $t\le t^-$.
Similar bounds can be found for $\mr$.
\end{nota}

\subsection{Some fundamental inequalities for $\lb^*(2\,\G,\,p-\G^2)$}
\label{3.subsecine}
Recall that Theorem \ref{2.teorlb*} associates the value $\lb^*(2\,\G,\,p-\G^2)$
to \eqref{3.ecuini}: $\lb^*(2\,\G,\,p-\G^2)$ is the bifurcation point in $\lb$ of
$x'=-(x-\G(t))^2+p(t)+\lb$. We will establish some interesting facts concerning
this value under different assumptions
on $\G$ and $p$ which will be clarified in the statement of each result.
Hypothesis \ref{3.hipo} is not in force in this subsection.
\par
Our first \lq\lq comparison" result relates $\lb^*(0,q)$ to $\lb^*(2\,\G,\,p-\G^2)$ for certain
functions $q$.
Recall that the construction of the hull $\W_p$ of a BPUC function
$p$, referred to in Section \ref{2.subsecAR}, is detailed
in Appendix \ref{appendix}. The function $p$ is {\em recurrent\/} when every
orbit of the flow on its hull is dense. It is well-known that every almost periodic
function is recurrent. In addition, the hull of any BUPC function contains
recurrent functions.
We say that a function $q\in\W_p$ {\em belongs to the alpha limit}
(resp.~{to the omega limit}) of $p$ if there exists a sequence $(t_n)_{n\ge 1}$
with limit $-\infty$ (resp.~$+\infty$) such that
$q=\lim_{n\to\infty}p_{t_n}$ on $\W_p$ (where $p_t(s):=p(s+t)$ for $t,s\in\R$).
\begin{prop}\label{3.propprec}
Let $\Gamma, p\colon\R\to\R$ be BPUC functions and let $\G$
have finite asymptotic limits.
Assume that $q\colon\R\to\R$ belongs to the alpha limit or to the omega limit of
$p$. Then, $\lb^*(0,q)\le \lb^*(2\,\G,\,p-\G^2)$.
In particular, if $p$ is recurrent, then $\lb^*(0,p)\le\lb^*(2\,\G,p-\G^2)$.
\end{prop}
\begin{proof}
We fix $q$ and $\G$ as in the statement, and denote $\lb^*:=\lb^*(2\,\G,p-\G^2)$.
Theorem~\ref{2.teorlb*} ensures the existence of
a globally bounded solution $\mb$ of
$y'=-(y-\G(t))^2+p(t)+\lb^*$.
Our goal is to check the existence of a bounded solution of
$x'=-x^2+q(t)+\lb^*$: this and Theorem~\ref{2.teorlb*}
prove that $\lb^*(0,q)\le\lb^*$.
\par
Let us work in the case of existence of $(t_n)\uparrow\infty$ such that
$q=\lim_{n\to\infty}p_{t_n}$ in $\W_p$.
Then $\mb_{t_n}(t):=\mb(t+t_n)$ solves
$y'=-(y-\G_{t_n}(t))^2+p_{t_n}(t)+\lb^*$, where $\G_{t_n}(t):=\G(t+t_n)$.
We can assume without
restriction the existence of $\lim_{n\to\infty}\mb_{t_n}(0)=:\mb_0$.
Clearly, $\lim_{n\to\infty}(-2\,\G_{t_n},p_{t_n}-\G_{t_n}^2+\lb^*)=
(-2\gamma_+,q-\gamma_+^2+\lb^*)$ in the common hull $\W_{-2\G,\,p-\G^2+\lb^*}$.
Therefore, Theorem~\ref{A.teorCont}
guarantees that the sequence of functions $(\mb_{t_n})_{n\ge 1}$
converges uniformly on compact sets
as $n\to\infty$ to the solution $\mb_{\gamma_+}$ of
$ y' =-(y-\gamma_+)^2+q(t)+\lb^*$
with $\mb_{\gamma_+}(0)=\mb_0$. In particular,
$\mb_{\gamma_+}$ is defined on the whole $\R$ and bounded.
Hence, $b:=\mb_{\gamma_+}-\gamma_+$ is a bounded
solution of $x'=-x^2+q(t)+\lb^*$, and the assertion is proved.
\par
The proof is analogous if
$q=\lim_{n\to\infty}p_{t_n}$ in $\W_p$ for $(t_n)\downarrow
-\infty$, working now with $\gamma_-$ instead of $\gamma_+$. The
last assertion is a trivial consequence of the first one.
\end{proof}
The next result compares the values of $\lb^*(2\,\G,p-\G^2)$ for two
different functions $\G$ under some conditions including
the nondecreasing character of their difference.
\begin{teor}\label{3.teorcrec}
Let $\G_1,\G_2,p\colon\R\to\R$ be BPUC functions
with $\G_2-\G_1$ nondecreasing, and let
$\lb_i:=\lb^*(2\,\G_i,\,p-(\G_i)^2)$ be the values provided by
Theorem {\rm \ref{2.teorlb*}}.
\begin{itemize}
\item[\rm(i)] If $\G_2-\G_1$ is continuous, then $\lb_1\le\lb_2$.
If, in addition, $\G_2-\G_1$ is absolutely continuous and nonconstant on
a nondegenerate interval, and $\lb_1=\lb_2$, then
$y'=-(y-\G_1(t))^2+p(t)+\lb_1$ has infinitely many bounded solutions (but no hyperbolic ones),
and the same happens for all the equations $y'=-(y-\G_\mu(t))^2+p(t)+\lb_\mu$
for $\mu\in(0,1)$, where $\G_\mu:=\mu\,\G_1+(1-\mu)\,\G_2$ and
$\lb_\mu:=\lb^*(2\,\G_\mu,\,p-(\G_\mu)^2)$.
\item[\rm(ii)] Assume that $\G_1$ and $\G_2$ have finite asymptotic limits.
Then, $\lb_1\le\lb_2$.
\end{itemize}
\end{teor}
\begin{proof}
(i) As recalled in Remark \ref{2.notaBV}, the continuous nondecreasing function
$\G_2-\G_1$ is of bounded variation, and hence there exists
$(\G_2-\G_1)'(t)\ge 0$ for Lebesgue-a.a.~$t\in\R$.
Let $\mb_2$ be a bounded solution
of $y'=-(y-\G_2(t))^2+p(t)+\lb_2$. Then, the bounded continuous function
$b_2:=\mb_2-(\G_2-\G_1)$, which is of bounded variation and has nonincreasing
singular part on every compact interval of $\R$ (see Remark \ref{2.notaBV}),
satisfies $b_2'(t)=-(b_2(t)-\G_1(t))^2+p(t)+\lb_2
-(\G_2-\G_1)'(t)\le-(b_2(t)-\G_1(t))^2+p(t)+\lb_2$
for almost all $t\in\R$. Theorem~\ref{2.teoruno}(v)
guarantees the existence of at least one bounded solution of
$x'=-(x-\G_1(t))^2+p(t)+\lb_2$. Therefore, Theorem \ref{2.teorlb*}
ensures that $\lb_1\le\lb_2$, which is the
first assertion in (i).
\par
If, in addition, $\G_2-\G_1$ is absolutely continuous and nonnonconstant
on an interval $[s,t]$, with $s<t$, it follows from $(\G_2-\G_1)(t)-(\G_2-\G_1)(s)=
\int_s^t(\G_2-\G_1)'(l)\,dl$ (see Remark \ref{2.notaBV}) that
there exists $t_0\in\R$ such that $(\G_2-\G_1)'(t_0)>0$. Therefore,
Theorem~\ref{2.teoruno}(v) ensures that
$x'=-(x-\G_1(t))^2+p(t)+\lb_2$ has more than one bounded solution.
The fact that $\lb_1=\lb_2$ implies infinitely
many bounded nonhyperbolic solutions for $x'=-(x-\G_1(t))^2+p(t)+\lb_1$
follows hence from Theorem \ref{2.teorlb*}, as explained in
Remark \ref{2.notabif}.  Finally, if we define $\G_\mu$ and $\lb_\mu$ as in
the statement, the initial assertion of (i) shows
that $\lb_1\le\lb_\mu\le\lb_2$ for any $\mu\in[0,1]$.
If  $\mu\in(0,1]$ and $\G_2-\G_1$ is nonconstant, so is $\G_2-\G_\mu$,
which is also absolutely continuous on compact intervals of $\R$.
Therefore, the argument used for $\G_1$ allows us to show the last assertion
for all these functions $\G_\mu$.
\smallskip\par
(ii) Let us fix $\ep>0$. Our goal is to prove that
$\lb_1\le \lb_2+\ep$, which ensures (ii).
Let $\kappa$ be a common bound for
$\n{\G_1}$ and $\n{\G_2}$. Theorem~\ref{2.teorconlb1} provides a constant
$\delta_\ep=\delta_\ep(\ep,\kappa)>0$ such that if $\wit \G_1$ and $\wit \G_2$ are BPUC functions
norm-bounded by $\kappa$ such that $\big\|\wit\G_1-\wit\G_2\big\|\le\delta_\ep$, then
$|\lb^*\big(2\,\wit\G_1,\,p-(\wit\G_1)^2\big)-
\lb^*\big(2\,\wit\G_2,\,p-(\wit\G_2)^2\big)|<\ep/2$. We call
$\gamma_i^\pm:=\lim_{n\to\pm\infty}\G_i(t)$, and look for a common
$t_\ep>0$ such that $|\G_i(t)-\gamma_i^\pm|\le\delta_\ep/2$ if $\pm t\ge t_\ep$
for $i=1,2$, assuming without restriction that
$\G_i(t)$ is continuous at $\pm t_\ep$ for $i=1,2$. Let
us define the BUPC functions $\G_{i,\ep}^\infty$ for $i=1,2$ by
\[
 \G_{i,\ep}^\infty(t):=\left\{\!\begin{array}{ll}
 \G_i(-t_\ep)&\quad\text{if $\;t<-t_\ep$}\,,\\
 \G_i(t)&\quad\text{if $\;-t_\ep\le t<t_\ep$}\,,\\
 \G_i(t_\ep)&\quad\text{if $\;t\ge t_\ep$}\,,\end{array}\right.
\]
so that
\begin{equation}\label{3.pri}
 |\lb^*(2\,\G_i,\,p-\G_i^2\big)-
 \lb^*(2\,\G^\infty_{i,\ep},\,p-(\G^\infty_{i,\ep})^2\big)|<\ep/2
 \quad\text{for $\;i=1,2$}\,.
\end{equation}
Now we take the smallest (finite) ordered set $\{a_0,\ldots,a_m\}$ composed by
the points of $(-t_\ep,t_\ep)$ at which either $\G_1$ or $\G_2$ are not continuous and
by $a_0:=-t_\ep$ and $a_m:=t_\ep$.
Recall that $\G_i$ is right-continuous on $a_j$ for all $j=0,\ldots,m$ and $i=1,2$.
We call $h:=\inf_{j\in\{0,\ldots,m-1\}}(a_{j+1}-a_j)>0$. For all $n\in\N$ and for
$i=1,2$, we define $\Lambda_{i,\ep}^n\colon[-t_\ep,t_\ep]\to\R$ as follows: if
$t\in [a_j,a_{j+1}-h/n)$, then
$\Lambda_{i,\ep}(t):=\G_i(t)$, whereas if $t\in[a_{j+1}-h/n,a_{j+1})$, then
$\Lambda_{i,\ep}(t):=\G_i(a_{j+1})+(a_{j+1}-t)(n/h)\big(\G_i(a_{j+1}-h/n)-\G_i(a_{j+1})\big)$.
We complete the definition to the whole line as follows:
\[
 \G_{i,\ep}^n(t):=\left\{\!\begin{array}{ll}
 \G_i(-t_\ep)&\quad\text{if $\;t<-t_\ep$}\,,\\
 \Lambda_{i,\ep}^n(t)&\quad\text{if $\;-t_\ep\le t<t_\ep$}\,,\\
 \G_i(t_\ep)&\quad\text{if $\;t\ge t_\ep$}\,,\end{array}\right.
\]
Clearly, each function $\G_{i,\ep}^n$ is continuous on $\R$ and belongs to
$BPUC(\R,\R)$, and
$\lim_{n\to\infty}\G_{i,\ep}^n(t)=\G^\infty_{i,\ep}(t)$ for all $t\in\R$.
In particular,
Lebesgue's dominated convergence theorem ensures
that the sequence $(\G_{i,\ep}^n)_{n\ge 1}$ converges to $\G_{i,\ep}^\infty$ in
$L^1_{\rm loc}(\R,\R)$; i.e., $\lim_{n\to\infty}\int_a^b|\G_{i,\ep}^n(t)-
\G_{i,\ep}^\infty(t)|\,dt=0$ whenever $a<b$.
In addition, $\G_{1,\ep}^n-\G_{2,\ep}^n$ is nondecreasing
for all $n\in\N$: it coincides with the function
$(\G_1-\G_2)_\ep^n$ which we obtain by the same procedure
starting with $\G_1-\G_2$,
and this procedure provides a nondecreasing function.
Hence, according to (i), $\lb^*(2\,\G^n_{1,\ep},\,p-(\G^n_{1,\ep})^2)
\le \lb^*(2\,\G^n_{2,\ep},\,\,p-(\G^n_{2,\ep})^2)$.
Our next purpose is showing that $\lim_{n\to\infty}\lb^*(2\,\G^n_{i,\ep},\,p-(\G^n_{i,\ep})^2)
=\lb^*(2\,\G^\infty_{i,\ep},\,p-(\G^\infty_{i,\ep})^2)$ for $i=1,2$, which yields
\[
 \lb^*(2\,\G^\infty_{1,\ep},\,p-(\G^\infty_{1,\ep})^2)
 \le \lb^*(2\,\G^\infty_{2,\ep},\,\,p-(\G^\infty_{2,\ep})^2)\,.
\]
In turn, this inequality and \eqref{3.pri}
prove $\lb_1\le\lb_2+\ep$ and complete the proof.
\par
Since the proof is the same for both values of $i$, we fix one
and omit the subindex.
Let us call $\lb_\ep(n):=\lb^*(2\,\G^n_\ep,\,p-(\G^n_\ep)^2)$ for
$n\in\N\cup\{\infty\}$, i.e., the index associated to
$y'=-\big(y-\G_\ep^n(t)\big)^2+p(t)$ by
Theorem \ref{2.teorlb*} for $n\in\N\cup\{\infty\}$. We must prove:
\begin{itemize}
\item[{\bf 1}] given $\lb<\lb_\ep(\infty)$, there exists $n_1$ such that
$\lb\le\lb_\ep(n)$ for all $n\ge n_1$,
\item[{\bf 2}] given $\lb>\lb_\ep(\infty)$, there exists $n_2$ such that
$\lb\ge\lb_\ep(n)$ for all $n\ge n_2$.
\end{itemize}
\par
Let us check {\bf 1}. Reasoning by contradiction,
we assume the existence of $\bar\lb<\lb_\ep(\infty)$ and a subsequence
$(\G_\ep^k)_{k\ge 1}$ of $(\G_\ep^n)_{n\ge 1}$ such that $\bar\lb>\lb_\ep(k)$
for all $k\ge 1$. Theorem~\ref{2.teorlb*}(i)
ensures the existence of a bounded solution $\mb_\ep^k$ of
$y'=-\big(y-\G_\ep^k(t)\big)^2+p(t)+\bar\lb$
for $k\ge 1$. The existence of a common bound for $\n{\G_\ep^k}$ for
all $k\ge 1$ ensures the existence of
$m>0$ and $\rho>0$ such that
$-m^2+2\,|\G_\ep^k(t)|\,m+|p(t)-(\G_\ep^k)^2(t)+\bar\lb|<-\rho$ for all $t\in\R$ and
$k\ge 1$. Hence, $\big\|\mb_\ep^k\big\|\le m$ for any $k\ge 1$:
see Theorem~\ref{2.teoruno}(iv). Now we take a new subsequence
$(\G_\ep^j)_{j\ge 1}$ of $(\G_\ep^k)_{k\ge 1}$ such that there exists
$y_{0}:=\lim_{j\to\infty}\mb_\ep^j(0)$.
Theorem~\ref{A.teorLC} ensures that the solution
$y^\infty_\ep(t,0,y_{0})$ of $y'=-\big(y-\G^\infty_\ep(t)\big)^2+p(t)+\bar\lb$
coincides with $\lim_{j\to\infty}\mb^j_\ep(t)$ for any $t$ in its maximal interval of
definition; therefore, it is bounded by $m$ (and hence globally defined).
This and Theorem~\ref{2.teorlb*}(i) contradict $\bar\lb<\lb_\ep(\infty)$.
Thus, {\bf 1} is proved.
\par
Let us now sketch the idea to prove {\bf 2}.
We fix $\bar\lb>\lb_\ep(\infty)$, so that the equation
\begin{equation}\label{3.ecunl}
 y'=-\big(y-\G_\ep^n(t)\big)^2+p(t)+\bar\lb
\end{equation}
corresponding to $n=\infty$ has an attractor-repeller pair
$(\wma_\ep^\infty,\wmr_\ep^\infty)$. We will
check that, if $n$ is large enough, then there exist the functions
$\ma_\ep^n$ and $\mr_\ep^n$ associated to
\eqref{3.ecunl}$_\ep^n$ by Theorem~\ref{2.teoruno}, they are respectively
defined at least on the intervals $(-\infty,t_{\ep}]$ and $[t_{\ep},\infty)$,
and they satisfy $\ma_\ep^n(t_{\ep})\ge\mr^n_\ep(t_{\ep})$.
As explained in Remark~\ref{2.notabasta},
this proves the existence of a bounded solution, and hence that
$\bar\lb\ge\lb_\ep(n)$, as {\bf 2} asserts.
\par
Observe that, outside the interval $[-t_\ep,t_\ep]$,
the coefficients of the equations~\eqref{3.ecunl}$_\ep^n$
are common for any $n\ge 1$ .
We can repeat the proof of Theorem~\ref{3.teorexit}(i), working with
the attractor-repeller pair $(\wma_\ep^\infty,\wmr_\ep^\infty)$
of \eqref{3.ecunl}$_\ep^\infty$ instead of $(\wma_-,\wmr_-)$, and
with time $-t_\ep$ instead of $t^-$. In this way we prove that,
for any $n\ge 1$, $\ma_\ep^n$ is defined at least on $(-\infty,-t_{\ep}]$,
where it coincides with $\wma_\ep^\infty$. Analogously,
for any $n\ge 1$, $\mr_\ep^n$ is defined at least on $[t_\ep,\infty)$,
where it coincides with $\wmr_\ep^\infty$.
We call $\rho:=$ $\min_{t\in[-t_\ep,t_\ep]}
\big(\wma_\ep^\infty(t)-\wmr_\ep^\infty(t)\big)>0$.
Theorem~\ref{A.teorLC} provides $n_2$ such that, for $n\ge n_2$,
\[
 \max_{t\in[-t_\ep,t_\ep]}
 \big|y_\ep^n(t,-t_\ep,\wma_\ep^\infty(-t_\ep))-
 y_\ep^\infty(t,-t_\ep,\wma_\ep^\infty(-t_\ep))\big|\le\rho\,,
\]
where $y_\ep^n(t,s,y_0)$ is the
solution of \eqref{3.ecunl}$_\ep^n$ with value $y_0$ at $t=s$.
Hence,
\[
 \min_{t\in[-t_\ep,t_\ep]}
 \big(y_\ep^n(t,-t_\ep,\wma_\ep^\infty(-t_\ep))-
 \wmr_\ep^\infty(t)\big)\ge 0\,.
\]
Altogether, we conclude that, if $t\in[-t_\ep,t_\ep]$ and $n\ge n_2$, then
\[
 \ma_\ep^n(t)=y_\ep^n(t,-t_\ep,\ma_\ep^n(-t_\ep))
 =y_\ep^n(t,-t_\ep,\wma_\ep^\infty(-t_\ep))\ge
 \wmr_\ep^\infty(t):
\]
the lower bound ensures that $\ma_\ep^n$ is also
defined on $[-t_\ep,t_\ep]$. Taking $t=t_\ep$ in the previous formula
provides the sought-for inequality and ensures {\bf 2}.
\end{proof}
\begin{coro}\label{3.corocrec}
Let $p\colon\R\to\R$ be a BPUC function.
\begin{itemize}
\item[\rm(i)] Let $\G^+,\G^-\colon\R\to\R$ be bounded, uniformly continuous,
and nondecreasing, and define $\G:=\G^+-\G^-$.
Then, $\lb^*(-2\,\G^-\!,\,p-(\G^-)^2)\le \lb^*(2\,\G,\,p-\G^2)\le
\lb^*(2\,\G^+\!,\,p-(\G^+)^2)$.
\item[\rm(ii)] Let $\G\colon\R\to\R$ be nondecreasing, and
either be a $BPUC$ function and have finite asymptotic limits or
be bounded and uniformly continuous. Then, $\lb^*(-2\,\G,\,p-\G^2)\le \lb^*(0,p)\le
\lb^*(2\,\G,\,p-\G^2)$. Moreover,
\begin{itemize}
\item[-] $\lb^*(-2\,\G,\,p-\G^2)=\lb^*(0,p)$ if $p$
is recurrent and $\G$ has finite asymptotic limits.
\item[-] Assume also that $\G$ is continuous, and
absolutely continuous and nonconstant on
a nongenenerate compact interval of $\R$. If
$\lb^*(0,p)=\lb^*(-2\,\G,\,p-\G^2)$, then
$y'=-(y+\G(t))^2+p(t)+\lb^*(-2\,\G,\,p-\G^2)$
has infinitely many bounded solutions;
and $\lb^*(0,p)<\lb^*(2\,\G,\,p-\G^2)$ if
$x'=-x^2+p(t)+\lb^*(0,p)$ has just one bounded solution.
\end{itemize}
\end{itemize}
\end{coro}
\begin{proof}
Theorem \ref{3.teorcrec}(i) ensures (i).
The first (or second) inequality in (ii) follows from Theorem
\ref{3.teorcrec} applied to $\G_1:=0$ and $\G_2:=\G$
(or $\G_1:=-\G$ and $\G_2:=0$). The assertion in (ii) concerning
a recurrent $p$ follows from Proposition \ref{3.propprec}, and the last assertions
follow also from Theorem \ref{3.teorcrec}(i).
\end{proof}
\begin{coro}
Let $p\colon\R\to\R$ be a BPUC function, and assume that
$x'=-x^2+p(t)$ does not have bounded solutions. Then,
the equation $y'=-(y-\G(t))^2+p(t)$ has no
bounded solutions in the following cases:
\begin{itemize}
\item[\rm(a)] if $p$ is recurrent and
the function $\G\colon\R\to\R$ is BPUC and has finite asymptotic
limits;
\item[\rm(b)] or if the function $\G\colon\R\to\R$ is nondecreasing and
either is $BPUC$ and has finite asymptotic
limits or is bounded and uniformly continuous.
\end{itemize}
Assume now that $x'=-x^2+p(t)$ has an attractor-repeller pair
and the conditions of \rm(b).
Then, the equation $y'=-(y+\G(t))^2+p(t)$ has an attractor-repeller pair.
\end{coro}
\begin{proof}
Assume the lack of bounded solutions.
In case (a), the result is an easy consequence of
Proposition \ref{3.propprec} and Theorem \ref{2.teorlb*}:
the lack of bounded solutions
for $x'=-x^2+p(t)$ means $\lb^*(0,p)>0$, so that
$\lb^*(2\,\G,p-\G^2)>0$, and hence $y'=-(y-\G(t))^2+p(t)$ has no
bounded solutions. The same arguments and Corollary
\ref{3.corocrec}(ii) prove case (b), as well as the last assertion.
\end{proof}
\subsection{Tipping induced by a local increment of the transition function}
\label{3.subsecsize}
The results
already proved allow us to analyze the existence of {\em tipping values of $c$}
(see Definition \ref{3.defiVIT} below) for the parametric family of equations
\begin{equation}\label{3.ecufin}
 y'=-\big(y-c\,\G(t)\big)^2+p(t)
\end{equation}
for $c\in\R$ under more restrictive conditions on $\G$ and $p$ which
we will describe in due time. We will represent by~\eqref{3.ecufin}$_c$
the equation corresponding to a fixed $c$. Observe that the corresponding future
and past equations also depend on the value of the multiplicative
parameter $c$.
\par
Our tipping analysis studies the change of the global dynamics as $c$ varies
under some assumptions involving the existence of an strictly increasing point
for $\G$. This dynamics
corresponds to \hyperlink{CA}{\sc cases A}, \hyperlink{CB}{B} or \hyperlink{CC}{C} of Definition
\ref{3.deficasos}. Recall that Theorem~\ref{2.teorhyp} shows that \hyperlink{CA}{\sc case A}
is equivalent to the existence of an attractor-repeller pair.
With the aim of talking about {\em occurrence of tipping\/}
when an attractor-repeller pair \lq\lq persists for a while
and then disappears", we define:
\begin{defi} \label{3.defiVIT}
The point $c_0\in\R$ is a {\em tipping value for the
family \eqref{3.ecufin}$_c$} if the equation \eqref{3.ecufin}$_c$ is in
\hyperlink{CA}{\sc case A} for $c$
in an open interval of endpoint $c_0$, but not at $c_0$.
\end{defi}
Theorem \ref{3.teorexit} and Remark \ref{3.notaABC} provide
more details concerning the three dynamical situations under Hypothesis \ref{3.hipo}
and the conditions assumed on $\G$ and $p$ at the beginning of Section \ref{3.sec}.
But these conditions will not in force unless otherwise indicated.
Theorem \ref{2.teorlb*} establishes a one-to-one relation between
the dynamical case of \eqref{3.ecufin}$_c$ and the sign at $c$ of the map
\begin{equation}\label{3.deflbt}
 \wih\lb\colon\R\cup\{\pm\infty\}\to\R\,,\quad c\mapsto
 \wih\lb(c):=\lb^*\big(2\,c\,\G,\,p-c^2\G^2\big)\,,
\end{equation}
given by the value associated to \eqref{3.ecufin}$_c$ by this theorem;
that is, the bifurcation point in $\lb$ of
$x'=-(x-c\,\G(t))^2+p(t)+\lb$. More precisely,
\hyperlink{CA}{\sc case A} (resp.~\hyperlink{CB}{\sc case B}, resp.~
\hyperlink{CC}{\sc case C}) occurs if and only if $\wih\lb(c)$
is strictly negative (resp.~null, resp.~strictly positive).
The next result implies that, as one might
expect, if~\eqref{3.ecufin}$_c$ undergoes a tipping at $c_0$
then \eqref{3.ecufin}$_{c_0}$ is in \hyperlink{CB}{\sc case B}.
\begin{prop}\label{3.proplabh}
Let $\G,p\colon\R\to\R$ be BPUC functions,
and let $\wih\lb$ be the map defined by \eqref{3.deflbt}. Then,
\begin{itemize}
\item[\rm (i)] for every $\kappa>0$ there exists $m_\kappa>0$ such that,
if $c_1,c_2\in[-\kappa,\kappa]$, then $|\wih\lb(c_1)-\wih\lb(c_2)|\le m_\kappa|c_1-c_2|$.
In particular, $\wih\lb$ is continuous and locally Lipschitz on $\R$.
\item[\rm(ii)] If, in addition, $\G$ is $C^1$ and $\n{\G'}:=
\sup_{t\in\R}|\G'(t)|<\infty$, then $|\wih\lb(c_1)-\wih\lb(c_2)|\le \n{\G'}\,|c_1-c_2|$
for all $c_1,c_2\in\R$. That is, under these conditions, $\wih\lb$ is Lipschitz on $\R$.
\end{itemize}
\end{prop}
\begin{proof}
Assertions (i) follow easily from Theorem \ref{2.teorconlb1}. Under the hypothesis
of (ii), for each $c\in\R$, the (bounded) change of
variable $x=y-c\,\G(t)$ takes \eqref{3.ecufin}$_c$ to
\begin{equation}\label{3.ecufintran}
 x'=-x^2+p(t)-c\,\G'(t)\,,
\end{equation}
without changing its dynamics: \hyperlink{CA}{\sc cases A}, \hyperlink{CB}{B}
or \hyperlink{CC}{C} are preserved. From this point, we check (ii) by repeating
the argument of the proof of \cite[Theorem 4.13(ii)]{lnor}.
\end{proof}
\begin{prop}\label{3.propcrec}
Let $p\colon\R\to\R$ be a BPUC function.
\begin{itemize}
\item[\rm(i)] Assume that $\G\colon\R\to\R$ is $C^1$ and that there exists a point $t_0$ at which it is
strictly increasing. Then,
there exists a value $c_0>0$ such that
\eqref{3.ecufin}$_c$ is in \hyperlink{CC}{\sc case C} for all $c\ge c_0$.
Moreover, $\lim_{c\to\infty}\wih\lb(c)=\infty$.
\item[\rm(ii)] Assume that Hypothesis {\rm \ref{3.hipo}} holds,
that $\G\colon\R\to\R$ is nonincreasing, and that
either is $BPUC$ and has finite asymptotic
limits or is bounded and uniformly continuous. Then,
\eqref{3.ecufin}$_c$ is in \hyperlink{CC}{\sc case A} for all $c\ge 0$.
\end{itemize}
\end{prop}
\begin{proof}
(i) To avoid extra technical difficulties in the proof, we assume that
$\G'(t)\ge\delta>0$ for all $t\in[0,1]$. The general case can
be proved by adapting the argument we will follow.
For each $c\in\R$, the (bounded) change of
variable $x=y-c\,\G(t)$ takes \eqref{3.ecufin}$_c$ to
\eqref{3.ecufintran}$_c$, preserving its global dynamics. We look for $c_0>0$ such that
$c_0\G'(t)\ge\pi^2+p(t)$ for all $t\in[0,1]$, and observe that the
same inequality holds for all $c\ge c_0$. Then, if $c\ge c_0$,
the solution $x_c(t,0,x_0)$ of \eqref{3.ecufintran}$_c$ with value
$x_0$ at $t=0$ satisfies $x_c(t,0,x_0)\le\pi\tan(-\pi t+\arctan(x_0/\pi))$
(which is the solution of $x'=-x^2-\pi^2$ with value $x_0$ at $t=0$)
for all the values of $t\in[0,1]$ for which they are defined.
(As usual, we take $\arctan(x_0/\pi)\in(-\pi/2,\pi/2)$.)
Since $-\pi+\arctan(x_0/\pi)<-\pi/2$, there exists $t_0\in[0,1]$
such that $\lim_{t\to t_0^-}\tan(-\pi t+\arctan(x_0/\pi))=-\infty$.
Consequently, $x_c(t,0,x_0)$ is unbounded for any $x_0\in\R$ and
$c\ge c_0$, which proves the first assertion in (i).
\par
Let us now take $k>0$. The previous property provides $c_k>0$
such that $\lb^*\big(2\,c\,\G,\,p+k-c^2\G^2\big)>0$ for all
$c\ge c_k$, and hence Theorem \ref{2.teorlb*}(v) ensures that
$\lb^*\big(2\,c\,\G,\,p-c^2\G^2\big)>k$ for all
$c\ge c_k$. This proves the last assertion in (i).
\smallskip\par
(ii) It follows from Corollary \ref{3.corocrec}(ii) and Hypothesis
\ref{3.hipo} that
$\lb^*(2\,c\,\G,p-c^2\G^2)\le \lb^*(0,p)<0$ whenever $c\ge 0$,
which proves (ii).
\end{proof}
\begin{prop}\label{3.proptip}
Let $p\colon\R\to\R$ be a BPUC function.
Assume that Hypothesis {\rm\ref{3.hipo}} holds, and that $\G\colon\R\to\R$
has finite asymptotic limits and is $C^1$, nondecreasing, and nonconstant.
Then there exists exactly a tipping value $\wih c$, which is strictly positive.
\end{prop}
\begin{proof}
Hypothesis \ref{3.hipo} ensures $\wih\lb(0)<0$, and Proposition \ref{3.propcrec}
provides at least a value of $c>0$ with $\wih\lb(c)>0$. The continuity of
$\wit\lb$ established by Proposition \ref{3.proplabh}(i) shows the existence
of a minimum $c_1>0$ with $\wih\lb(c_1)=0$.
Let us assume for
contradiction the existence of $c_2>c_1$ with $\wih\lb(c_2)=0$. By applying
Theorem \ref{3.teorcrec}(i) to $\G_1=c_1\,\G$ and $\G_2=c_2\,\G$, we deduce that
$y'=-(y-c_1\,\G(t))^2+p(t)$ has infinitely many
bounded solutions but no hyperbolic ones. But this contradicts the information provided by
Remark \ref{3.notaABC} in \hyperlink{CB}{\sc case B}.
\end{proof}
To close this section we point out that the tipping analysis just performed
can also be understood as a bifurcation analysis depending on $c$:
Proposition \ref{3.proptip}(ii) establishes conditions
under which the $c$-parametric family \eqref{3.ecufin} follows a global
saddle-node nonautonomous bifurcation pattern (see also Remark \ref{2.notabif}).
\section{Rate-induced tipping in the continuous case}\label{4.sec}
In the rest of the paper, $\G\colon\R\to\R$ represents a continuous map with finite
asymptotic limits $\gamma_\pm:=\lim_{t\to\pm\infty}\G(t)$, and
$p\colon\R\to\R$ is a BPUC function.
One of the main goals of the paper is to analyze the possibility of
occurrence of {\em rate-induced} tipping for the
one-parametric family of equations
\begin{equation}\label{4.ecuini}
 y'=-\big(y-\G_c(t)\big)^2+p(t)\,,\quad\text{with $\G_c(t):=\G(c\,t)$}
\end{equation}
for $c\in\R$ (which will be referred to as~\eqref{4.ecuini}$_c$ if $c$ is fixed).
The parameter $c$ is the {\em rate}. For $c>0$, $\G_c$ is often understood as
a {\em transition} from $\gamma_-$ to $\gamma_+$ as time increases, and $c$ determines
the velocity of this transition.
Note that the function $\G^-(t):=\G(-t)$ for $t\in\R$ maintains the same
properties required to $\G$, and $\G^-(ct)=\G(-ct)$ for every $c\in\R$.
Therefore, the analysis of~\eqref{4.ecuini}$_c$ for $c<0$ is implicitly contained in
the analysis of~\eqref{4.ecuini}$_c$ for any $\G$ and $c>0$. However,
we will formulate several properties also for $c<0$, to provide a better
understanding of the global picture.
\par
In our rate-induced tipping analysis for \eqref{4.ecuini}, a fundamental role is played
by the Carath\'{e}o\-dory equations
\begin{equation}\label{4.ecucara+}
 y'=-(y-\G_\infty(t))^2+p(t)\,,
 \quad \text{where}\;\G_\infty(t):=\left\{\begin{array}{ll}
 \gamma_-&\quad\text{if $t< 0$}\,,\\
 \gamma_+&\quad\text{if $t\ge 0$}\,,\end{array}\right.
\end{equation}
and
\begin{equation}\label{4.ecucara-}
 y'=-(y-\G_{-\infty}(t))^2+p(t)\,,
 \quad \text{where}\;\G_{-\infty}(t):=\left\{\begin{array}{ll}
 \gamma_+&\quad\text{if $t<0$}\,,\\
 \gamma_-&\quad\text{if $t\ge 0$}\,.\end{array}\right.
\end{equation}
Note that \eqref{4.ecucara+} and~\eqref{4.ecucara-} can be respectively understood
as the limiting systems of \eqref{4.ecuini}$_c$ as $c\to\infty$ and $c\to-\infty$.
(We will describe this limiting behaviour more precisely
in Section~\ref{5.sec}.) From now on, \eqref{4.ecuini}$_\infty$
and \eqref{4.ecuini}$_{-\infty}$ represent \eqref{4.ecucara+} and \eqref{4.ecucara-}.
Note also that $\G_{\pm\infty}\in BPUC_{\R-\{0\}}(\R,\R)$
(see Subsection \ref{2.subsecAR}).
\par
Following the ideas explained in Subsection \ref{3.subsecsize},
our tipping analysis studies the change of the global dynamics, determined by
\hyperlink{CA}{\sc cases A}, \hyperlink{CB}{B} or \hyperlink{CC}{C}:
\begin{defi} \label{4.defiRIT}
The point $c_0\in\R$ is a {\em tipping rate for the
family \eqref{4.ecuini}$_c$} if the equation \eqref{4.ecuini}$_c$ is in
\hyperlink{CA}{\sc case A} for $c$
in an open interval of finite endpoint $c_0$, but not at $c_0$.
A tipping rate $c_0$ is {\em transversal}
if there is an open interval containing $c_0$ such that, for values of $c$ at one
side of $c_0$, the equation~\eqref{4.ecuini}$_c$ is in
\hyperlink{CA}{\sc case A}
whereas at the other side is in \hyperlink{CC}{\sc case C}.
In the case of existence of a (transversal) tipping point $c_0$, we have
a ({\em transversal}) {\em rate-induced tipping at $c_0$.}
\end{defi}
Observe that a transversal tipping can be understood as
a local saddle-node bifurcation phenomenon occurring as the parameter $c$ varies.
In Subsection \ref{4.subseclb}, we will explain why we use the word {\em transversal\/}
in Definition \ref{4.defiRIT}. We now anticipate that it is related to the properties of the
map $c\mapsto\lb_*(c)=\lb^*(2\,\G_c,p-\G_c^2)$ for fixed $\G$ and $p$,
where $\lb^*(2\,\G_c,p-\G_c^2)$ is
the value associated to \eqref{4.ecuini}$_c$ by Theorem \ref{2.teorlb*}
for $c\in\R\cup\{\pm\infty\}$; that is, the bifurcation point in $\lb$ of
$x'=-(x-\G_c(t))^2+p(t)+\lb$. In particular, the sign of $\lb_*(c)$
determines the dynamics of \eqref{4.ecuini}$_c$. Observe that, unlike the
situation in Subsection \ref{3.subsecsize}, Theorem \ref{2.teorconlb1}
does not imply immediately the continuity of the map $\lb_*$ on the extended
real line. But we will prove this continuity in Section \ref{5.sec}.
Therefore, as in Subsection \ref{3.subsecsize},
if~\eqref{4.ecuini} undergoes a rate-induced tipping at $c_0$,
then \eqref{4.ecuini}$_{c_0}$ is in \hyperlink{CB}{\sc case B}.
\par
Observe that the future and past equations are given for
any $c>0$ by
\begin{equation}\label{4.ecu+}
 y' =-(y-\gamma_+)^2+p(t)
\end{equation}
and \begin{equation}\label{4.ecu-}
 y' =-(y-\gamma_-)^2+p(t)\,,
\end{equation}
while the roles of~\eqref{4.ecu+} and~\eqref{4.ecu-} are
interchanged for $c<0$.
We consider also the equation~\eqref{4.ecuini}$_0$, namely
\begin{equation}\label{4.ecu0}
 y' =-(y-\gamma_0)^2+p(t)
\end{equation}
for $\gamma_0:=\G(0)$. Note that the global dynamics of these three equations
is \lq\lq identical", since all of them are obtained from
\begin{equation}\label{4.ecup}
 x'=-x^2+p(t)
\end{equation}
by trivial changes of variables: see Remark \ref{3.notahipo}.
\par
For the reader's convenience, we complete this initial part of the section
by repeating the fundamental
Hypothesis \ref{3.hipo}:
\begin{hipo}\label{4.hipo}
The equation \eqref{4.ecup}
has an attractor-repeller pair $(\wit a,\wit r)$.
\end{hipo}
\begin{notas}\label{5.notasigno}
1.~Under this condition, all the information provided
by Theorem \ref{3.teorexit} and Remark \ref{3.notaABC} applies to
\eqref{4.ecuini}$_c$ for any $c\in\R\cup\{\pm\infty\}$. But one must have in
mind that the future and past equations, and hence the corresponding
attractor-repeller pairs, depend on the sign of $c$: $(\wit a +\gamma_+,\wit r+\gamma_+)$
is the future (resp.~past) pair for $c>0$ (resp.~$c<0$),
$(\wit a +\gamma_-,\wit r+\gamma_-)$
is the past (resp.~future) pair for $c>0$ (resp.~$c<0$),
and $(\wit a +\gamma_0,\wit r+\gamma_0)$
is the future and past pair for $c=0$.
\par
2.~Proposition \ref{5.proplim}, proved in the next section,
shows that part of the dynamical properties described in
Theorem \ref{3.teorexit} and Remark \ref{3.notaABC}
also hold when Hypothesis \ref{4.hipo} is substituted by the
existence of attractor-repeller pair for \eqref{4.ecucara+} or
\eqref{4.ecucara-}.
\end{notas}
\subsection{The bifurcation curve $\lb_*(c)$}\label{4.subseclb}
Let us define
\begin{equation}\label{4.deflb*}
 \lb_*\colon\R\cup\{\pm\infty\}\to\R\,,\quad c\mapsto \lb_*(c):=
 \lb^*\big(2\,\G_c,\,p-(\G_c)^2\big)\,
\end{equation}
and recall the relation between the dynamical situation of \eqref{4.ecuini}$_c$
(\hyperlink{CA}{\sc case A}, {\sc \hyperlink{CB}{B}}, or {\sc \hyperlink{CC}{C}})
and the sign of $\lb_*(c)$ (negative, null, or positive).
In particular, Hypothesis~\ref{4.hipo} can be reformulated as $\lb_*(0)<0$:
see Remark ~\ref{3.notahipo}. This hypothesis is not in force for the next
result, already mentioned, and proved by Theorem~\ref{5.teorcon} in the
next section.
\begin{teor}\label{4.teorcon}
Let $\lb_*\colon\R\cup\{\pm\infty\}\to\R$ be defined by~\eqref{4.deflb*}. Then,
\begin{itemize}
\item[\rm \rm(i)] the map $\lb_*$ is bounded: it takes values in $[-||p||,||p||+||\G||^2]$.
\item[\rm \rm(ii)] the map $\lb_*$ is continuous on the extended real line.
\end{itemize}
\end{teor}
So, according to Definition \ref{4.defiRIT}, $c_0\in\R$
is a tipping rate if $\lb_*(c_0)=0$ and there is $\delta_0>0$
such that $\lb_*(c)<0$ either for $c\in(c_0-\delta_0,c_0)$ or for
$c\in(c_0,c_0+\delta_0)$; and a transversal tipping rate
if, in addition, $\lb_*(c)>0$ either for $c\in(c_0-\delta_0,c_0)$
or for $c\in(c_0,c_0+\delta_0)$. Hence,
the graph of the continuous map $\lb_*$ crosses the horizontal axis
transversally at a transversal tipping-rate $c_0$.
Although one might expect this situation to be the most frequent one,
other types of tipping are also possible.
\par
Assuming Hypothesis \ref{4.hipo}, the continuity of $\lb_*$
established in Theorem \ref{4.teorcon} allows us to determine the dynamical
case of \eqref{4.ecuini}$_c$ for small values of the rate, and also for large
ones under additional conditions.
This is what the next theorem states. Its
scope will be clearer in Subsection \ref{4.subsec_partial_tipping}.
\begin{teor}\label{4.teorc_infty}
Assume Hypothesis~{\rm \ref{4.hipo}}, and let $(\wit a_\lb, \wit r_\lb)$
be the attractor-repeller pair for $x'=-x^2+p(t)+\lb$ for $\lb>\lb_*(0)$
(with $(\wit a_0,\wit r_0)=(\wit a,\wit r)$). Then,
\begin{itemize}[itemsep=2pt]
\item[\rm \rm(i)] there exists $c_0>0$ such that~\eqref{4.ecuini}$_c$
is in \hyperlink{CA}{\sc case A} for $c\in(-c_0,c_0)$.
\item[\rm \rm(ii)] If $\wit a(0)-\wit r(0)>\gamma_+-\gamma_-$, then
the equation \eqref{4.ecucara+} has an attractor-repeller pair, and hence there
exists a minimum  $c_M\ge0$ such that~\eqref{4.ecuini}$_c$
is in \hyperlink{CA}{\sc case A} for $c>c_M$.
\item[\rm \rm(iii)] If $\wit a(0)-\wit r(0)<\gamma_+-\gamma_-$, then
the equation \eqref{4.ecucara+} has no bounded solutions, and hence
there is a minimum $c_M^*>0$ such that~\eqref{4.ecuini}$_c$
is in \hyperlink{CC}{\sc case C}
for $c>c_M^*$. In this case, $\lb_*(\infty)=\lb_\infty$, where
$\lb_\infty>0$
is the unique value of the parameter such that
$\wit a_{\lb_\infty}(0)-\wit r_{\lb_\infty}(0)=\gamma_+-\gamma_-$.
\item[\rm \rm(iv)] If $\wit a(0)-\wit r(0)>\gamma_--\gamma_+$, then
the equation \eqref{4.ecucara-} has an attractor-repeller pair, and hence
there exists  a maximum $c_m\le0$ such that~\eqref{4.ecuini}$_c$
is in \hyperlink{CA}{\sc case A} for $c<c_m$.
\item[\rm \rm(v)] If $\wit a(0)-\wit r(0)<\gamma_--\gamma_+$, then
the equation \eqref{4.ecucara-} has no bounded solutions, and hence
there is a maximum $c_m^*<0$ such that~\eqref{4.ecuini}$_c$
is in \hyperlink{CC}{\sc case C} for $c<c_m^*$.
In this case, $\lb_*(-\infty)=\lb_{-\infty}>0$, where
$\lb_{-\infty}$ is the unique value of the parameter such that
$\wit a_{\lb_{-\infty}}(0)-\wit r_{\lb_{-\infty}}(0)=\gamma_--\gamma_+$.
\end{itemize}
\end{teor}
\begin{proof}
This result is included in Theorem \ref{5.teorc_infty}, proved in the next
section.
\end{proof}
\begin{nota}\label{4.notaGp}
It follows from Theorems \ref{4.teorc_infty} and \ref{3.teorexit} that
the equation \eqref{4.ecucara+} (or~\eqref{4.ecucara-})
has only a bounded solution if and only
if $\wit a(0)-\wit r(0)=\gamma_+-\gamma_-$
(or $\wit a(0)-\wit r(0)=\gamma_--\gamma_+$).
But $\lb_*(\infty)=0$ (or $\lb_*(-\infty)=0$)
does not lead us to any conclusion for large
values of $c$ (or $-c$).
Observe also that, since
$\wit a$ and $\wit r$ depend just on $p$ while $\gamma^+$ and $\gamma^-$ depend just on
$\G$, a suitable choice of $\G$ once $p$ is fixed (or the converse)
determines the dynamical situation of \eqref{4.ecucara+}: tipping occurs
when $\wit a(0)-\wit r(0)<\gamma_+-\gamma_-$ (which is not possible if
$\gamma_+<\gamma_-$), and there is tracking whenever $\gamma^+-\gamma^-$
is large enough. Analogous conclusions hold for \eqref{4.ecucara-}.
\end{nota}
\begin{nota}
Proposition \ref{5.propcrecch} applied to the case $h=0$
establishes some inequalities regarding $\lb_*(c)$ which
provide valuable information about the corresponding
dynamical case for \eqref{4.ecuini}$_c$ depending on that for
\eqref{4.ecuini}$_0=\;$\eqref{4.ecup}.
\end{nota}
\subsection{Partial and total tipping on the hull}\label{4.subsec_partial_tipping}
Partial tipping and total tipping are phenomena introduced in~\cite{alas} in the context of
two-dimensional asymptotically autonomous systems.
In particular, the attractors of
the future and past systems in \cite{alas} are compact sets, each given by the
trajectory of an orbitally asymptotically stable solution. The associated pullback
attractor for the nonautonomous system is hence a compact nonautonomous set which
is not a singleton. Upon the variation of the rate, it is shown that the associated pullback
attractor can break up in the sense that some of the trajectories limiting at the limit
cycle of the past limit system also limit at the limit cycle of the future system,
but others fail to do so, and {\em partial tipping} occurs. If all the trajectories
which determine the pullback attractor do not limit to the limit cycle of the future
system, then a {\em total tipping\/} happens.
\par
While the current state of the art does not
allow us to pose the same question in the context of two-dimensional asymptotically
nonautonomous systems, the results of the previous section do allow us to address
a different phenomenon which can still be regarded as an instance of partial and total
tipping.
\par
The key point is that, in the nonautonomous case, our Hypothesis \ref{4.hipo}
of existence of an attractor-repeller pair for the equation $x'=-x^2+p(t)$, and hence for the
future and past equations $y'=-(y-\gamma_\pm)^2+p(t)$, means the existence
of two hyperbolic copies of the base for the corresponding skew-product flow defined on the
hull $\W_p$ of $p$ (described in Subsection \ref{3.subsecine}). The proof of this assertion
is the fundamental point in the proof of \cite[Theorem 3.5]{lnor}, where
the interested reader can find a more detailed explanation of the meaning
of hyperbolic copy of the base. What is interesting for us, now,
is that this property ensures that each equation $x'=-x^2+q(t)$ given by $q\in\W_p$,
as well as the corresponding future and past equations, has an attractor-repeller pair
(given by the corresponding sections of the hyperbolic copies of the base).
So, a natural question arises: for a given value of $c>0$ (or $c=\infty$),
are all the equations $y'=-(y-\G_c(t))^2+q(t)$, where $\G$ is fixed and $q$
varies in the hull $\W_p$, in the same dynamical case?
We will talk about {\em partial tipping on the hull} when
\hyperlink{CA}{\sc cases A} and \hyperlink{CC}{C}
coexist for different functions in the hull for a given value of $c$, and about
{\em total tipping on the hull} when the dynamics is always in \hyperlink{CC}{\sc case C}.
The global occurrence of \hyperlink{CA}{\sc case A} is {\em total tracking on the hull}.
\par
The next example has a double purpose: to illustrate a simple way to determine
the dynamical situation of \eqref{4.ecucara+}$\,=\,$\eqref{4.ecucara+}$_\infty$,
and hence that of \eqref{4.ecuini}$_c$
for large enough $c$; and to show a situation of partial tipping on the hull.
Let us define
\begin{equation}\label{4.eqpt}
 \G(t):=\frac{2}{\pi}\arctan(t)\quad\text{and}\quad p(t):=0.962-\sin(t/2)-\sin(\sqrt{5}\,t).
\end{equation}
The choice of $\G$ and $p$ is not coincidental: it permits a direct connection
to the numerical analysis carried out in~\cite{lnor}, which features the problem given
by the same $\G$ and $p(t):=0.895-\sin(t/2)-\sin(\sqrt{5}\,t)$.
Theorem~\ref{2.teorlb*}(v) guarantees that the bifurcation curve $\lb_*$ of the
equation~\eqref{4.ecuini}, defined by
\eqref{4.deflb*} for the chosen $\G$ and $p$ in~\eqref{4.eqpt}, is a
vertical translation (of $-0.067$) of that depicted in Figure 8 of~\cite{lnor}.
As justified in \cite{lnor}, we can assume that Hypothesis \ref{4.hipo} holds.
We can also assume that, for all $c\in(0,50]\cup\{\infty\}$, we are able to
approximate beyond machine precision the (possibly locally defined) solutions
$\ma_{c}$ and $\mr_{c}$ associated to~\eqref{4.ecuini} by Theorem~\ref{2.teoruno}.
A detailed clarification supporting this
last assumption is given in Appendix \ref{appendix2}.
\begin{figure}[]
\caption{Characterization of the dynamics for the differential equation
\eqref{4.ecucaras} for $s\in[-40,40]$. In the upper panel the attractor
$\wit a$ (in red) and repeller $\wit r$ (in blue) of $x'=-x^2+p(t)$.
In the lower panel, the curves $\lb^\infty(s)$
(solid green curve in the lower panel) and $d^\infty(s)=2-\wit a(s)+\wit r(s)$
(magenta in the lower panel). The (common) points $s$ on which they are
strictly positive (i.e., the points $s$ for which $\wit a(s)$ and $\wit r(s)$
are close enough) are highlighted in thick red on the
axis $y=0$ in both panels. These are the points for which
\eqref{4.ecucaras} has no bounded solutions,
The complementary of the closure of this set, given by the points for which
$\lb^\infty(s)$ and $d^\infty(s)$ are strictly negative, is composed by the points
for which \eqref{4.ecucaras} has an attractor-repeller pair.
}
\includegraphics[trim={1cm 9.2cm 0.5cm 8.8cm},clip,width=\textwidth]{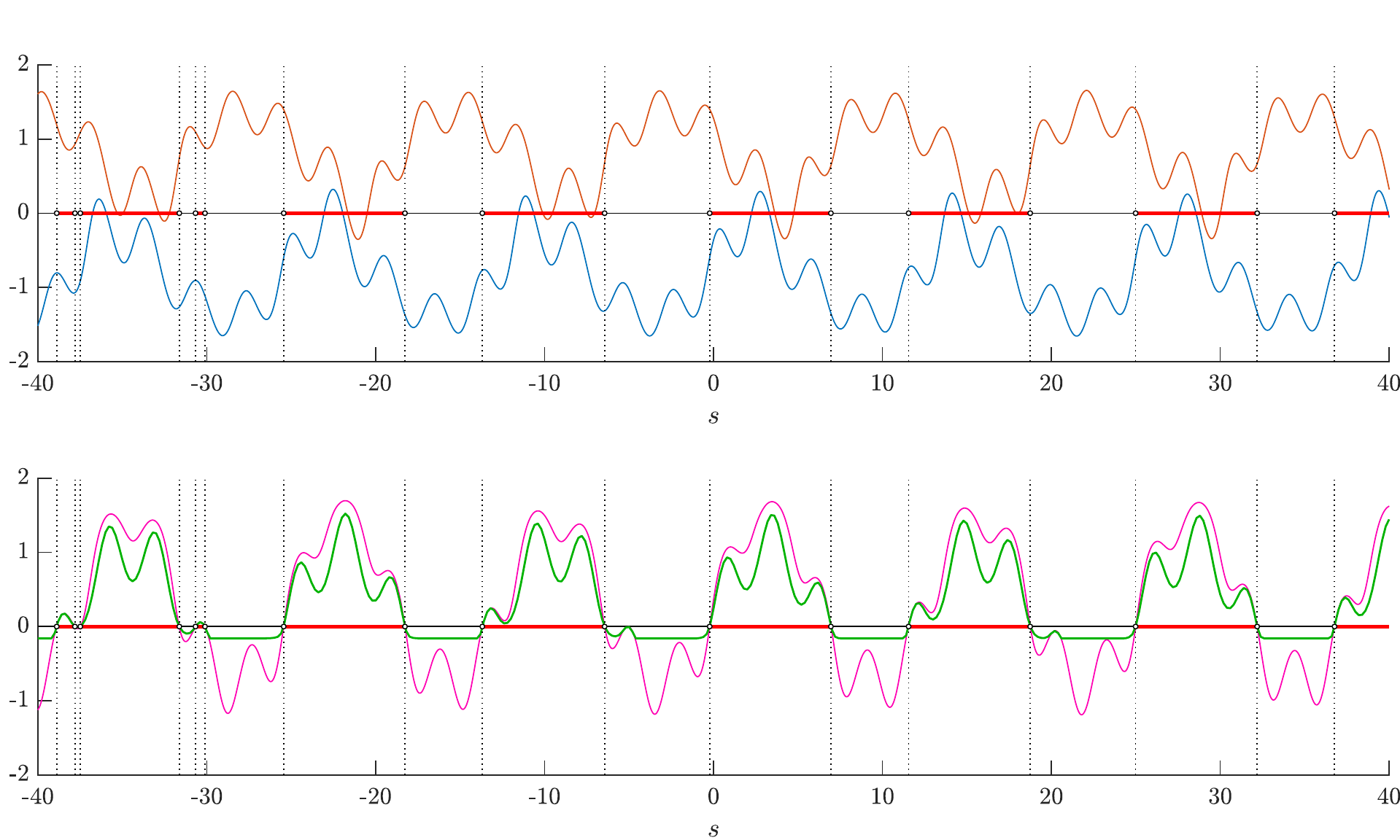}
\label{fig:theorem3_2}
\vspace{-0.5cm}
\end{figure}
\par
As said in Subsection \ref{2.subsecAR}, since $p$ is a quasiperiodic function,
the corresponding hull $\W_p$ is constructed as the closure of the set of the shifts
$p_s(t):=p(t+s)$ in the uniform topology.
For this reason (as we will explain later),
instead of working with the whole hull, it
suffices to our purposes working with the shifts of $p$.
Therefore, we consider the equations
\begin{equation}\label{4.ecuinis}
 y'=-\big(y-\G_c(t)\big)^2+p_s(t)
\end{equation}
and their limits as $c\to\infty$,
\begin{equation}\label{4.ecucaras}
 y'=-(y-\G_\infty(t))^2+p_s(t)\,,
\end{equation}
for $s\in\R$. Let us call $\wit a_s(t):=\wit a(t+s)$
and $\wit r_s(t):=\wit r(t+s)$. It is easy to check that
$(\wit a_s,\wit r_s)$ is the attractor-repeller pair of $x’=-x^2+p_s(t)$.
Theorem~\ref{4.teorc_infty}(ii)\&(iii) (see also their proofs) reveal that, for a given value
of $s$, \eqref{4.ecucaras} is in \hyperlink{CA}{\sc case A} (resp.~\hyperlink{CC}{\sc case C}),
and hence the same happens with \eqref{4.ecuini}$_c$ for large enough $c$,
if and only if
\[
 d^\infty(s):=\gamma_+-\gamma_--\wit a_s(0)+\wit r_s(0)=2-\wit a(s)+\wit r(s)
\]
is strictly negative (resp.~positive).
That is, if the distance from $\wit r(s)$ to $\wit a(s)$ is large enough,
then \eqref{4.ecucaras} and all the equations \eqref{4.ecuinis}$_c$ for large
enough $c$ (depending on $s$) have an attractor-repeller pair which, according to
Remark \ref{3.notaABC}, connects
that of the past equation $y’=-(y+1)^2+p_s(t)$ to that of the future equation
$y’=-(y-1)^2+p_s(t)$; i.e, $(\wit a_s-1,\wit r_s-1)$ to $(\wit a_s+1,\wit r_s+1)$.
And if the distance is small, then the tracking is lost and tipping occurs: there
are no longer bounded solutions. This guarantees the existence of at least a tipping
rate $c_0>0$ for the $c$-parametric family \eqref{4.ecuinis}, as
Theorem \ref{4.teorc_infty}(iii) ensures.
\par
\begin{figure}[]
\caption{Numerical simulation of the bifurcation map $\lb_*(c,s)$ of \eqref{4.ecuinis}$_c$
(surface with gradient color). The red grid identifies the plane
$\lb=0$. Consequently, the points of the surface below it correspond to {\sc {Case A}},
the points above to {\sc {Case C}}, and the points of intersection to {\sc {Case B}}. The curve
in green is the graph of the bifurcation curve
$\lb^\infty(s)$ of \eqref{4.ecucaras}
(see also Figure~\ref{fig:theorem3_2}). Theorem~\ref{5.teorcon}
guarantees the convergence
of $\lb_*(c,s)$ to $\lb^\infty(s)$ as $c$ increases. The figure
indicates how fast this convergence is.}
\includegraphics[trim={1cm 9.6cm 1cm 10.2cm},clip,width=\textwidth]{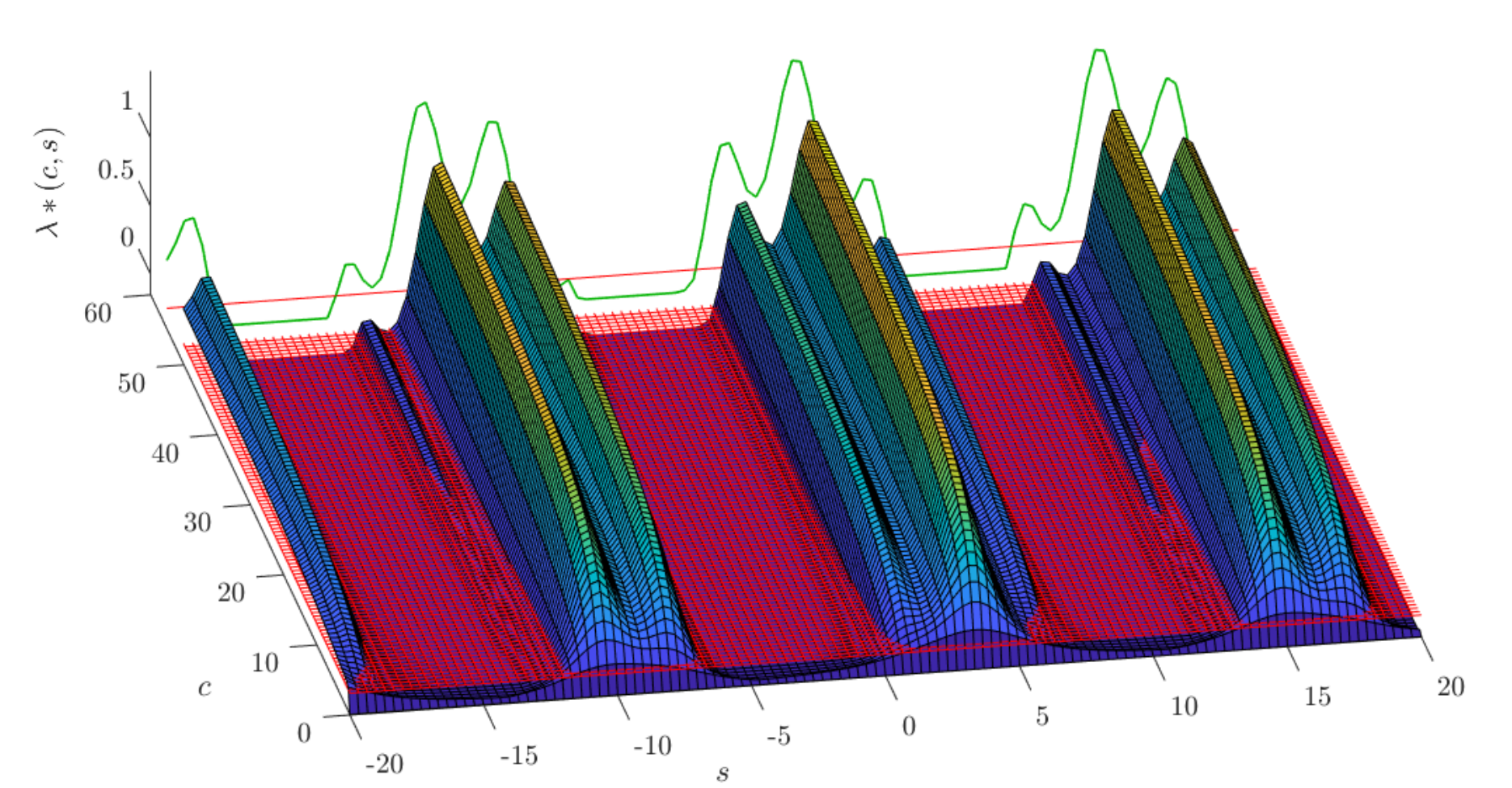}
\label{fig:surface_c_s}
\vspace{-0.5cm}
\end{figure}
\par
We point out that the argument of the proofs of Theorem~\ref{4.teorc_infty}(ii)\&(iii)
relies on showing that $d^\infty(s)$ is strictly positive or negative whenever
$\lb^\infty(s)$ is strictly positive or negative, where
$\lb^\infty(s):=\lb^*(2\,\G_\infty,\,p_s-(\G_\infty)^2)$
is the bifurcation value associated  to \eqref{4.ecucaras} by Theorem~\ref{2.teorlb*}.
This fact provides two methods to identify the values of
$s\in\R$ at which \eqref{4.ecucaras} is in \hyperlink{CA}{\sc cases A} or
\hyperlink{CC}{\sc C} (with tracking or tipping: see Remark \ref{3.notaABC}):
one can numerically calculate $\lb^\infty(s)$, which is rather computationally expensive;
or calculate $d^\infty(s)$, which is a considerably more economic alternative.
In the upper panel of Figure~\ref{fig:theorem3_2},
the attractor-repeller pair $(\wit a, \wit r)$ of $x'=-x^2+p(t)$ is depicted on
the plane $(s,y)$ for $s\in[-40,40]$, and the values of $s\in\R$ for which $d(s)>0$
are highlighted in thick red on the axis $y=0$. The lower panel shows the graphs of
$\lb^\infty$ (in green) and $d^\infty(s)$ (in magenta). Of course,
the two curves have the same signs. We recall once more that, when this sign is positive
(resp.~negative), we can assure the tipping (resp. the tracking) for \eqref{4.ecuinis}$_c$
if $c>0$ is large enough.
\par
\begin{figure}[]
\caption{Total tipping on the hull for
 $y'=-\big(y-2\,\G_c(t)\big)^2+p_s(t)$ for large enough $c$.}
\includegraphics[trim={1cm 9.6cm 1cm 10.2cm},clip,width=\textwidth]{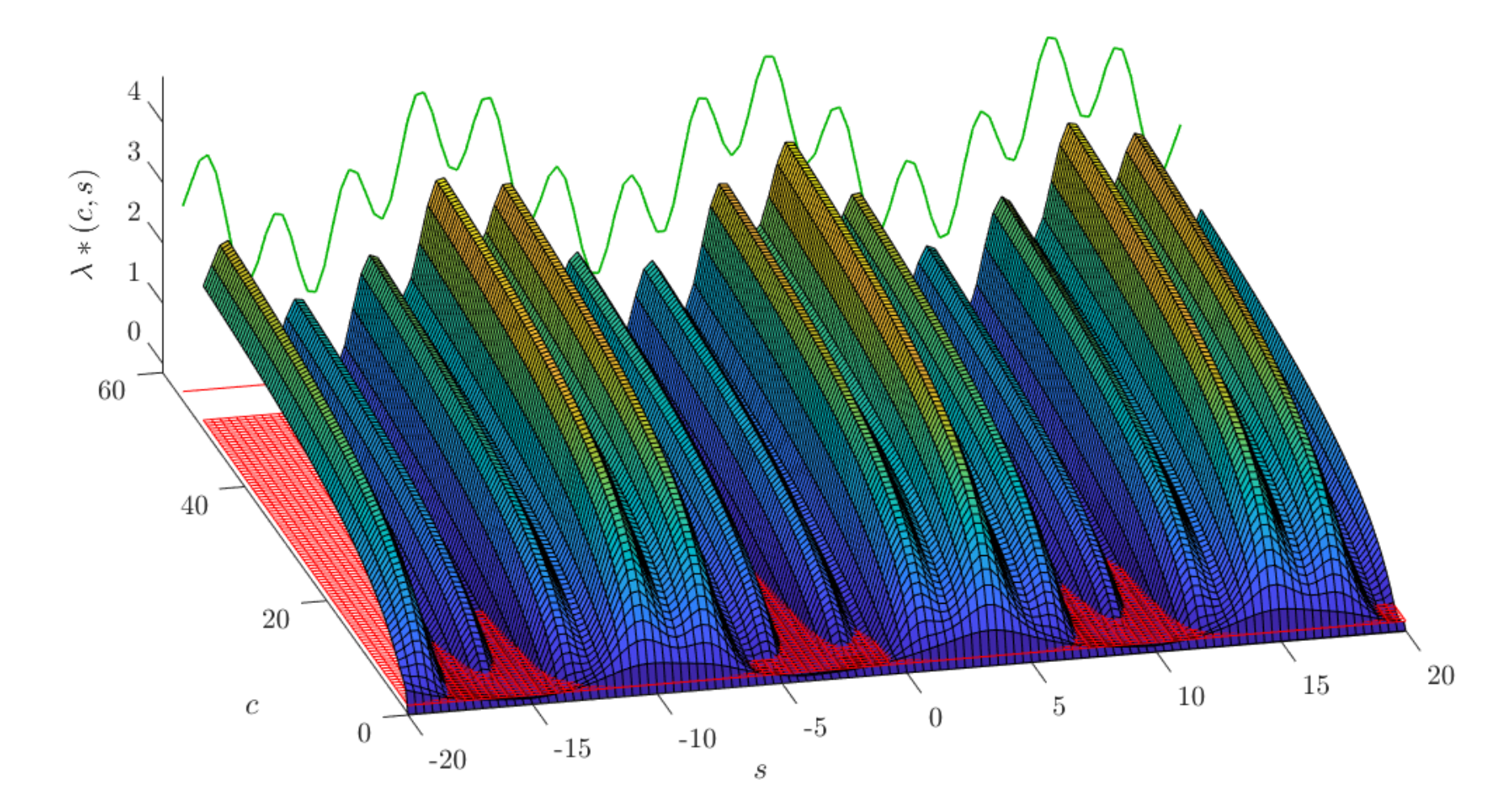}
\label{fig:total tipping}
\vspace{-0.7cm}
\end{figure}
In Figure \ref{fig:surface_c_s}, we show how the bifurcation curve
$\lb^\infty(s)$ (green curve) of $y'=-(y-\G_\infty(t))^2+p_s(t)$ depending on the
variation of $s\in[-20,20]$ seems to be rapidly approached by the bifurcation curve
$\lb_*(c,s)$ of \eqref{4.ecuinis}$_c$ as $c$ increases:
$\lb^\infty(s)$ is very similar to $\lb_*(50,s)$.
\par
Coming back to our notion of partial tipping, what we have in this example is the following.
For a small value of $c_0>0$, we have tracking for all the equations \eqref{4.ecuinis}$_{c_0}$ for
$s\in\R$. More precisely, if $c_0$ is small enough, we have $\lb_*(c_0,s)<-\ep$ for all
$s\in\R$ and an $\ep>0$. From the hull definition and the continuity of
$s\mapsto\lb_*(c_0,s)$ guaranteed by Theorem \ref{2.teorconlb1},
it can be easily deduced that this means tracking for all the equations corresponding
to $c_0$ and any $q$ in the hull of $p$: we have total tracking on the hull.
This means that the whole hyperbolic
copies of the base existing for the future and past families of equations
on the hull are connected
by the hyperbolic families existing for $c_0$. But at a certain value $c_1>c_0$, $\lb(c,s)$
is no longer negative for all $s\in\R$, and for $c_2>c_1$,
it takes positive and negative values at non degenerate intervals.
The functions $\wit a_s$ and $\wit r_s$ for which $\lb_*(c_2,s)$
is negative, which are contained in the hyperbolic copies of the base for the past family,
approach the hyperbolic families of the future as time increases (by the action of
\eqref{4.ecuinis}$_{c_2}$). But the function $\wit a_s$ for which $\lb_*(c_2,s)$ is positive
gives rise to unbounded solutions (always under the action of
\eqref{4.ecuinis}$_{c_2}$). Therefore, the global connection is lost. This is the phenomenon
which we have called partial tipping on the hull. Observe that the continuity of the
function $s\to\lb_*(c_2,s)$ guaranteed by Theorem \ref{2.teorconlb1} ensures that
\hyperlink{CB}{\sc case B} also occurs for some $s\in\R$ in this situation.
(Incidentally, observe also that Figure \ref{fig:surface_c_s}
shows the existence of many values of $s\in\R$ such that tracking
occurs for the all the systems in $c$-parametric family \eqref{4.ecuinis}, since
$\lb(c,s)<0$ for all $c\ge 0$.)
\par
Figure \ref{fig:surface_c_s} also
indicates that partial tipping persists for $c\in(c_1,\infty)\cup\{\infty\}$.
The simple modification of changing $\G$ by $2\,\G$ gives rise
to an always positive $d^\infty$ (see the next paragraph to
understand this phenomenon),
and hence to an always positive $\lb^\infty$.
This fact combined with the previous arguments of continuity of
$s\mapsto\lb_*(c,s)$ for a fixed $c$
means the existence of a large enough value of
$c_3$ (for $2\Gamma$)
such that the tipping is total on the hull for $c>c_3$.
See Figure \ref{fig:total tipping}.
\par
We want to insist in the fact that Theorem \ref{4.teorc_infty} is the key to talk
about partial and total tipping (or total tracking) on the hull for the family
\eqref{4.ecucaras}, which corresponds to $c=\infty$ (and hence also for
the family \eqref{4.ecuinis}$_c$ if $c$ is large enough).
As explained in Remark \ref{4.notaGp}, a suitable choice of $\G$ once $p$ is fixed
can determine this dynamical situation: partial tipping on the hull occurs when
$s\mapsto\wit a(s)-\wit r(s)$ takes values greater and smaller than
$\gamma^+-\gamma^-$; when $\gamma_+-\gamma_-<0$, there is total tracking
on the hull for \eqref{4.ecucaras}; and, a large enough value of
$\gamma^+-\gamma^-$ guarantees the occurrence of total tipping.
\par
The phenomenon that we have described admits also a
different interpretation: the change of variable $s+t=l$
transforms \eqref{4.ecuinis} and \eqref{4.ecucaras} into
\begin{equation}\label{4.ecuinifinal}
 y'=-(y-\G_c(l-s))^2+p(l),
\end{equation}
and
\begin{equation}\label{4.ecucarafinal}
 y'=-(y-\G_\infty(l-s))^2+p(l)\,,
\end{equation}
respectively. Observe that \eqref{4.ecucarafinal}
is obtained from \eqref{4.ecuinifinal} by taking limits as $c\to\infty$.
Therefore, one can read Figures \ref{fig:theorem3_2} and \ref{fig:surface_c_s} as a
characterization of the dynamical scenario for
\eqref{4.ecuinifinal}$_c$ depending on $s$ for $c$ sufficiently large.
In particular, a time shift of the connecting function $\G$
can change the scenario, from the occurrence to the absence of rate-induced tipping.
\par
We close this section by recalling that, recently, other notion of partial tipping
has been described for some switched predator-prey models, in \cite{altw}:
tipping as the climate varies occurs or not depending on the initial point of the phase space.
The model is given by a Carath\'{e}odory equation, which can be understood as a limit of equations
with bounded and uniformly continuous
coefficients, as \eqref{4.ecucaras} is the limit of \eqref{4.ecuinis}$_c$ as $c\to\infty$.
In this way, our analysis of partial tipping on the hull for \eqref{4.ecucaras} is, to some extent,
related to that of \cite{altw}.
\section{Approaching $\G$ by piecewise uniformly continuous functions}\label{5.sec}
Throughout this whole section, and unless otherwise indicated,
$\G\colon\R\to\R$ will be a (bounded and uniformly)
continuous function such that there
exist the real limits $\gamma_\pm:=\lim_{t\to\pm\infty}\G(t)$.
From this map, we define $\G_{\pm\infty}$ as in equations
\eqref{4.ecucara+} and~\eqref{4.ecucara-}, and $\G_c$ as in equation
\eqref{4.ecuini}. We also define, for $c\in\R$
and $h>0$,
\[
\begin{split}
 \G_c^0(t):=\G_c(t)\,,\ \ \qquad&\G_{\pm\infty}^0(t):=\G_{\pm\infty}(t)\,,\\
 \G_c^h(t):=\G(cjh) \quad\;
 \text{if }\;t\in\,&[\,jh\,,(j+1)h\,)\text{ for }j\in\Z\,,\\
 \G^h_\infty(t):=\left\{\!\begin{array}{ll}
 \gamma_-&\quad\text{if $t< 0$}\,,\\
 \gamma_0&\quad\text{if $0\le t<h$}\,,\\
 \gamma_+&\quad\text{if $t\ge h$}\,,
 \end{array}\right.\qquad&\G^h_{-\infty}(t):=\left\{\!\begin{array}{ll}
 \gamma_+&\quad\text{if $t< 0$}\,,\\
 \gamma_0&\quad\text{if $0\le t<h$}\,,\\
 \gamma_-&\quad\text{if $t\ge h$}\,.\end{array}\right.
 \end{split}
\]
Note that, if $h>0$, then $\G_c^h\in BPUC_{\R\!-\!\Delta_h}(\R,\R)$ if $c\in\R$
for $\Delta_h:=\{jh\mid j\in\Z\}$, and that $\G_{\pm\infty}^h\in BPUC_{\R-\{0,h\}}(\R,\R)$.
We fix (also for the whole section) a BPUC function
$p\colon\R\to\R$ and consider the equations
\begin{equation}\label{5.ecutro}
 y'=-(y-\G_c^h(t))^2+p(t)\,.
\end{equation}
We will represent by \eqref{5.ecutro}$_{c,h}$ the equation~\eqref{5.ecutro}
for fixed $c\in\R\cup\{\pm\infty\}$ and $h\ge 0$, and by $y_{c,h}(t,s,y_0)$
its maximal solution with value $y_0$ at $t=s$. Note that
\eqref{5.ecutro}$_{0,h}$,
\eqref{5.ecutro}$_{c,0}$,~\eqref{5.ecutro}$_{\infty,0}$  and
\eqref{5.ecutro}$_{-\infty,0}$ respectively coincide with
\eqref{4.ecu0},~\eqref{4.ecuini}$_c$,~\eqref{4.ecucara+} and \eqref{4.ecucara-},
that~\eqref{5.ecutro}$_{\pm\infty,h}$ play the role of limit equations
for~\eqref{5.ecutro}$_{c,h}$ when $c\to\pm\infty$ for any fixed $h\ge 0$,
that~\eqref{5.ecutro}$_{c,0}$ plays the role of limit
equation for~\eqref{5.ecutro}$_{c,h}$ when $h\to 0$ for any fixed
$c\in\R\cup\{\pm\infty\}$, and that \eqref{4.ecu+} and
\eqref{4.ecu-} are the future and past equations of
 \eqref{5.ecutro}$_{c,h}$ for all $c\in\R-\{0\}$ and $h\ge 0$.
In this formulation, $c$ is again
the {\em rate} of the {\em transition function} $\G_c$. The notions of tipping rate and
transversal tipping rate of Definition \ref{4.defiRIT} are extended
without changes to the newly presented context of families \eqref{5.ecutro}$_{c,h_0}$
(given by piecewise uniformly
continuous transition functions $\G_c^h$) for a fixed $h_0\ge 0$.
\begin{nota}\label{5.notaABC}
Theorem \ref{3.teorexit} establishes a fundamental consequence
of the existence of an attractor-repeller pair
for~\eqref{4.ecup} for the dynamics
induced by~\eqref{5.ecutro}$_{c,h}$:
the existence of the special functions $\ma_{c,h}$ and $\mr_{c,h}$
provided by Theorem~\ref{2.teoruno} for any $c\in\R\cup\{\pm\infty\}$
and for any $h\ge 0$. In particular, under Hypothesis~\ref{4.hipo},
the description made in Remarks \ref{3.notaABC} and \ref{5.notasigno}.1
of the dynamics in \hyperlink{CA}{\sc case A} (or tracking), \hyperlink{CB}{B} and
\hyperlink{CC}{C} (or tipping) for $h=0$ and
$c\in(\R-\{0\})\cup\{\pm\infty\}$ is also valid for any $h>0$.
\end{nota}
As indicated in Remark \ref{5.notasigno}.2 for the continuous case,
the information provided in the previous remark can be partially
extended to some situations in which Hypothesis \ref{4.hipo} does not hold, but
instead an attractor-repeller pair for one of the equations \eqref{5.ecutro}$_{\pm\infty,0}$
exists. It establishes the existence of a local pullback attractor
and a local pullback repeller of \eqref{5.ecutro}$_{c,h}$ for $\pm c>0$ and $h\ge 0$,
respectively connecting with the attractor and repeller
of \eqref{5.ecutro}$_{\pm\infty,0}$ as times decreases and increases,
as well as the behavior of the rest of (at least) half-bounded solutions.
If both equations have an attractor-repeller pair, the comments in Remarks
\ref{5.notaABC} also apply.
\begin{prop}\label{5.proplim}
\begin{itemize}
\item[\rm(i)] Assume that \eqref{4.ecucara+}$\,=\,$\eqref{5.ecutro}$_{\infty,0}$
has an attractor-repeller pair,
$(\wma_\infty,\wmr_\infty)$.
Then, there exist the functions $\ma_{c,h}$ and $\mr_{c,h}$
associated to~\eqref{5.ecutro}$_{c,h}$ by Theorem~{\rm \ref{2.teoruno}} for
$c\in(0,\infty)\cup\{\infty\}$ and $h\ge 0$,
and they satisfy: $\lim_{t\to-\infty}|\ma_{c,h}(t)-\wma_\infty(t))|=0$,
$\lim_{t\to-\infty}|y_{c,h}(t,s,y_0)-\wmr_\infty(t)|=0$ whenever
$y_0<\ma_{c,h}(s)$, $\lim_{t\to+\infty}|\mr_{c,h}(t)-\wmr_\infty(t)|=0$,
and $\lim_{t\to+\infty}|y_{c,h}(t,s,y_0)-\wma_\infty(t)|=0$ whenever
$y_0>\mr_{c,h}(s)$.
\item[\rm(ii)] Assume that \eqref{4.ecucara-}$\,=\,$\eqref{5.ecutro}$_{-\infty,0}$
has an attractor-repeller pair, $(\wma_{-\infty},\wmr_{-\infty})$.
Then, there exist the functions $\ma_{c,h}$ and $\mr_{c,h}$
associated to \eqref{5.ecutro}$_{c,h}$ by Theorem~{\rm \ref{2.teoruno}} for
$c\in\{-\infty\}\cup(-\infty,0)$
and $h\ge 0$, with
$\lim_{t\to-\infty}|\ma_{c,h}(t)-\wma_{-\infty}(t)|=0$,
$\lim_{t\to-\infty}|y_{c,h}(t,s,y_0)-\wmr_{-\infty}(t)|=0$ whenever
$y_0<\ma_{c,h}(s)$,
$\lim_{t\to+\infty}|\mr_{c,h}(t)-\wmr_{-\infty}(t)|=0$, and
$\lim_{t\to+\infty}|y_{c,h}(t,s,y_0)-\wma_{-\infty}(t)|=0$ whenever
$y_0>\mr_{c,h}(s)$.
\item[\rm(iii)] In both cases, the solutions $\ma_{c,h}$ and $\mr_{c,h}$
are respectively locally pullback attractive and locally
pullback repulsive.
\item[\rm(iv)] Assume that \eqref{4.ecucara-} and \eqref{4.ecucara+}
have attractor-repeller pairs. If $\ma_{c,h}$ and $\mr_{c,h}$ are globally defined
and different, then they are uniformly separated, and hence
$(\wma_{c,h},\wmr_{c,h}):=(\ma_{c,h},\mr_{c,h})$ is an attractor-repeller pair for
\eqref{5.ecutro}$_{c,h}$. Consequently, if the equation \eqref{5.ecutro}$_{c,h}$
does not have hyperbolic solutions, it has at most one bounded solution.
\end{itemize}
\end{prop}
\begin{proof}
We assume the existence of an attractor-repeller pair
$(\wma_\infty,\wmr_\infty)$ for~\eqref{4.ecucara+}, and take $c\in(0,\infty)
\cup\{+\infty\}$ and $h\ge 0$. We repeat the arguments of the
proof of Theorem~\ref{3.teorexit}(i)\&(ii) taking $(\wma_\infty,\wmr_\infty)$
instead of $(\wma_-,\wmr_-)$ as starting point, and
working just on $(-\infty,0]$. This requires to
make now use of Proposition~\ref{2.proppers} instead of
Proposition~\ref{2.proppersiste}, which provides
\eqref{3.att3} just for $0\ge t\ge s$. In this way, we prove the existence
of $\ma_{c,h}$, as well as the limit behavior as time decreases of all the
functions which are bounded at $-\infty$. To check the rest of the assertions
in (i), we work on $[0,\infty)$.
The proof of (ii) is analogous, and
the proofs of (iii) and (iv) repeat those of (iii) and (iv) in
Theorem~\ref{3.teorexit}.
\end{proof}

\subsection{The bifurcation map $\lb_*(c,h)$}\label{5.subseclb}
For each $c\in\R\cup\{\pm\infty\}$ and $h\ge 0$, we represent by $\lb_*(c,h)$ the value
of the parameter associated to
\eqref{5.ecutro}$_{c,h}$ by Theorem~\ref{2.teorlb*}; that is, the bifurcation value of
$y'=-(y-\G_c^h(t))^2+p(t)+\lb$,
\begin{equation*}\label{4.deflb*ch}
 \lb_*(c,h):=\lb^*\left(2\,\G^h_c,p-({\G^h_c})^2\right).
\end{equation*}
\begin{teor}\label{5.teorcon}
Let $\lb_*\colon(\R\cup\{\pm\infty\})\times[0,\infty)\to\R$ be defined as above. Then,
\begin{itemize}
\item[\rm(i)] the map $\lb_*$ is bounded: it takes values in $[-||p||,||p||+||\G||^2]$.
\item[\rm(ii)] The map $\lb_*$ is jointly continuous.
\end{itemize}
\end{teor}
\begin{proof}
(i) This assertion is a consequence of the first statement
of Theorem~\ref{2.teorlb*}, since $\n{\G}\ge\n{\G_c^h}\ge\n{\G_{\pm\infty}^h}$
if $c\in\R-\{0\}$ and $h\ge 0$, and $\n{\G}\ge\n{\G_0^h}$
if $h\ge 0$.
\smallskip\par
(ii) We will prove (ii) in three steps.
\smallskip\par
{\sc Step 1}. First we will check the continuity at
$(c_0,h_0)$ with $c_0\in\R\cup\{\pm\infty\}-\{0\}$ and $h_0\ge 0$.
We begin by assuming $c_0>0$.
Let us take a sequence $((c_k,h_k))_{k\ge 1}$ with limit $(c_0,h_0)$,
and assume without restriction that $c_k\ge c_0/2$ for all $k\ge 1$.
What we do in the following paragraphs reproduces the ideas of the
proof of Theorem \ref{3.teorcrec}(ii), where we also prove a continuity
property for the bifurcation function. The reader is referred there for
the details which we omit here.
\par
We will associate below a suitable time $t_\ep>0$ to any $\ep>0$. Once fixed,
we represent the index associated by Theorem~\ref{2.teorlb*} to
\[
 y'=-\big(y-(\G^h_c)^\ep(t)\big)^2+p(t)\,, \quad \text{where}
\quad (\G^h_c)^\ep(t):=\left\{\!\begin{array}{ll}
 \gamma_-&\quad\text{if $t<-t_\ep$}\,,\\
 \G_c^h(t)&\quad\text{if $-t_\ep\le t< t_\ep$}\,,\\
 \gamma_+&\quad\text{if $t\ge t_\ep$}\end{array}\right.
\]
as $\lb_\ep(c,h)$. To check that $\lim_{k\to\infty}\lb_*(c_k,h_k)=\lb_*(c_0,h_0)$ we will
prove that, given $\ep>0$, there exists $t_{\ep}>0$ such that
\begin{itemize}
\item[\bf 1] $|\lb_*(c_k,h_k)-\lb_\ep(c_k,h_k)|<\ep$
for any $k$, including $k=0$; and,
\item[\bf 2] $\lim_{k\to\infty}\lb_\ep(c_k,h_k)=\lb_\ep(c_0,h_0)$.
\end{itemize}
\par
Let us fix $\ep>0$ and prove {\bf 1}. Theorem~\ref{2.teorconlb1} ensures that, for each
$k\ge 1$, there exists $\delta_\ep>0$ such that, if
$\big\|\G_{c_k}^{h_k}-(\G_{c_k}^{h_k})^\ep\big\|\le\delta_\ep$,
then $|\lb_*(c_k,h_k)-\lb_\ep(c_k,h_k)|<\ep$. The goal is finding
$t_\ep$ and hence $\delta_\ep$ such that this bound works for all $k\ge 0$.
We look for $\kappa>0$ and
$\eta>0$ such that $c_k\ge\kappa$ and $h_k\le\eta$ for any $k\ge 0$,
and for $t_\ep=t_\ep(\ep,\kappa,\eta)>\eta$ such that
$|\G(t)-\gamma_-|\le\delta_\ep$ if $t\le-\kappa\,t_\ep$ and
$|\G(t)-\gamma_+|\le\delta_\ep$ if $t\ge\kappa\,t_\ep$.
If $t\le-t_\ep$, then $c_kt\le-\kappa\,t_\ep$, so that
$|\G^0_{c_k}(t)-\gamma_-|\le\delta_\ep$; and if, in addition, $h_k>0$ and
$t\in[jh_k,(j+1)h_k)$, then $jh_kc_k\le c_kt\le-\kappa\,t_\ep$, so that
$|\G^{h_k}_{c_k}(t)-\gamma_-|\le\delta_\ep$.
If $t\ge t_\ep$,
then $c_k\,t\ge\kappa\,t_\ep$, so that
$|\G^0_{c_k}(t)-\gamma_+|\le\delta_\ep$; and if, in addition, $h_k>0$ and
$t\in[jh_k,(j+1)h_k)$, then $jh_kc_k\ge c_kt\le\kappa\,t_\ep$, so that
$|\G^{h_k}_{c_k}(t)-\gamma_+|\le\delta_\ep$.
Hence, the time $t_{\ep}$ is fixed
once $\ep>0$ is fixed, and {\bf 1} is proved.
\par
Let us fix $\ep>0$, which determines $t_\ep>0$ and hence $\lb_\ep$.
To prove {\bf 2}, it suffices to check that
\begin{itemize}
\item[{\bf 2.1}] given $\lb<\lb_\ep(c_0,h_0)$, there exists $k_0$ such that if
$k\ge k_0$ then $\lb\le\lb_\ep(c_k,h_k)$,
\item[{\bf 2.2}] given $\lb>\lb_\ep(c_0,h_0)$, there exists $k_0$ such that
if $k\ge k_0$ then $\lb\ge\lb_\ep(c_k,h_k)$.
\end{itemize}
\par
The proof of {\bf 2.1} reproduces without changes that of point {\bf 1}
of Theorem \ref{3.teorcrec}(ii).
The same happens with the idea to prove {\bf 2.2} and point {\bf 2}
of that theorem. We take $\bar\lb>\lb_\ep(c_0,h_0)$, so that the equation
\begin{equation}\label{5.ecuchlb}
 y'=-\big(y-(\G_{c}^{h})^\ep(t)\big)^2+p(t)+\bar\lb\,,
\end{equation}
corresponding to $(c,h)=(c_0,h_0)$
has an attractor-repeller pair $(\wma_0,\wmr_0)$. The goal is to check
the existence of at least one bounded solution for \eqref{5.ecuchlb}$_{c_k,h_k}^\ep$
if $k$ is large enough, what we
achieve by checking the existence corresponding functions
$\ma_\ep^k$ and $\mr_\ep^k$ given by Theorem~\ref{2.teoruno} at least on
$(-\infty,t_{\ep}]$ and $[t_{\ep},\infty)$, with
$\ma_\ep^k(t_{\ep})\ge\mr_\ep^k(t_{\ep})$.
\par
Observe that the coefficients of the equations~\eqref{5.ecuchlb}$^\ep_{c_k,h_k}$
are common for any $k\ge 0$ outside the interval $[-t_\ep,t_\ep]$. This fact allows
us to repeat the procedure followed to prove {\bf 2} in Theorem \ref{3.teorcrec}(ii)
in order to check the previous assertion.
This completes {\sc step 1} for $c_0>0$, and the proof is
analogous if $c_0<0$.
\smallskip\par
{\sc Step 2}. We will prove that $\lim_{c\to 0}\lb_*(c,0)=\lb_*(0,0)$.
Recall that $\G_c^0(t)=\G(c\,t)$. Let us assume first that $\G$ is $C^1$ with
$\G'\colon\R\to\R$ bounded. For each
$c\in\R$, the change of variables $x=y-\G(c\,t)$ takes equation
$y'=-(y-\G(c\,t))^2+p(t)+\lb$ to
$x'=-x^2-c\,\G'(c\,t)+p(t)+\lb$, and transforms
bounded solutions in bounded solutions. Therefore, the role of
$\lb_*(c,0)$ does not change. Let us take
$c\ne 0$, and let $b_{c}$ be a bounded solution for
$x'=-x^2-c\,\G'(c\,t)+p(t)+\lb_*(c,0)$. Then
$b_{c}'(t)\le -b_{c}^2(t)+p(t)+|c|\!\n{\G'}+\lb_*(c,0)$,
so that Theorem~\ref{2.teoruno}(v) and Theorem~\ref{2.teorlb*}(i) ensure that
$\lb_*(0,0)\le |c|\!\n{\G'}+\lb_*(c,0)$; that is,
$\lb_*(0,0)-\lb_*(c,0)\le |c|\!\n{\G'}$.
Now, let $b_0(t)$ be a bounded solution of
$x'=-x^2+p(t)+\lb_*(0,0)$, so that
$b_0'(t)\le-x^2+p(t)+|c|\!\n{\G'}+\lb_*(0,0)$.
Reasoning as before, we get $\lb_*(c,0)-\lb_*(0,0)\le |c|\!\n{\G'}$.
Consequently, $|\lb_*(c,0)-\lb_*(0,0)|\le |c|\!\n{\G'}$, which
proves the assertion in this case.
\par
Still in {\sc step 2}, we look for a sequence $(\G_n)_{n\ge 1}$ of bounded $C^1$
functions with bounded derivatives and such
that $\lim_{n\to\infty}\G_n=\G$ uniformly on $\R$.
This can be easily done since $\G$ is asymptotically constant and
$C^1(\mI,\R)$ is dense in $C(\mI,\R)$ for any compact interval $\mI$.
We represent by $\lb_n^*(c,0)$ the parameter
associated to the equation $y'=-(y-\G_n(c\,t))^2+p(t)$ by Theorem~\ref{2.teorlb*}.
Then $|\lb_*(c,0)-\lb_*(0,0)|\le
|\lb_*(c,0)-\lb_n^*(c,0)|+|\lb_n^*(c,0)-\lb_n^*(0,0)|+
|\lb_n^*(0,0)-\lb_*(0,0)|$.
Let us take $\ep>0$. Note that $\sup_{t\in\R}|\G_n(c\,t)-\G(c\,t)|\le
\n{\G_n-\G}$ for any $c\in\R$ (in fact they are equal if $c\ne 0$).
Theorem~\ref{2.teorconlb1} provides
$n_0\in\N$ such that $\n{\G_{n_0}-\G}$ is small enough
as to guarantee that $|\lb^*_{n_0}(c,0)-\lb_*(c,0)|\le\ep/3$
for any $c\in\R$. Besides, we have proved in the previous paragraph that
$|\lb_{n_0}^*(c,0)-\lb_{n_0}^*(0,0)|\le |c|\!\n{\G_{n_0}'}$.
Let us take $c_0>0$ such that if $|c|\le c_0$ then
$|c|\!\n{\G_{n_0}'}\le\ep/3$. Altogether, we have
$|\lb_*(c,0)-\lb_*(0,0)|\le\ep$ if $|c|\le c_0$, and
this completes the second step.
\smallskip\par
{\sc Step 3}. Note now that $\lb_*(0,h_0)=\lb_*(0,0)$, since
$\G_0^{h_0}=\G_0^0\equiv\G(0)$. Therefore, in the third and
last step we will prove that, if the sequence $((c_k,h_k))_{k\ge 1}$ tends to
$(0,h_0)$, with $h_0\ge 0$, then $\lim_{k\to\infty}\lb_*(c_k,h_k)=
\lb_*(0,0)$. We write $|\lb_*(c_k,h_k)-\lb_*(0,0)|\le
|\lb_*(c_k,h_k)-\lb_*(c_k,0)|+|\lb_*(c_k,0)-\lb_*(0,0)|$.
We have proved in the second step that
$\lim_{k\to\infty}|\lb_*(c_k,0)-\lb_*(0,0)|=0$. In addition,
$\lim_{k\to\infty} \n{\G_{c_k}^{h_k}-\G_{c_k}^0}=0$,
since $|c_k h_k j -c_k t|\le|c_kh_k|\to 0$ if $t\in[jh_k,(j+1)h_k)$
for a $j\in\Z$ and $\G$ is uniformly continuous. Hence,
Theorem~\ref{2.teorconlb1} ensures that
$\lim_{k\to\infty}|\lb_*(c_k,h_k)-\lb_*(c_k,0)|=0$,
and this completes the proof of (ii).
\end{proof}
\begin{notas}
1.~Observe that Theorem \ref{2.teorlb*}
and the definition of tipping rate given at the beginning of Section \ref{5.sec}
(which repeats Definition \ref{4.defiRIT}) ensure that, for a fixed value of $h_0\ge 0$,
$c_0$ is a tipping rate for the $c$-parametric family of equations~\eqref{5.ecutro}$_{c,h_0}$
if $\lb_*(c_0,h_0)=0$ and there is $\delta>0$ such that $\lb_*(c,h_0)<0$ either for
$c\in(c_0-\delta,c_0)$ or for $c\in(c_0,c_0+\delta)$; and that the tipping
rate is transversal if, in addition, $\lb_*(c,h_0)>0$ either for $c\in(c_0-\delta,c_0)$ or for
$c\in(c_0,c_0+\delta)$. This characterization combined with the just proved
joint continuity of $\lb_*$ shows that,
at a tipping rate $c_0$, the graph of the continuous map $c\mapsto\lb_*(c,h_0)$
reaches the horizontal axis coming from negative values to the left side or to
the right side of $c_0$; and it crosses the horizontal axis transversally at $c_0$
if the tipping rate is transversal.
\par
2.~Assume that the family \eqref{4.ecuini}$_c$$\,=\,$\eqref{5.ecutro}$_{c,0}$
has a transversal tipping rate at $c_0$, passing
from \hyperlink{CA}{\sc case A} to \hyperlink{CC}{\sc C} as $c$ increases.
This means the existence of $\delta>0$ such that $\lb_*(c,0)<0$ for $c\in(c_0-\delta,c_0)$
and $\lb_*(c,0)>0$ for $c\in(c_0,c_0+\delta)$. In particular, $\lb_*(c_0,0)=0$.
The continuity of $\lb_*(c,h)$ ensures the existence of $h_0>0$ such that
$\lb_*(c_0-\delta,h)>0$ and $\lb_*(c_0+\delta,h)<0$ for every $h\in[0,h_0]$.
Therefore, $c(h):=\min\{c\in(c_0-\delta,c_0+\delta)\,|\;\lb_*(c,h)=0\}$
is a (non necessarily transversal) tipping rate of the $c$-parametric family
\eqref{5.ecutro}$_{c,h}$, and in addition $\lim_{h\to 0^+}c(h)=c_0$. In consequence,
every transversal tipping rate of \eqref{4.ecuini}$_c$ can be approximated by
tipping rates of the piecewise continuous transition equations~\eqref{5.ecutro}$_{c,h}$
as $h\to 0^+$.
The other type of rate-induced transversal tipping leads to the same conclusion.
\end{notas}
The next theorem, which includes and extends Theorem \ref{4.teorc_infty},
combines Hypothesis \ref{4.hipo} with the continuity of $\lb_*$ in order
to analyze the dynamical case of \eqref{5.ecutro}$_{c,h}$
for small and large values of $|c|$ and a fixed $h\ge 0$.
We represent by $x(t,s,x_0)$ the solution of
$x'=-x^2+p(t)$ with value $x_0$ at $t=s$.
\begin{teor}\label{5.teorc_infty}
Assume Hypothesis~{\rm \ref{4.hipo}}, and let $(\wit a_\lb, \wit r_\lb)$
be the attractor-repeller pair for $x'=-x^2+p(t)+\lb$ for $\lb>\lb_*(0)$. Let us fix
$h\ge 0$. Then,
\begin{itemize}[itemsep=2pt]
\item[\rm(i)] there exists $c_{0,h}>0$ such that~\eqref{5.ecutro}$_{c,h}$
is in \hyperlink{CA}{\sc case A} for $c\in(-c_{0,h},c_{0,h})$.
\item[\rm(ii)] If there exists $x(h,0,\wit a(0)+\gamma_--\gamma_0)>\wit r(h)+\gamma_+-\gamma_0$,
then the equation~\eqref{5.ecutro}$_{\infty,h}$ has an attractor-repeller pair
$(\wma_{\infty,h},\wmr_{\infty,h})$. In this case, there
exists a minimum $c_{M,h}\ge0$ such that~\eqref{5.ecutro}$_{c,h}$
is in \hyperlink{CA}{\sc case A} for $c>c_{M,h}$.
\item[\rm (iii)] If $x(h,0,\wit a(0)+\gamma_--\gamma_0)$ does not exits, or if
$x(h,0,\wit a(0)+\gamma_--\gamma_0)<\wit r(h)+\gamma_+-\gamma_0$,
then the equation~\eqref{5.ecutro}$_{\infty,h}$ has no bounded solutions.
In this case, there is a minimum $c_{M,h}^*>0$ such that~\eqref{5.ecutro}$_{c,h}$
is in \hyperlink{CC}{\sc case C} for $c>c_{M,h}^*$. In addition, if
$\wit a(0)+\gamma_-<\wit r(0)+\gamma_+$, then
$\lb_*(\infty,0)=\lb_\infty$, where $\lb_\infty>0$
is the unique value of the parameter such that
$\wit a_{\lb_\infty}(0)-\wit r_{\lb_\infty}(0)=\gamma_+-\gamma_-$.
\item[\rm \rm(iv)] If there exists $x(h,0,\wit a(0)+\gamma_+-\gamma_0)>
\wit r(h)+\gamma_--\gamma_0$, then the equation~\eqref{5.ecutro}$_{-\infty,h}$
has an attractor-repeller pair $(\wma_{-\infty,h},\wmr_{-\infty,h})$. In this case, there
exists a maximum $c_{m,h}\le0$ such that~\eqref{5.ecutro}$_{c,h}$
is in \hyperlink{CA}{\sc case A} for $c<c_{m,h}$.
\item[\rm \rm(v)] If $x(h,0,\wit a(0)+\gamma_+-\gamma_0)$ does not exits, or if
$x(h,0,\wit a(0)+\gamma_+-\gamma_0)>\wit r(h)+\gamma_--\gamma_0$
then the equation~\eqref{5.ecutro}$_{-\infty,h}$ has no bounded solutions.
In this case, there is a maximum $c_{m,h}^*>0$ such that~\eqref{5.ecutro}$_{c,h}$
is in \hyperlink{CC}{\sc case C} for $c<c_{m,h}^*$.
In addition, if $\wit a(0)+\gamma_+<\wit r(h)+\gamma_-$,
then $\lb_*(-\infty,0)=\lb_{-\infty}>0$, where
$\lb_{-\infty}$ is the unique value of the parameter such that
$\wit a_{\lb_{-\infty}}(0)-\wit r_{\lb_{-\infty}}(0)=\gamma_--\gamma_+$.
\end{itemize}
\end{teor}
\begin{proof}
(i) Hypothesis~\ref{4.hipo} ensures that $\lb_*(0,h)<0$ for any $h\ge 0$,
and hence (i) is a trivial consequence of the continuity of $\lb_*$.
\smallskip\par
(ii) The goal is proving that the functions $\ma_{\infty,h}$ and $\mr_{\infty,h}$
associated to \eqref{5.ecutro}$_{\infty,h}$ by Theorem \ref{2.teoruno}
form an attractor-repeller pair if $x(h,0,\wit a(0)+\gamma_--\gamma_0)>
\wit r(h)+\gamma_+-\gamma_0$. In these conditions,
Theorem \ref{2.teorlb*} ensures $\lb_*(\infty,h)<0$, and hence
the continuity established in Theorem \ref{4.teorcon} provides the
value $c_{M,h}$ of statement (ii).
\par
The existence of $(\wit a,\wit r)$ ensures that of the
attractor-repeler pairs $(\wma_\pm,\wmr_\pm):=(\wit a+\gamma_\pm,\wit r+\gamma_\pm)$
for the future and past equations $y'=-(y-\gamma_\pm)^2+p(t)$: see Remark \ref{3.notahipo}.
Let us prove that $\ma_{\infty,h}(t)=\wma_-(t)$ for all $t\le 0$.
First, we observe that the existence of $\ma_{\infty,h}$ on $(-\infty,0]$ is
guaranteed by the existence of a solution of~\eqref{5.ecutro}$_{\infty,h}$
bounded on $(-\infty,0]$, which is the case of
$y_{\infty,h}(t,0,\wma_-(0))$: it coincides with $\wma_-(t)$ for $t\le 0$,
since $\G_{\infty}^h(t)=\gamma_-$ for $t<0$. This fact also proves
that $\ma_{\infty,h}(t)\ge\wma_-(t)$ for $t\le 0$. To prove the
converse inequality, we take $s\le 0$ and $y_0>\ma_{\infty,h}(s)$.
Then, the solution $y_-(t,s,y_0)$ of $y'=-(y-\gamma_-)^2+p(t)$,
which coincides with $y_{\infty,h}(t,s,y_0)$ for $t\le 0$, is unbounded
as $t$ decreases, so that $y_0>\wma_-(s)$ and hence
$\ma_{\infty,h}(s)\ge\wma_-(s)$.
\par
The same argument allows us check that
$\mr_{\infty,h}(t)=\wmr_+(t)$ for $t\ge h$. Note also that
for those values of $t\in[0,h]$ for which $\ma_{\infty,h}(t)$ exists,
it coincides with $x(t,0,\wma_-(0)-\gamma_0)+\gamma_0$.
These previous properties and the existence and
inequality assumed in (ii) ensure that $\ma_{\infty,h}(h)=
x(h,0,\wma_-(0)-\gamma_0)+\gamma_0>\wmr_+(h)=\mr_{\infty,h}(h)$.
According to Remark \ref{2.notabasta},
\eqref{4.ecucara+} has at least two bounded solutions, and hence
the information provided by Remark \ref{3.notaABC} ensures
that it has an attractor-repeller pair.
This completes the proof of (ii).
\smallskip\par
(iii) Under the assumptions on (iii), and according to the proof of (ii),
either $\ma_{\infty,h}(h)$ does not exist or we have
$\ma_{\infty,h}(h)<\mr_{\infty,h}(h)$. This precludes the existence of globally
bounded solutions for \eqref{5.ecutro}$_{\infty,h}$, and hence Theorem \ref{2.teorlb*}
ensures that $\lb_*(\infty)>0$. This fact combined with the hypothesis
$\lb_*(0)<0$ and the continuity of $\lb_*$
ensures the existence of a maximum $c_M^*>0$ with $\lb_*(c_M^*)=0$,
which proves the first assertion in (iii).
\par
Assume now that $h=0$ and that $\wit a(0)+\gamma_-<\wit r(0)+\gamma_+$.
Theorem~\ref{2.teorlb*}(ii) ensures the existence
of a unique value $\lb_\infty>0$ of the parameter with
$\wit a_{\lb_\infty}(0)-\wit r_{\lb_\infty}(0)=\gamma_+-\gamma_-$.
We repeat the arguments of the proof of (ii)
taking as starting point the attractor-repeller pair
$(\wit a_{\lb_\infty},\wit r_{\lb_\infty})$ of $x'=-x^2+p(t)+\lb_\infty$
in order to conclude that the functions $\bar a_{\infty,0}$ and $\bar r_{\infty,0}$
associated to the equation $y'=-(y-\G_\infty^0(t))^2+p(t)+\lb_\infty$ by Theorem \ref{2.teoruno}
satisfy $\bar a_{\infty,0}(0)=\wit a_{\lb_\infty}(0)+\gamma_-=
\wit r_{\lb_\infty}(0)+\gamma_+=\bar r_{\infty,0}(0)$. This ensures that
$y'=-(y-\G_\infty^0(t))^2+p(t)+\lb_\infty$ has a unique bounded solution,
which in turn yields $\lb_*(\infty)=\lb_\infty$, as asserted.
\smallskip\par
(iv)\&(v) The arguments to prove these properties are the analogues of those
previously used.
\end{proof}
\begin{nota}
Theorems \ref{5.teorc_infty} and \ref{3.teorexit} show that
the equation \eqref{5.ecutro}$_{\infty,h}$ (or~\eqref{5.ecutro}$_{-\infty,h}$)
has only a bounded solution if and only if there exists
$x(h,0,\wit a(0)+\gamma_--\gamma_0)>\wit r(h)+\gamma_+-\gamma_0$
(or $x(h,0,\wit a(0)+\gamma_+-\gamma_0)>\wit r(h)+\gamma_--\gamma_0$).
But $\lb_*(\infty,h)=0$ (or $\lb_*(-\infty,h)=0$)
does not take us to any conclusion for large values of $c$ (or $-c$).
\end{nota}
We complete this subsection by adapting to equations
\eqref{5.ecutro} part of the information obtained in Subsection
\ref{3.subsecine}. Observe that the value $\lb_*(0)$ appearing in the next
statement coincides with $\lb_*(0,h)$ for any $h\ge 0$, since
$\G_0^h(t)\equiv\gamma_0:=\Gamma(0)$ and dynamics of $x'=-x^2+p(t)$ and
$y'=-(y-\gamma_0)^2+p(t)$ are identical.
\begin{prop}\label{5.propcrecch}
Let $\lb_*(0):=\lb^*(0,p)$ be the value associated
to $x'=-x^2+p(t)$ by Theorem {\rm \ref{2.teorlb*}}.
\begin{itemize}
\item[\rm(i)] If $p$ is  recurrent, then
$\lb_*(0)\le\lb_*(c,h)$ for all $c\in\R$  and $h\ge0$.
\item[\rm(ii)] If $\G$ is nondecreasing (resp.~nonincreasing), then
\[
 \qquad \lb_*(-c,h)\le\lb_*(0)\le\lb_*(c,h), \quad (\text{resp.~}\lb_*(c,h)
 \le\lb_*(0)\le\lb_*(-c,h))
\]
for all $c>0$ and $h\ge0$. If, in addition, $p$ is recurrent, then:
$\lb_*(0)=\lb_*(-c,h)$ (resp.~$\lb_*(0)=\lb_*(c,h)$) for all $c>0$ and $h\ge0$;
and $\lb_*(0)<\lb_*(c,0)$ (resp.~$\lb_*(c,0)<\lb_*(0)$)
if the equation $x'=-x^2+p(t)+\lb^*(0,p)$ has just one bounded solution.
\end{itemize}
\end{prop}
\begin{proof}
Statement (i) is a direct consequence of Proposition \ref{3.propprec}.
The properties of (ii) follow from Corollary \ref{3.corocrec}(ii),
having in mind that $\G_c^h$ and $-\G_{-c}^k$ are nondecreasing
for $c>0$ if $\G$ is nondecreasing, and that $\G_{-c}^h$ and $-\G_c^k$
are nondecreasing for $c>0$ if $\G$ is nonincreasing.
\end{proof}
The bifurcation map $\lb_*(c,h)$ of~\eqref{5.ecutro}$_{c,h}$
when $\G$ and $p$ are given by \eqref{4.eqpt} and $c,h\in[0,5]$ is
depicted in Figure \ref{fig:surface_c_h}.
Besides the joint continuity of $\lb_*$,
observe that the section map $h\mapsto\lb^*(c,h)$ is
not increasing for a fixed $c>0$, unlike the map $h\mapsto\G^h_c$.
\begin{figure}[]
\caption{Numerical simulation of the bifurcation map $\lb_*(c,h)$
of~\eqref{5.ecutro}$_{c,h}$ for $\G$ and $p$ as in~\eqref{4.eqpt} and
$c,h\in[0,5]$. On the left: the gradient surface represents the graph
of $\lb_*(c,h)$; the red grid identifies the plane $\lb_*=0$: the points
of the surface below this plane correspond to {\sc {Case A}}, the points
above to {\sc {Case C}}, and the points of intersection to {\sc {Case B}}.
The red dashed line is the graph of the bifurcation curve
$\lb_*(c,0)=\lb_*(c)$ of the family~\eqref{4.ecuini}, whereas the solid
green line is the graph of the  bifurcation curve $\lb_*(\infty,h)$
of~\eqref{5.ecutro}$_{\infty,h}$, represented at $c=6$ for convenience.
On the right: a projection of the same picture on the plane $c=0$.
}
\includegraphics[trim={1.8cm 8.3cm 2.5cm 8cm},clip,width=0.5\textwidth]{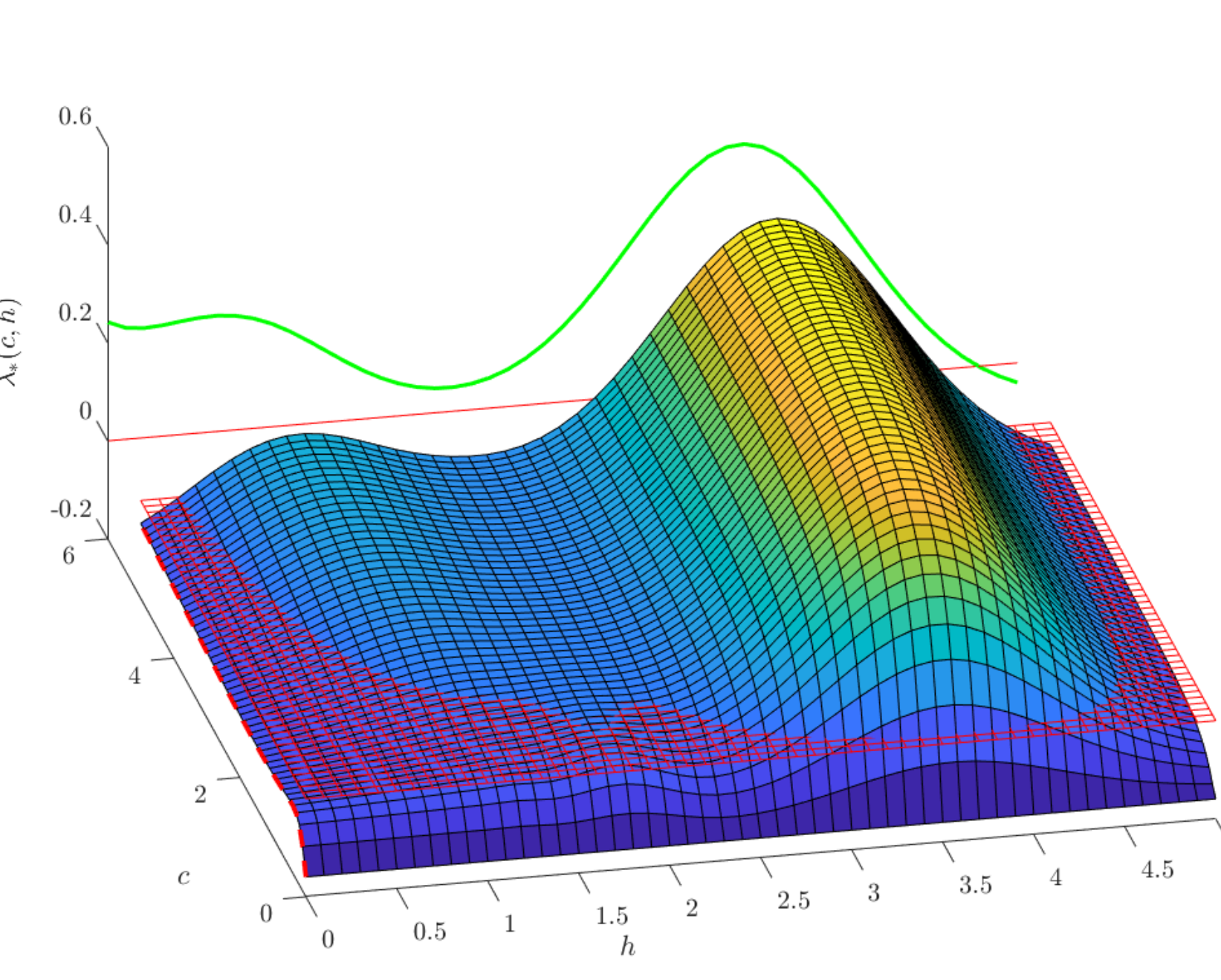}
\includegraphics[trim={2cm 8.3cm 2.8cm 8cm},clip,width=0.45\textwidth]{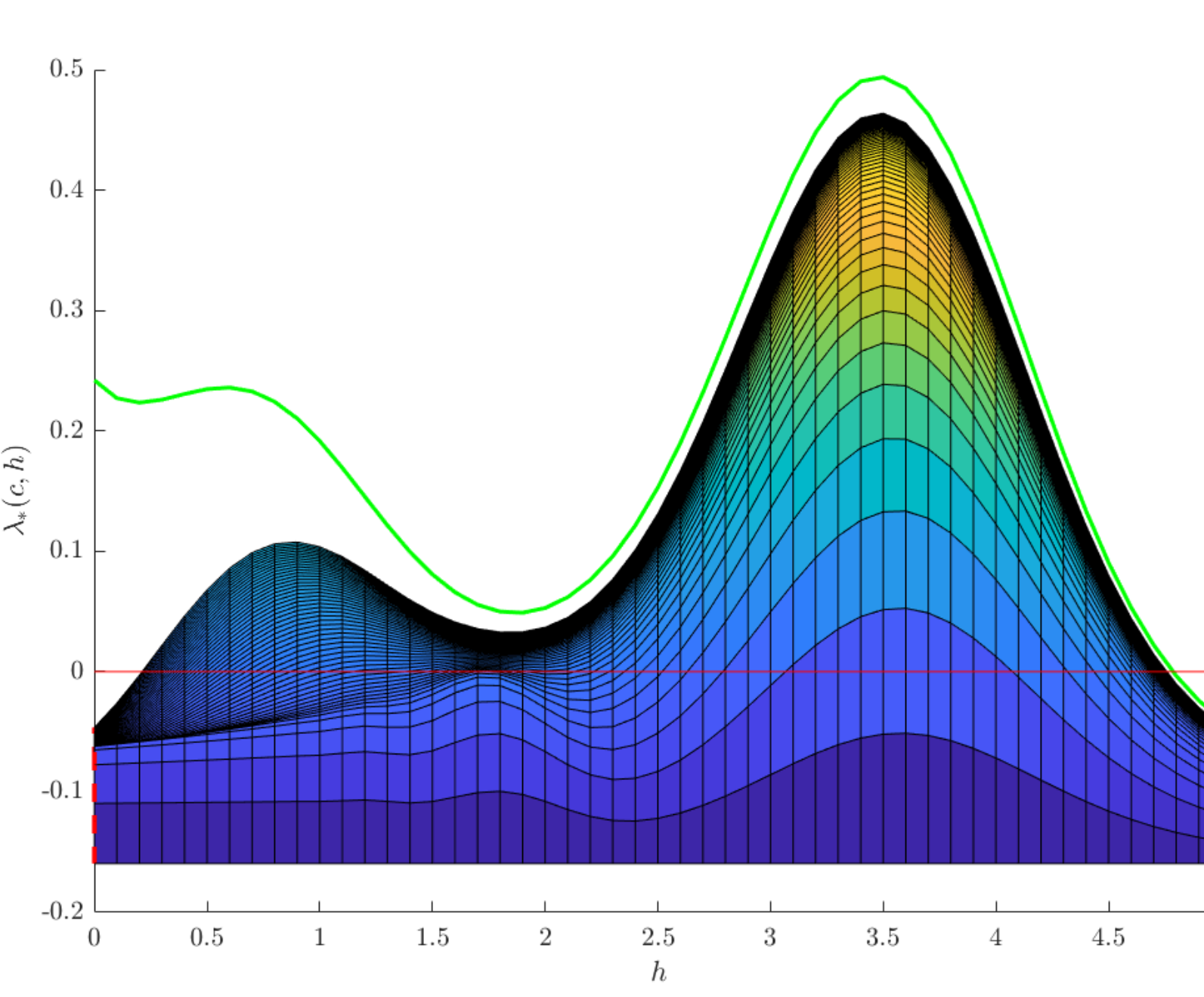}
\label{fig:surface_c_h}
\end{figure}

\appendix
\section{Compactness of the Hull and continuity of the flow}\label{appendix}
Hereby, we recall some facts on nonautonomous equations of the
type \eqref{2.ecucon},
\begin{equation}\label{A.ecucon}
 x'=-x^2+q(t)\,x+p(t)\,,
\end{equation}
where $p$ and $q$ are BPUC functions: see Definition \ref{2.defiBPUC}.
The first objective is proving Theorems \ref{2.teorhyp} and \ref{2.teorlb*},
which rely on Theorem \ref{A.teorCont} (in turn based on Theorem \ref{A.teorHull}).
The second one is to prove Theorem \ref{A.teorLC}, a result on continuous variation of the solutions
with respect to the coefficients which we have used several times.
\par
Let $\Delta=\{a_j\in\R\,|\; j\in\Z\}\subset\R$ be a disperse set (see Subsection \ref{2.subsecAR}).
Recall that $q\colon\R\to\R$ belongs to $BPUC_\Delta(\R,\R)$ if and only it is
right-continuous and
\begin{itemize}
\item[\hypertarget{c1}{{\bf c1}}] there is $c>0$ such that $|q(t)|<c$ for all $t\in\R$;
\item[\hypertarget{c2}{{\bf c2}}] for all $\ep>0$, there is $\delta=\delta(\ep)>0$
such that,
if $t_1,t_2\in (a_j,a_{j+1})$ for some $j\in\Z$ and $t_2-t_1<\delta$, then $|q(t_2)-q(t_1)|<\ep$.
\end{itemize}
Recall also a function $q\in BPUC(\R,\R)$ is a finite sum of a finite number of functions
$q_i\in BPUC_{\Delta_i}(\R,\R)$, for possibly different disperse sets $\Delta_i$.
It is clear that
$BPUC(\R,\R)\subset L^\infty(\R,\R)\subset L^1_{\rm loc}(\R,\R)$.
Recall that $L^1_{\rm loc}(\R,\R)$ is a
complete metric space for the distance defined by
\[
 d(q_1,q_2):=\sum_{k=1}^\infty\frac{1}{2^k}\frac{\int_{-k}^k|q_1(t)-q_2(t)|\,dt}
 {1+\int_{-k}^k|q_1(t)-q_2(t)|\,dt}\;.
\]
In addition, for every $q\in BPUC(\R,\R)$ and $s\in\R$, the shift
$q_s\colon\R\to\R$ defined by $q_s(st):=q(t+s)$, belongs to $L^\infty(\R,\R)$ and
has norm $\n{q}$. We define
\[
 \W_q:=\text{closure}_{L^1_{\rm loc}(\R,\R)}\{q_t\,|\;t\in\R\}\subset
 L^1_{\rm loc}(\R,\R)\cap L^\infty(\R,\R)\,.
\]
The set $\W_q$ is the {\em hull} of $q$ (in $L^1_{\rm loc}(\R,\R)$).
Theorem \ref{A.teorHull} shows that $\W_q$ is a compact metric space, and
that the time-translation map
\[
 \sigma\colon\R\times\W_q\to\W_q\,,\;(t,\w)\mapsto \w_t\,,
 \qquad \text{with $\w_t(s):=\w(t+s)$}
\]
defines a (real) continuous flow on $\W_q$.
Recall that being a flow means that
$\sigma_0=\text{Id}$ and $\sigma_{s+t}=\sigma_t\circ\sigma_s$
for each $s,t\in\R$, where $\sigma_t(\w):=\sigma(t,\w)$.
\begin{teor}\label{A.teorHull}
Let $q\colon\R\to\R$ belong to $BPUC(\R,\R)$. Then
its hull $\W_q$ is a compact subset of $L^1_{\rm loc}(\R,\R)$,
and $\sigma$ defines a continuous flow.
\end{teor}
\begin{proof}
Let $\Delta\subset\R$ be a disperse set, and let us take $q\in BPUC_\Delta(\R,\R)$.
We will first prove the compactness of $\W_q$ in this case.
Note that, since $q$ is bounded, there exists a common bound for
all the shifts $q_s$. Therefore, and according to \cite[Theorem 1]{sell},
to prove the compactness of $\W_\delta$
it suffices to show that given $\ep>0$ and a compact interval
$\mI\subset\R$ there exists
$\delta=\delta(\ep,q,\mI)>0$ such that, for any $s\in\R$,
\begin{equation}\label{A.intI}
 \int_\mI|q_{s}(t+\tau)-q_s(t)|\,dt<\ep\quad\text{whenever }|\tau|<\delta\,.
\end{equation}
Let us write $\Delta=\{a_j\in\R\mid j\in\Z\}$ and define
$h:=\inf_{j\in\Z}(a_{j+1}-a_j)>0$.
We fix $\ep>0$, a non-degenerate interval $\mI=[r_1,r_2]$, and $s\in\R$,
and look for $a_{k+1},a_{k+m-1}\in\Delta$ (depending on $s$), if they
exist, such that $[a_{k+1},a_{k+m-1}]\subseteq(r_1+s,r_2+s)\subset[a_k,a_{k+m}]$.
We set $\wit a_k:=r_1+s$, $\wit a_j:=a_j$ for $j=k+1,\ldots,k+m-1$ and
$\wit a_{k+m}:=r_2+s$, so that $\wit a_{j+1}-\wit a_{j}\le r_2-r_1$ for
$j=k,\ldots,k+m-1$. Then,
\[
 \int_\mI|q_{s}(t+\tau)-q_s(t)|\,dt=\int_{\mI+s}|q(t+\tau)-q(t)|\,dt
 \le\sum_{j=k}^{k+m-1}\int_{\wit a_j}^{\wit a_{j+1}}|q(t+\tau)-q(t)|\,dt\,.
\]
Condition \hyperlink{c1}{\bf c1} gives $\delta_1=\delta_1(\ep,q,\mI)>0$
such that, if $0\le\tau<\delta_1$ and $d\in\R$, then
\[
 \int_{d-\tau}^{d}|q(t+\tau)-q(t)|\,dt<\frac{\ep}{2m}\:.
\]
In addition, \hyperlink{c2}{\bf c2} provides $\delta_2=\delta_2(\ep,q,\mI)\in(0,h]$
such that, if $0<\tau<\delta_2$, then
\[
 |q(t+\tau)-q(t)|<\frac{\ep}{2 m (r_2-r_1)} \quad\text{for all $j\in\Z$ and all
 $t\in(a_j,a_{j+1}-\tau)$}\,.
\]
Let us call $\delta:=\min(\delta_1,\delta_2)$ and fix
$\tau\in[0,\delta)$. For $j\in\{k+1,\ldots,k+m-2\}$ (if there is any),
\[
\begin{split}
 \int_{\wit a_j}^{\wit a_{j+1}}\!\!|q(t+\tau)-q(t)|\,dt
 &=\int_{a_j}^{a_{j+1}-\tau}\!\!|q(t+\tau)-q(t)|\,dt+
 \int_{a_{j+1}-\tau}^{a_{j+1}}\!\!|q(t+\tau)-q(t)|\,dt\\
 &<\frac{\ep}{2 m (r_2-r_1)}\:(a_{j+1}-\tau-a_j)+\frac{\ep}{2m}\le\frac{\ep}{m}\:.
\end{split}
\]
In the case or cases $j=k$ and $j=k+m-1$, the length of the interval
$[\wit a_j,\wit a_{j+1}]$ can be less than $h$, and hence greater than $\tau$.
We proceed in the same way as before if
$\wit a_j\le \wit a_{j+1}-\tau$, getting the bound $\ep/m$; if this is not the case,
we forget about $\int_{\wit a_j}^{\wit a_{j+1}-\tau}$
(which is negative), getting $\ep/(2m)$ as a bound.
It follows that \eqref{A.intI} holds for $0\le\tau<\delta$.
To work with $\tau<0$, we write $\int_{\wit a_j}^{\wit a_{j+1}}=
\int_{\wit a_j}^{\wit a_{j}-\tau}
+\int_{\wit a_j-\tau}^{\wit a_{j+1}}$ and use the same arguments.
This proves the compactness of $\W_q$ in the case $q\in BPUC_\Delta(\R,\R)$.
\par
To extend the result to the general BPUC case, it is enough to observe that
\eqref{A.intI} holds for $q=q_1+\cdots +q_n$ if it holds for every $q_i$.
It is also easy to check that $\sigma$ defines a flow on $\W_q$, and its continuity
follows from \cite[Theorem III.11]{sell2}.
\end{proof}
The function $f(x,t):=-x^2+q(t)\,x+p(t)$ giving rise to \eqref{A.ecucon} is
Lipschitz Carath\'{e}odory whenever $q,p\in BPUC(\R,\R)$.
Recall that a function $f\colon\R\times\R\to\R$ is said to be {\em Lipschitz
Carath\'{e}odory}, which we represent by $f\in\mathfrak{LC}(\R,\R)$, if
\begin{itemize}
\item[-] $f$ is Borel measurable,
\item[-] for every compact interval $\mI\subset\R$ there exists a
function $m^\mI\in L^1_{\rm loc}(\R,\R)$ such that $|f(t,x)|\le m^\mI(t)$ for any $x\in \mI$
and almost every $t\in\R$,
\item[-] for every compact interval $\mI\subset\R$ there exists a function
$l^\mI\in L^1_{\rm loc}(\R,\R)$ such that $|f(t,x_1)-f(t,x_2)|\le l^\mI(t)|x_1-x_2|$ for any
$x_1,x_2\in \mI$ and almost every $t\in\R$.
\end{itemize}
We endow the set $\mathfrak{LC}(\R,\R)$ with the
$\mathcal T_{\Q}$ topology, which is generated by the countable family of seminorms
\[
 n_{[r_1,r_2],s}(f)=\int_{r_1}^{r_2}|f(t,s)|\,dt\quad
 \text{for $\:r_1,\,r_2,\,s\in\Q\:$ with $\;r_1<r_2$}\,,
\]
and for which $\mathfrak{LC}(\R,\R)$
is a locally convex metric space: see e.g.~\cite{sell}.
\par
The results on existence, uniqueness, and basic properties of the solutions of
the initial value problems of equations $x'=f(t,x)$ for
$f\in\mathfrak{LC}(\R,\R)$ are classical: see e.g.~\cite[Chapter 2]{cole}.
But working with $f\in\mathfrak{LC}(\R,\R)$ with the
topology $\mathcal T_\Q$ does not allow to define a {\em continuous\/}
flow from the solutions of the equation via the hull procedure,
which is required to apply techniques from topological
dynamics. One needs to restrict to a suitable subset of functions $f$.
Theorem \ref{A.teorCont} deals with this question when
$f(x,t)$ takes the shape $-x^2+q(t)\,x+p(t)$ with $q,p\in BPUC_\Delta(\R,\R)$.
An in-depth analysis of different topologies and additional conditions
in different spaces of Carath\'{e}odory functions giving rise to
continuous flows appears in \cite{lono2}.
\par
Before stating Theorem \ref{A.teorCont}, we fix some notation.
Given $q,p\in BPUC(\R,\R)$, we can consider the joint hull of $(q,p)$,
\[
 \W_{q,p}:=\text{closure}_{L^1_{\rm loc}(\R,\R^2)}\{(q_t,p_t)\mid t\in\R\}
 \subset L^1_{\rm loc}(\R,\R^2)\,.
\]
For every $(\bar q,\bar p)\in \W_{q,p}$ and $x_0\in\R$,
the map $t\mapsto x(t,\bar q,\bar p,x_0)$
is the unique maximal solution of $x'=-x^2+\bar q(t)\,x+\bar p(t)$ with $x(0)=x_0$;
and $\U\subseteq\R\times\W_{q,p}\times\R$ is the domain of the function $x$.
\begin{teor}\label{A.teorCont}
Let $q,p\colon\R\to\R$ belong to $BPUC(\R,\R)$.~Then,~$\mU$ is an open set,~and
\[
 \Phi\colon\U\subseteq\R\times\W_{q,p}\times\R\to \W_{q,p}\times\R,\quad
 \big(t,(\bar q,\bar p),x_0\big)\mapsto \big((\bar q_t,\bar p_t),
 x(t,\bar q,\bar p,x_0)\big)
\]
defines a continuous local flow on $\W_{q,p}\times\R$,
which is $C^1$ in $x_0$.
\end{teor}
\begin{proof}
The flow properties follow easily from the uniqueness of
the solutions. Since $\W_{q,p}\subset\W_q\times\W_p$, it is compact, and
the continuity of the flow translation follows from those of
$\W_q$ and $\W_p$. Hence, it remains to check the continuity of the
second component of $\Phi$.
\par
Set $f(\w,x):=-x^2+q(t)\,x+p(t)$, and let $\mathrm{Hull}_{\mathcal T_\Q}(f)$ be
the closure in $\mathfrak{LC}(\R,\R)$ endowed with $\T_\Q$
of the set $\{f_t\mid t\in\R\}$ of time-translations of $f$.
There exists a one-to-one correspondence between
$\mathrm{Hull}_{\T_\Q}(f)$ and $\W_{q,p}$: for $(\wit q,\wit p)\in\W_{q,p}$
we define $\wit f(t,x):=-x^2+\wit q(t)\,x+\wit p(t)$ and check that it belongs to
$\mathrm{Hull}_{\T_\Q}(f)$; and given
$\wit f\in \mathrm{Hull}_{\T_\Q}(f)$,
we define $\wit q(t):=\wit f(t,0)$ and $\wit p(t):=\wit f(t,1)+1-\wit q(t)$,
and check that $(\wit q,\wit p)\in\W_{q,p}$.
In addition, it is not hard to check that the bijection $(\wit q,\wit p)\mapsto\wit f$ is
continuous, and hence a homeomorphism between both compact spaces.
\par
On the other hand, the expression of $f$ makes it easy to check that,
for every $j\in\N$, there is $\kappa_j>0$ such that, for a.e.~$t\in\R$,
\[
 \quad|f(t,x_1)-f(t,x_2)|<\kappa_j|x_1-x_2|\quad\text{whenever }x_1,x_2\in[-j,j]\,.
\]
In this situation, \cite[Theorem 5.9(i)]{lono2} (which is formulated for
$\mathrm{Hull}_{\T_\Q}(f)$) shows that
$\mU$ is open as well as the continuity of the second component of $\Phi$.
The $C^1$ character in $x_0$ follows in a standard way from the fact that
the derivative solves the corresponding variational equation
(see \cite[Theorem 2.3.1]{brpi}) combined with the just established
continuity.
\end{proof}
\noindent{\em Proofs of Theorems {\rm \ref{2.teorhyp}} and {\rm \ref{2.teorlb*}}}.
Once established the continuity of the flow induced by
\eqref{2.ecucon} on $\W_{q,p}\times\R$, in Theorem \ref{A.teorCont},
we can repeat the arguments leading to the (long and complex) proof of
\cite[Theorem 3.5]{lnor}, which deals with the analogous properties in the case of
bounded and uniformly continuous functions $q,p$. These arguments allow us to prove
(a)$\Rightarrow$(b), as well as points (i) and (ii) of Theorem \ref{2.teorhyp}
when (a) (and hence (b)) holds.
The proof of the analogue of \cite[Theorem 3.5]{lnor}
requires to apply the first approximation theorem
to a scalar equation of the type $z'=(-2\,\wit b(t)+q(t))\,z-z^2$,
where $\wit b(t)$ is a bounded continuous function, but with
$q\in L^\infty(\R,\R)$ instead of bounded continuous. This is not a problem:
the proof of Theorem III.2.4 in~\cite{hale} works without
changes for this situation.
\par
The assertion (b)$\Rightarrow$(c) of Theorem \ref{2.teorhyp} is trivial, and
(c)$\Rightarrow$(a) can be deduced, for instance, of Theorem \ref{2.teorlb*}(iv),
whose proof is independent of Theorem \ref{2.teorhyp}. The whole proof
of Theorem \ref{2.teorlb*} repeats that of \cite[Theorem 3.6]{lnor}.
The results required there for bounded and uniformly continuous
coefficients $q,p$ have been established in
this paper for $q,p\in BPUC(\R,\R)$.
\qed
\medskip\par
Now we consider, as in the previous sections,
a BPUC function $p\colon\R\to\R$ and a continuous function
$\G\colon\R\to\R$ with finite asymptotic limits $\gamma_\pm:=\lim_{t\to\pm\infty}\G(t)$.
And we define $\G_c^h$ for all $c\in\R\cup\{\pm\infty\}$ and $h\ge 0$
as at the beginning of Section~\ref{5.sec}.
\begin{teor} \label{A.teorLC}
Let us define $f_c^h(t,y):=-(y-\G_c^h(t))^2+p(t)$
for $c\in\R\cup\{\pm\infty\}$ and $h\ge 0$. Then, $f_c^h\in\mathfrak{LC}(\R,\R)$.
In addition, let us denote by $\tilde y_{c,h}(\cdot,0,y_0)$ the solution of
$y'=f_c^h(t,y)$ with $y(0)=y_0$. If the sequence $((c_k,h_k))$ in
$\R\times[0,\infty)$ converges to $(c_0,h_0)$, with $c_0\in\R\cup\{\pm\infty\}$
and $h_0\in(0,\infty)$, and the sequence $(y_k)$ in $\R$ converges to $y_0\in\R$, then
\[
 \lim_{k\to\infty}\tilde y_{c_k,h_k}(t,0,y_k)=\tilde y_{c_0,h_0}(t,0,y_0)
\]
uniformly in $t$ varying in
any compact interval contained in the maximal interval of definition of
$\tilde y_{c_0,h_0}(\cdot,0,y_0)$.
\end{teor}
\begin{proof}
Note first that, if $j\in\N$ and $y_1,y_2\in[-j,j]$, then
\[
 |f_c^h(t,y_1)-f_c^h(t,y_2)|\le (2j+2\n{\G})|y_1-y_2|
\]
for all $c\in\R\cup\{\pm\infty\}$ and $h\ge 0$. This ensures the third of the
conditions ensuring that $f_c^h\in\mathfrak{LC}(\R,\R)$, and the two first ones
are easier to check. Now, let us take a sequence $((c_k,h_k))$ in $\R\times[0,\infty)$
with limit $(c_0,h_0)$, as in the statement.
According to Theorem 5.8(i) in~\cite{lono2}, to prove the last
assertion it suffices to check that
\[
 \lim_{k\to\infty}\int_{r_1}^{r_2}|f_{c_k}^{h_k}(t,y)-f_{c_0}^{h_0}(t,y))|\,ds=0
\]
for all $r_1,r_2,y\in\Q$ with $r_1<r_2$. Clearly, this limiting behavior is guaranteed by
\[
 \lim_{k\to\infty}\int_{r_1}^{r_2}\big|\,\G_{c_k}^{h_k}(t)-\G_{c_0}^{h_0}(t)\big|\,dt=0
 \quad\;\text{for $\,r_1,r_2\in\Q\,$ with $\,r_1<r_2$}\,,
\]
which is hence the property to be proved. In turn, this property follows from
the dominated convergence theorem, since $\G$ is bounded and
$\lim_{k\to\infty}\G_{c_k}^{h_k}(t)=\G_{c_0}^{h_0}(t)$ for almost every $t\in\R$.
The detailed proof of this last assertion is a nice exercise, for which we give some
hints. Given $h>0$ and $t\in\R$, if $j^t$ is the unique integer number with
$t\in[j^th,(j^t+1)\,h)$, we have $j^th\in(t-h,t\,]$. This is the key point to prove that
in the case $c_0=0$ and $h_0\ge 0$, as well as in the case $c_0\in\R-\{0\}$ and $h_0=0$,
the convergence holds for every $t\in\R$. If $c_0\in\R-\{0\}$ and $h_0>0$, the convergence
holds when $t\ne j\,h_0$ for every $j\in\Z$. For $c=\pm\infty$ and $h_0=0$, it holds for $t\ne 0$.
And finally, for $c=\pm\infty$ and $h_0>0$, it holds for every $t\ne h_0$.
\end{proof}
\section{Clarification on the numerical analysis}\label{appendix2}
Hereby, we clarify the way in which we obtain the figures in
Sections \ref{4.subsec_partial_tipping} and \ref{5.subseclb},
corresponding to the differential
equation~\eqref{5.ecutro}$_{c,h}$ (equal to
\eqref{4.ecuini}$_c$ for $h=0$) for
\[
 \G(t):=\frac{2}{\pi}\arctan(t)\,,
 \quad p(t):=0.962-\sin(t/2)-\sin(\sqrt{5}\,t),\quad c\ge0, \quad\text{and}\quad h\ge 0\,.
\]
All the involved equations have been numerically integrated using the
MATLAB function {\tt ode45} with double precision and
the options on the relative and absolute
tolerance respectively set to {\tt RelTol=1e-9} and {\tt AbsTol=1e-9}.
The numerical method used  to compute  $\lb_*$ is based on the bisection
idea outlined in~\cite{lnor}, to which we refer the reader for further details.
In that example, $\G$ is the same, and $p(t):=0.892-\sin(t/2)-\sin(\sqrt{5}\,t)$.
We point out here that the section $\lb_*(c,0)$ of the bifurcation map
$\lb_*$ (see Subsection \ref{5.subseclb}) corresponding to our
present example coincides with $\wit\lb(c)-(0.962-0.895)=\wit\lb(c)-0.067$,
where $\wit\lb(c)$ is the function corresponding to the example in \cite{lnor}:
see Theorem \ref{2.teorlb*}(v). In particular, the detailed justification
given in \cite{lnor}, taken for valid also in what follows,
allows us to ensure that $\lb_*(c,0)<0$ for $c\in[0,0.25]$ (at least).
Hence, Hypothesis \ref{4.hipo} is fulfilled for our coefficients $\G$ and $p$.
\par
We work under the next fundamental assumption, which is
based on a consistent numerical evidence and
which we will explain later: if $\wit\G$ belongs to $BPUC_\Delta(\R,\R)$ for a
disperse set $\Delta$ and satisfies $\|\wit\G\|\le 0.1$,
then the equation $x'=-(x-\wit\G(t))^2+p(t)$ has an attractor-repeller pair.
In these conditions, Hypothesis \ref{3.hipo} and Theorem~\ref{3.teorexit}
guarantee the existence of the (possibly locally defined) solutions $\ma_{c,h}$ and
$\mr_{c,h}$ of \eqref{5.ecutro}$_{c,h}$ for $c\in\R\cup\{\pm\infty\}$ and $h\ge 0$
described in Theorem~\ref{2.teoruno}.
In addition, since the constant $m=3.4$ satisfies the condition required in
Theorem~\ref{2.teoruno} for all the differential equations~\eqref{5.ecutro}$_{c,h}$
(as deduced from $\n{\G_c^h}\le 1$ and $\n{p}\le 3$),
we know that $\ma_{c,h}(t)<3.5$ and $\mr_{c,h}(t)>-3.5$ on their respective domains,
and that any bounded solution, if it exists, takes values in $(-3.5,3.5)$.
\par
We already know that $\lb_*(c,0)<0$ for $c\in[0,0.25]$. Therefore,
$\ma_{c,0}$ and $\mr_{c,0}$ are globally defined hyperbolic solutions for
$c\in[0,0.25]$. We will now check that $\ma_{c,h}$ and $\mr_{c,h}$
are respectively defined on $(-\infty,-35]$ and $[35,\infty)$ whenever
$c\ge 0.25$ (including $c=\infty)$ and $h\in(0,6]$. Let us define
\[
 (\G_c^h)^-(t):=\left\{\begin{array}{ll}
 \G_c^h(t)&\quad\text{if $t<-35$}\,,\\
 \G_c^h(-35)&\quad\text{if $t\ge -35$}\,,\end{array}\right.
 \quad\;\;
 (\G_c^h)^+(t):=\left\{\begin{array}{ll}
 \G_c^h(35)&\quad\text{if $t<35$}\,,\\
 \G_c^h(t)&\quad\text{if $t\ge35$}\,,\end{array}\right.
\]
and observe that $(\G_c^h)^-(t)\in(-1,-0.9)$ and $(\G_c^h)^+(t)\in(0.9,1)$.
In the case of $c=\infty$ and $h\in(0,6]$, these assertions are trivial.
In the remaining cases, they follow
from these facts: given $h>0$ and $t\in\R$, if $j^t$ is the unique
integer number with $t\in[j^th,(j^t+1)\,h)$, then $j^th\in(t-h,t\,]\subseteq(t-6,t]$;
$-1<(2/\pi)\arctan (c\,j^th)<(2/\pi)\arctan (c\,t)\le(2/\pi)\arctan (0.25\cdot(-35))<-0.9$
if $t\le -35$ and $c\ge 0.25$; and $1>(2/\pi)\arctan (c\,j^th)>
(2/\pi)\arctan (c\,(t-6))\ge(2/\pi)\arctan (0.25\cdot (35-6))>0.9$
if $t\ge 35$ and $c\ge 0.25$. Then,
$y'=-\big(y-(\G_c^h)^-(t)\big)^2+p(t)$ has an attractor-repeller pair
$(\wit a^-_{c,h},\wit r^-_{c,h})$, since the
trivial change of variables $x=y+1$ provides the equation
$x'=-\big(x-(1+(\G_c^h)^-(t))\big)^2+p(t)$, with $\|1+(\G_c^h)^-\|<0.1$.
As before, we have $-3.5< \wit r^-_{c,h}(t)<\wit a^-_{c,h}<3.5$ for all $t\in\R$.
In addition, by reviewing the proof of Theorem \ref{3.teorexit}, we observe that
$\ma_{c,h}(t)=\wit a^-_{c,h}(t)$ whenever $t\le-35$. This proves our assertion
concerning $\ma_{c,h}$. To prove it for $\mr_{c,h}$, we work with $(\G_c^h)^+\!$ and
with the change of variables $x=y-1$, using now that $\|1-(\G_c^h)^+\|<0.1$.
\par
\begin{figure}[]
\caption{Phase planes for two values of $(c,h)$. In the left figure, an attractor-repeller
pair exists (Case A), while in the right one there are no bounded solutions (Case C). The
observed behavior is similar for any value of $(c,h)\in(0,50]\times[0,6]$, being the collision
time always less than $20$.}
\includegraphics[trim={1.9cm 6cm 1.45cm 3.4cm},clip,width=\textwidth]{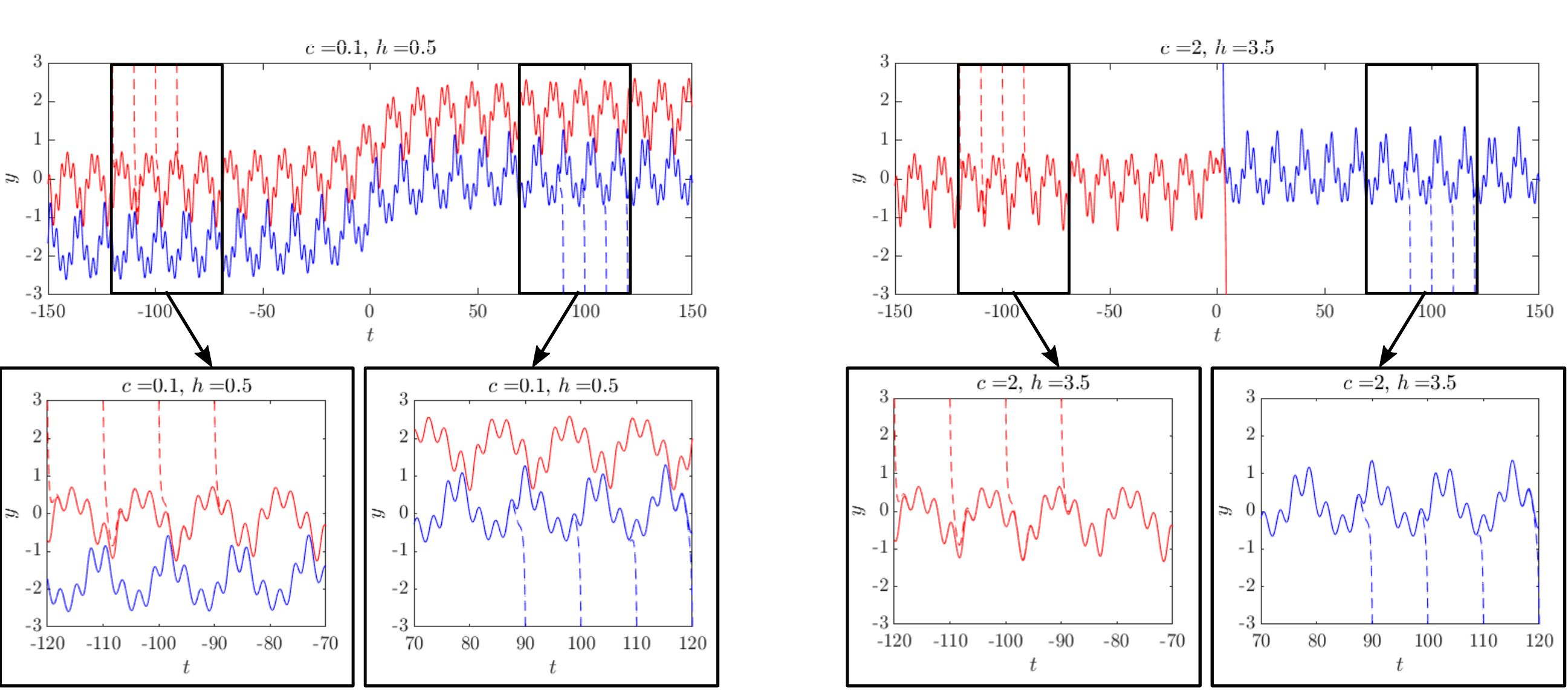}
\label{fig:zoom}
\vspace{-0.9cm}
\end{figure}

Our goal now is finding suitable pairs
(initial time, initial value) to reliably approximate $\ma_{c,h}$ and $\mr_{c,h}$
in the range of values $c\in(0,50]\cup\{\infty\}$ and $h\in[0,6]$
by finite integration. Recall that $\ma_{c,h}$ behaves like $\wit a^-_{c,h}$
on $(-\infty,-35]$, and hence it attracts exponentially fast solutions
starting above $\wit r^-_{c,h}$ as time increases: see Theorem \ref{2.teorhyp}.
Having in mind this fact, and trusting the simplicity of the numerical
integration that we are performing, we can say that the computer
does not distinguish $\ma_{c,h}(-35)$ from $y_{c,h}(-35,-500,3.5)$.
In fact, independently of the value of $(c,h)\in(0,50]\times\in[0,6]$, we observe that
the graph of any solution $y_{c,h}(t,s,3.5)$ with $s\le -85$ \lq\lq collides"~after
less that 20 units of time with the graph of $y_{c,h}(-35,-500,3.5)$: see Figure \ref{fig:zoom}.
And the same happens with $\mr_{c,h}(35)$ and $y_{c,h}(35,500,-3.5)$, so that the data
we are taking are very precautionary.
\par
The way to proceed is clear now. If we can continue the
solution $y_{c,h}(t,-500,3.5)$ at least until $t=35$, and observe that
$y_{c,h}(35,-500,3.5)>y_{c,h}(35,500,-3.5)$, then we
are in \hyperlink{CA}{\sc case A} (see Remark \ref{2.notabasta}).
If this is not the case, we will find $t_a>-35$ with
$y_{c,h}(t_a,-500,3.5)<-3.5$, which means that the graph of $y_{c,h}(t_a,-500,3.5)$
intersects that of any function taking values on $[-3.5,3.5]$, as is the case of any
possible bounded solution; therefore, there are no bounded solutions, and hence
the dynamics is given by \hyperlink{CC}{\sc case C}.
\par
\begin{figure}[]
\caption{Global dynamics of $x'=-x^2-0.2\,|x|-0.011+p(t)$. The behavior is analogous at any
interval of integration.
}
\includegraphics[trim={0.5cm 10.8cm 0.9cm 11cm},clip,width=\textwidth]{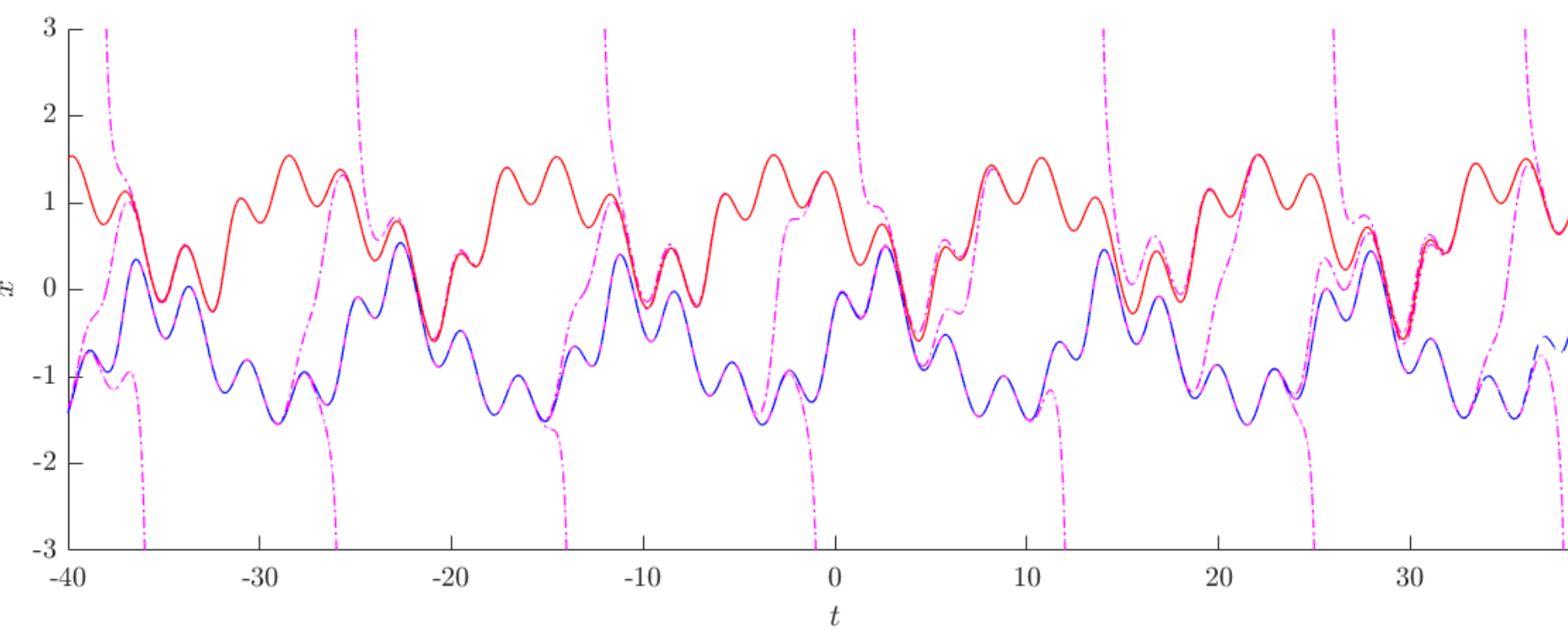}
\label{fig:ecunueva}
\vspace{-1cm}
\end{figure}
Let us justify our initial assumption. First, we check that the equation
$x'=-x^2-0.2\,|x|-0.011+p(t)$ has a bounded solution. In fact,
using the same MATLAB routine to
represent a large number of solutions of this equation, we observe that,
independently of the initial time, the numerical approximation of every
solution starting at an initial value greater than $3$ eventually falls
onto the graph of the function represented in solid red in
Figure \ref{fig:ecunueva}.
The analogous behavior is observed backwards in
time when computing solutions with initial value less than $-3$, which are
eventually mapped on the graph of the function represented in dashed
blue in Figure \ref{fig:ecunueva}. In addition, the solution corresponding to
any initial pair (initial time, initial value) between the graphs of both
functions falls onto the red curve as time increases and onto the blue curve as
time decreases. In other words, we observe numerically that the dynamics for the
(concave) equation $x'=-x^2-0.2\,|x|-0.011+p(t)$ is that of existence of
an attractor-repeller pair, which is more than required. Let $b$ be a
bounded solution, and take $\wit\G\in BPUC_\Delta(\R,\R)$ for a disperse set $\Delta$
with $\|\wit\G\|\le 0.1$. Then
$b'(t)=-b^2(t)-0.2\,|b(t)|-0.011+p(t)<-(b(t)-\wit\G(t))^2+p(t)$ for all $t\in\R\!-\!\Delta$,
and hence Theorem \ref{2.teoruno}(v) (see also Remark \ref{2.notaBV})
ensures that $x'=-(x-\wit\G(t))^2+p(t)$
has at least two different bounded solutions. According to Remark \ref{3.notaABC},
this ensures the existence of the attractor-repeller pair,
which was our initial assumption. This completes our explanation, and the appendix.

\end{document}